\documentclass[11pt]{article}
\usepackage{amssymb}
\usepackage{amsmath}
\usepackage{amsthm}
\newcommand{\cA}{{\cal A}}
\newcommand{\cB}{{\cal B}}
\newcommand{\cC}{{\cal C}}
\newcommand{\cD}{{\cal D}}
\newcommand{\cH}{{\cal H}}
\newcommand{\cE}{{\cal E}}

\newcommand{\cO}{{\cal O}}
\newcommand{\cL}{{\cal L}}

\newcommand{\cF}{{\cal F}}
\newcommand{\cK}{{\cal K}}
\newcommand{\cP}{{\cal P}}
\newcommand{\cQ}{{\cal Q}}
\newcommand{\cS}{{\cal S}}
\newcommand{\cT}{{\cal T}}
\newcommand{\cU}{{\cal U}}
\newcommand{\cV}{{\cal V}}

\newcommand{\cX}{{\cal X}}
\newcommand{\cY}{{\cal Y}}
\newcommand{\cZ}{{\cal Z}}

\renewcommand{\AA}{{\mathbb A}}

\newcommand{\GG}{{\mathbb G}}

\newcommand{\ZZ}{{\mathbb Z}}

\newcommand{\PP}{{\mathbb P}}

\newcommand{\SSS}{{\mathbb S}}

\newcommand{\gp}{\mathfrak{p}}
\newcommand{\gq}{\mathfrak{q}}

\newcommand{\gs}{\mathfrak{s}}

\newcommand{\gU}{\mathfrak{U}}
\newcommand{\gB}{\mathfrak{B}}

\newcommand{\on}{\operatorname}

\newcommand{\RCov}{\on{RCov}}

\newcommand{\Rep}{{\on{Rep}}}

\newcommand{\Qlb}{\mathbb{\bar Q}_\ell}
\newcommand{\Gm}{\mathbb{G}_m}

\newcommand{\A}{\mathbb{A}}

\newcommand{\toup}[1]{\stackrel{#1}{\to}}
\newcommand{\hook}[1]{\stackrel{#1}{\hookrightarrow}}

\newcommand{\getsup}[1]{\stackrel{#1}{\gets}}
\newcommand{\Sp}{\on{\mathbb{S}p}}
\newcommand{\Spin}{\on{\mathbb{S}pin}}
\newcommand{\GSpin}{\on{G\mathbb{S}pin}}
\newcommand{\GSp}{\on{G\mathbb{S}p}}
\newcommand{\IC}{\on{IC}}

\newcommand{\Hom}{\on{Hom}}

\newcommand{\Ext}{\on{Ext}}

\newcommand{\Sym}{\on{Sym}}
\newcommand{\SO}{\on{S\mathbb{O}}}

\newcommand{\Aut}{\on{Aut}}

\newcommand{\RG}{\on{R\Gamma}}

\newcommand{\sym}{\on{sym}}

\newcommand{\uPic}{\on{\underline{Pic}}}
\newcommand{\Bun}{\on{Bun}}

\newcommand{\Bunb}{\on{\overline{Bun}} }
\newcommand{\Bunt}{\on{\widetilde\Bun}}

\newcommand{\Spec}{\on{Spec}}

\newcommand{\Gr}{\on{Gr}}

\newcommand{\GL}{\on{GL}}

\newcommand{\PGL}{\on{PGL}}

\newcommand{\Eis}{{\on{Eis}}}
\newcommand{\pr}{\on{pr}}
\newcommand{\id}{\on{id}}
\newcommand{\QED}{$\square$} 
\newcommand{\Fq}{\mathbb{F}_q}  
\newcommand{\Fp}{\mathbb{F}_p}  % what for??????    

\newcommand{\iso}{{\widetilde\to}}

\newcommand{\comp}{\circ}
\newcommand{\Four}{\on{Four}}
\renewcommand{\H}{{\on{H}}}   %cohomologies
\newcommand{\R}{\on{R}\!}   %highest direct images
\newcommand{\DD}{\mathbb{D}}  %for duality
\newcommand{\D}{\on{D}}       %for derived categories     
\newcommand{\wt}{\widetilde}
\newcommand{\ov}[1]{\overline{#1}}
\newcommand{\select}[1]{{\it{#1}}}

\newcommand{\und}[1]{\underline{#1}}

\renewcommand{\div}{\on{div}}

\renewcommand{\P}{{\on{P}}}

\newcommand{\<}{\langle}
\renewcommand{\>}{\rangle}

\newcommand{\ev}{\mathit{ev}}

\newcommand{\Sph}{\on{Sph}}

\newcommand{\dimrel}{\on{dim.rel}}
\newcommand{\codim}{\on{codim}}

\newcommand{\SL}{\on{SL}}

\newcommand{\tboxtimes}{\,\tilde\boxtimes\,}

\newcommand{\Vect}{\on{Vect}}

\newcommand{\ra}{\rightarrow}
\newcommand{\la}{\leftarrow}

\newtheorem{Lm}{Lemma}
\newtheorem{Th}{Theorem}
\newtheorem{Pp}{Proposition}
\newtheorem{Cor}{Corollary}
\newtheorem{Con}{Conjecture}

\theoremstyle{remark}
\newtheorem{Rem}{Remark}

\theoremstyle{definition}
\newtheorem{Def}{Definition}

\newenvironment{Prf}{\par\noindent {\it Proof }}{\QED}

\begin{document}
\author{Sergey Lysenko}
\title{Geometric Whittaker models and Eisenstein series for $\wt\SL_2$}
\date{}
\maketitle
\begin{abstract}
\noindent{\scshape Abstract}\hskip 0.8 em 
Let $X$ be a smooth projective curve over an algebraically closed field of characteristic $>2$. Let $\Bunt_G$ be the stack of metaplectic bundles on $X$ of rank 2 from \cite{L1}. In this paper we study the derived category $\D_-(\Bunt_G)$ of genuine $\ell$-adic sheaves in the framework of the quantum geometric Langlands. We describe the corresponding Whittaker category, develop the theory of geometric Eisenstein series and calculate the most non-degenerate Fourier coefficients of these Eisenstein series. The existing constructions of automorphic sheaves for $\GL_n$ are based on using \select{Whittaker sheaves}. Our calculations lead to a conjectural characterization of the Whittaker sheaf for $\wt\SL_2$, though its existence is not clear. We also formulate a precise conjectural relation between the quantum Langlands functors and the theta-lifting functors for the dual pair $(\wt\SL_2, \PGL_2)$.
\end{abstract}

\bigskip

\centerline{\scshape 0. Introduction}

\bigskip\noindent
0.1 Let $X$ be a smooth projective curve over a finite field $k=\Fq$ with $q$ odd. Let $\Omega$ be the canonical line bundle on $X$ and $G=\SL(\cO_X\oplus\Omega)$, a group scheme over $X$. Let $F=k(X)$, $\AA$ the adeles of $F$, and $\cO\subset\AA$ the entire adeles. Let $\wt G(\AA)$ be the metaplectic extension of $G(\AA)$ by $\{\pm 1\}$. Write $\Bun_G$ for the moduli stack of $G$-torsors on $X$.

 The nonramified Langlands program for $\wt G(\AA)$ is the problem of analyzing the space of genuine functions in $L^2(G(F)\backslash \wt G(\AA)/G(\cO))$ as a module over the global Hecke algebra. In \cite{L1} we proposed a geometrization of this problem. We introduced a $\mu_2$-gerb $\Bunt_G$ over $\Bun_G$ and a category $\D_-(\Bunt_G)$, which is a geometric analog of the above space, it is acted on by the tensor category $\Rep(\SL_2)$ of Hecke functors at each point of $X$.

This paper arose as a trial to answer the following question. Can one construct automorphic sheaves in $\D_-(\Bunt_G)$ in a way similar to the case of $\GL_2$ via the Whittaker models? One could rather use the geometric theta-lifting functors from $\PGL_2$ to $\tilde G$, but the sheaves obtained by theta-lifting are expected to be irreducible perverse sheaves tensored by some vector spaces (see Conjecture~\ref{Con_2} in Section~4.4). Since it is not easy to get rid of these vector space factors, the above question looks natural.

We first describe the Whittaker category for $\wt G$ proving (\cite{G2}, Conjecture~0.4) in our particular case. Then we define and establish some basic properties of the geometric Eisenstein series for $\tilde G$. Our main result is a calculation of the most non-degenerate Whittaker coefficients of these Eisenstein series. This calculation is a geometric analog of a result of Maas \cite{M} of 1937. We show that these Whittaker coefficients are the central values of the standard L-function twisted by some quadratic characters. In particular, the first Whittaker coefficient of these Eisenstein series is not "one", but is given by the central value of the standard L-function.

 Based on this calculation and some other results, we propose a conjectural description of sheaves $\cS^d_E$ that should play a role of Whittaker sheaves for $\wt G$ for a given $\SL_2$-local system $E$. They should provide a "square root" of the central value of the standard L-function for $E$ twisted by quadratic characters. Let us also mention a question (posed to the author by Drinfeld) that could be important for the theory of motives: if $E$ is of motivic origin, is this square root also of motivic origin?

\medskip\noindent
0.2 {\scshape Notation} Let $k$ be an algebraically closed field of characteristic $p>2$. All the stacks we consider are defined over $k$. Let $X$ be a smooth projective connected curve. Write $\Omega$ for the canonical line bundle on $X$. Write $\Bun_m$ for the stack of rank $m$ vector bundles on $X$. 

 Fix a prime $\ell\ne p$. For a $k$-stack $S$ locally of finite type write $\D(S)$ for the category introduced in (\cite{LO}, Remark~3.21) and denoted $\D_c(S,\Qlb)$ in \select{loc.cit.} It should be thought of as the unbounded derived category of constructible $\Qlb$-sheaves on $S$. For $*=+,-,b$ we have the full subcategory $\D^*(S)\subset \D(S)$ denoted $\D_c^*(S,\Qlb)$ in \select{loc.cit.} 
Write $\D^{\ast}(S)_!\subset \D^{\ast}(S)$ for the full subcategory of objects which are extensions by zero from some open substack of finite type. Write $\D^{\prec}(S)\subset \D(S)$ for the full subcategory of complexes $K\in \D(S)$ such that for any open substack $U\subset S$ of finite type we have $K\mid_U\in \D^-(U)$. Let $\P(S)\subset \D(S)$ denote the category of perverse sheaves. Write $\DD: \D^b(S)\to \D^b(S)$ for the Verdier duality functor. 
 
 Fix a nontrivial character $\psi: \Fp\to\Qlb^*$ and denote by $\cL_{\psi}$ the corresponding Artin-Shreier sheaf on $\A^1$. Since we are working over an algebraically closed field, we systematically ignore the Tate twists. For a morphism of stacks $f: Y\to Z$ we denote by $\dimrel(f)$ the function of a connected component of $Y$ given by $\dim C-\dim C'$, where $C'$ is the connected component of $Z$ containing $f(C)$. 

 If $V\to S$ and $V^*\to S$ are dual rank $r$ vector bundles on $S$, we normalize the Fourier transform $\Four_{\psi}: \D^b(V)\to \D^b(V^*)$ by $\Four_{\psi}(K)=(p_{V*})_!(\xi^*\cL_{\psi}\otimes p_V^*K)[r]$, where $p_V,p_{V^*}$ are the projections, and $\xi: V\times V^*\to\A^1$ is the pairing. 

\medskip\noindent
0.3 {\scshape Main results}

\medskip\noindent
0.3.1 Let $n\ge 1$. Let $G=\Sp(M_n)$, where $M_n=\cO_X^n\oplus\Omega^n$ with a natural symplectic form $\wedge^2 M_n\to\Omega$, here $G$ is a group scheme on $X$. We write $G=G_n$ if we need to express the dependence on $n$. The stack $\Bun_G$ classifies $M\in\Bun_{2n}$ equipped with a symplectic form $\wedge^2 M\to\Omega$. Let $\cA$ be the line bundle on $\Bun_G$ whose fibre at $M\in\Bun_G$ is $\det\RG(X,M)$, it is $\ZZ/2\ZZ$-graded purely of degree zero.  Let $\Bunt_G$ be defined as in \cite{L1}. The stack $\Bunt_G$ classifies $M\in\Bun_G$, a $\ZZ/2\ZZ$-graded purely of degree zero 1-dimensional vector space $\cB$ equipped with a $\ZZ/2\ZZ$-graded isomorphism $\cB^2\,\iso\, \det\RG(X, M)$.

 Let $\epsilon$ be the 2-automorphism of $\Bunt_G$ acting on $(M,\cB)$ so that it acts trivially on $M$, and as $-1$ on $\cB$. By $\P_-(\Bunt_G)$ we denote the category of perverse sheaves on $\Bunt_G$ on which $\epsilon$ acts as $-1$, let $\D_-(\Bunt_G)$ denote the corresponding derived category. 
 
 Let $\bar\epsilon$ be the 2-automorphism of $\Bunt_G$ acting on $(M,\cB)$ so that it acts trivially on $\cB$ and as $-1$ on $M$. This makes sense, because $-1\in \Aut(M)$ acts trivially on $\det\RG(X,M)$. Write $\D_{=}(\Bunt_G)$ and $\D_{\mp}(\Bunt_G)$ for the full triangulated subcategory of $\D_-(\Bunt_G)$ of objects on which $\bar\epsilon$ acts as $-1$ and $1$ respectively.
The category $\D_-(\Bunt_G)$ decomposes as a direct product 
$$
\D_-(\Bunt_G)\,\iso\, \D_{=}(\Bunt_G)\times\D_{\mp}(\Bunt_G) 
$$
and similarly for $\P_-(\Bunt_G)$. This means that any $K\in \D_-(\Bunt_G)$ is of the form $K_-\oplus K_+$ with $K_-\in \D_{=}(\Bunt_G)$ and $K_+\in \D_{\mp}(\Bunt_G)$. Besides, $\Hom(K_-, K_+)=\Hom(K_+, K_-)=0$, where $\Hom$ is calculated in $\D_-(\Bunt_G)$.
 
  Write $\Rep(\Sp_{2n})$ for the category of finite-dimensional representations of $\Sp_{2n}$ over $\Qlb$. By \cite{L1}, $\Rep(\Sp_{2n})$ acts on $\D_-(\Bunt_G)$ by Hecke functors, see Sections~1.3-1.4. The category $\Rep(\Sp_{2n})$ is $\ZZ/2\ZZ$-graded by the action of the center of $\Sp_{2n}$. We check in Section~1.4.1 that the Hecke functors are compatible with these $\ZZ/2\ZZ$-gradings on $\Rep(\Sp_{2n})$ and on $\D_-(\Bunt_G)$. 
  
   Let $H=\SO_{2n+1}$ denote the split orthogonal group. The quantum Langlands conjecture from \cite{St} in our case specializes to the following (we do not precise the finiteness conditions on the corresponding derived categories).
   
\begin{Con} 
\label{Con_QL}
There is an equivalence $QL: \D(\Bun_H)\,\iso\, \D_-(\Bunt_G)$ commuting with the actions of $\Rep(\Sp_{2n})$ by Hecke functors. Under this equivalence the grading of $\D(\Bun_H)$ by the connected components of $\Bun_H$ corresponds to the grading of $\D_-(\Bunt_G)$ by the action of the 2-automorphism $\bar\epsilon$.
\end{Con}

  Let $P\subset G$ be the Siegel parabolic preserving $\Omega^n$. We fix an inclusion $\GL_n\hook{} P$ such that $g\in\GL_n$ acts as $g$ on $\Omega^n$ and as $^tg^{-1}$ on $\cO_X^n$. In this way we identify $\GL_n$ with the Levi quotient of $P$. For $n=1$ we also let $B=P$. 
  
  In Section~4.1 we define the geometric Eisenstein series functor $\Eis: \D(\Bun_n)\to \D_-(\Bunt_G)$ for the parabolic $P$. As in \cite{BG}, we check that 
$\Eis$ commutes with the Verdier duality (and maps pure complexes to pure ones). The quantum Langlands is expected to be compatible with the parabolic induction. For $P$ this is expressed as follows. 

 We need the geometric Eisenstein series functor $\Eis_{\GL_n}^H: \D(\Bun_n)\to\D(\Bun_H)$. Its definition from (\cite{BG}, Section~2.3.2) does not apply litterally, as in \cite{BG} it is assumed that the derived group is simply-connected. The modification to be done is as follows. Let $\bar H=\GSpin_{2n+1}$, let $L_{\bar H}$ be the preimage of $\GL_n\subset \SO_{2n+1}$ in $\bar H$, then $\Eis_{\GL_n}^H$ is characterised by the commutative diagram
$$
\begin{array}{ccc}
\D(\Bun_{L_{\bar H}}) & \toup{\Eis_{L_{\bar H}}^{\bar H}} &\D(\Bun_{\bar H})\\
\uparrow && \uparrow\\
\D(\Bun_n) & \toup{\Eis_{\GL_n}^H} & \D(\Bun_H),
\end{array}
$$
where the vertical arrows are the natural restriction functors. The definition of the Drinfeld compactifications (used to define the geometric Eisenstein series) for $G$ with $[G,G]$ not necessarily simply-connected will appear in \cite{SS}. 

\begin{Con} 
\label{Con_compatibility_parab_induction}
There is an automorphism $\delta: \Bun_n\,\iso\, \Bun_n$ such that the diagram is 2-commutative
$$
\begin{array}{ccc}
\D(\Bun_n) & \toup{\Eis} & \D_-(\Bunt_G)\\
\uparrow\lefteqn{\scriptstyle{\delta^*}} && \uparrow\lefteqn{\scriptstyle QL}\\
\D(\Bun_n) & \toup{\Eis^H_{\GL_n}} & \D(\Bun_H)
\end{array}
$$
If $n=1$ then $\delta(\cB)=\cB\otimes\Omega^{1/2}$ for some square root $\Omega^{1/2}$ of $\Omega$.
\end{Con}

 Since $\mu_2\subset G$ is central, $\Bun_{\mu_2}$ acts naturally on $\Bun_G$. In Section~5 we lift this action to an action on $\Bunt_G$, it preserves each of the categories $\D_{=}(\Bunt_G)$ and $\D_{\mp}(\Bunt_G)$. We also describe an action of $\Bun_{\mu_2}$ on $\D(\Bun_H)$ that should correspond to the latter one via the equivalence of Conjecture~\ref{Con_QL}. 
   
\medskip\noindent
0.3.2 For $m$ odd in (\cite{L1}, Section~7) we defined a morphism $\tilde\tau: \Bunt_{G_n}\times\Bun_{\SO_m}\to \Bunt_{G_{nm}}$. Let $\Aut$ denote the theta-sheaf on $\Bun_{G_n}$ introduced in (\cite{L1}, Definition~1). We have the diagram of projections
$$
\Bun_H\getsup{\gq}\Bunt_G\times\Bun_H\toup{\gp} \Bunt_G
$$
Define the theta-lifting functors $F_G: \D^-(\Bun_H)_!\to \D^{\prec}(\Bunt_G)$ and $F_H:\D^-(\Bunt_G)_!\to \D^{\prec}(\Bun_H)$ by
$$
F_G(K)=\gp_!(\tilde\tau^*\Aut\otimes \gq^*K)[-\dim\Bun_H+\dimrel(\tilde\tau)]
$$
$$
F_H(K)=\gq_!(\tilde\tau^*\Aut\otimes \gp^*K)[-\dim\Bun_G+\dimrel(\tilde\tau)]
$$

\noindent
0.3.3 Let $n=1$. In Section~1 we study the Whittaker category for $\wt\SL_2$ in a way similar to \cite{FGV}. We prove Lurie's Conjecture~0.4 from \cite{G2} in the particular case of $\wt\SL_2$ (corresponding to $G=\SL_2$ and $q$  a primitive 4-th root of unity in the notation of \cite{G2}). Namely, our Corollary~\ref{Cor_free_Whit_module} identifies the corresponding Whittaker category with $\Rep(\SL_2)$. We also prove that similarly to the case of a reductive group considered in \cite{FGV}
the objects of our Whittaker category are "clean" perverse sheaves (cf. Theorem~\ref{Th_2}). 

 Our description of the Whittaker category agrees with the uniqueness of the Whittaker functional for the representations of the metaplectic group obtained in \cite{S}. We did not check the compatibility of our result with the Casselman-Shalika formula for $\wt\SL_2$ from \cite{BFH}.
 
 In Sections~2,3 we define the analog of Whittaker functors for $\wt\SL_2$, similar to the Whittaker functors introduced by Gaitsgory in \cite{G} for $\GL_n$. 
Let  $\cQ_1$ be the stack classifying $(M,\cB)\in\Bunt_G$, $\cE\in\Bun_1$ with an inclusion of coherent sheaves $s: \cE\otimes\Omega\hook{} M$. Let $\cQ_{2,ex}$ be the stack classifying a point of $\cQ_1$ as above together with a section $s_2: \cE^2\to\cO$. We define a full triangulated subcategory $\D^{\cL}(\cQ_{2,ex})\subset \D(\cQ_{2,ex})$ singled out by some equivariance condition. 
We think of $\D^{\cL}(\cQ_{2,ex})$ as a "family" of Whittaker categories from Section~1 indexed by pairs $(\cE, s_2)$. The purpose of Section~3 is to prove Proposition~\ref{Pp_one}, which provides an equivalence 
$W_{ex}: \D(\cQ_1)\,\iso\, \D^{\cL}(\cQ_{2,ex})$ exact for the perverse t-structures. 
 
  The stack $\Bun_B$ classifies $\cE\in\Bun_1$ and an exact sequence on $X$
$$
0\to\cE\otimes\Omega\to M\to \cE^{-1}\to 0
$$ 
Let $\tilde\nu_B: \Bun_B\to \Bunt_G$ be the map sending this sequence to $(M,\cB)$, where $\cB=\det\RG(X, \cE\otimes\Omega)$ with the induced isomorphism $\cB^2\,\iso\, \det\RG(X,M)$. Let $\cS_B$ be the stack classifying $\cE\in\Bun_1$ and $s_2: \cE^2\to\cO_X$. Let 
\begin{equation}
\label{Four_from_Bun_B}  
\Four_{\psi}: \D(\Bun_B)\,\iso\, \D(\cS_B)
\end{equation}
denote the Fourier transform. Over the open substack of $\cQ_{2,ex}$ given by the property that $s:\cE\otimes\Omega\to M$ is a subbundle, the functor $W_{ex}$ is essentially given by (\ref{Four_from_Bun_B}). 

For $d\ge 0$ let $X^{(d)}$ denote the $d$-th symmetric power of $X$, write $^{rss}X^{(d)}\subset X^{(d)}$ for the open subscheme of reduced divisors. Write $\RCov^d$ for the stack classifying $\cE\in\Bun_1$, $D\in {^{rss}X^{(d)}}$ and $s_2: \cE^2\,\iso\, \cO_X(-D)$. So, $\RCov^d\subset \cS_B$ is an open substack. By abuse of notation, we write also $\bar\epsilon$ for the 2-automorphism of $\RCov^d$ acting as $-1$ on $\cE$. 
 
 For a rank one local system $V$ on $X$ write $AV$ for the automorphic local system on $\Bun_1$ corresponding to $V$ (cf. Section~4.3.1). 
 
  Let $E$ be an $\SL_2$-local system on $X$ assumed to be either irreducible or a direct sum of rank one local systems. Let $\Aut_E\in\D(\Bun_H)$ denote the Hecke eigensheaf corresponding to $E$ and normalized by requiring that its first Whittaker coefficient is "one". More precisely, for $E$ irreducible, this is the perverse sheaf denoted $\Aut_E$ in \cite{FGV}. For $E=V\oplus V^*$, where $V$ is a rank one local system on $X$, this is the geometric Eisenstein series $(AV)_{\Omega}\otimes \Eis_{\GL_1}^H(AV[\dim\Bun_1])$. Our purpose is to characterize the sheaf 
\begin{equation}
\label{sheaf_Whitt_we_seek_for}
\Four_{\psi}(\tilde\nu_B^* QL(\Aut_E))\mid_{\RCov^d}[\dimrel(\tilde\nu_B)-\dim\RCov^d]
\end{equation} 
on $\RCov^d$ for all $d\ge 0$, we think of it as a \select{Whittaker sheaf} for $\wt\SL_2$ by analogy with the case of $\GL_2$. Note that (\ref{sheaf_Whitt_we_seek_for}) is $\ZZ/2\ZZ$-graded by the action of $\bar\epsilon$. 
  
\medskip\noindent
0.3.4  The kernel of the Eisenstein series functor $\Eis: \D(\Bun_n)\to \D_-(\Bunt_G)$ is given by the perverse sheaf $\IC_{P,W}$ on $\Bunt_{\tilde P}$ introduced in Section~4.1.1. Assume $n=1$. We also write $\Bunb_{\tilde B}=\Bunt_{\tilde P}$. In Section~4 we prove Theorem~\ref{Th_3} that calculates the $*$-restrictions of $\IC_{B,W}$ to the natural stratification of $\Bunb_{\tilde B}$. 
This calculation adopts the method of \cite{BFGM} to the case of $\wt\SL_2$. We first calculate the local analog $\IC(W_{max})$ of this perverse sheaf over the Zastava spaces $Z^{\theta}$ in Proposition~\ref{Pp_restriction_zero_section} and derive Theorem~\ref{Th_3} as a consequence.

  In our particular case the Zastava space $Z^{\theta}$ is a vector bundle over $X^{(\theta)}$. We also study the Fourier transform $\Four_{\psi}(\IC(W_{max})$ of $\IC(W_{max})$ over the dual vector bundle denoted $\check{Z}^{\theta}$. A new feature compared to \cite{BFGM}, which probably does not happen for more general groups, is that $\Four_{\psi}(\IC(W_{max}))$ is described by an "explicit formula" given in Proposititon~\ref{Pp_zeta_2_summand_great}. Namely, in Section~4.2.3 we introduce a vector bundle $\check{Y}^{\theta}\to X^{(\theta)}$ and a morphism $\pi_{\check{Y}}:\check{Y}^{\theta}\to \check{Z}^{\theta}$ over $X^{(\theta)}$. We also define a group scheme $\GG^{\theta}_2$ on $X^{(\theta)}$ and a homomorphism $\zeta_2: \GG^{\theta}_2\to \mu_2$ in Section~4.2.5. We show that $\pi_{\check{Y}!}\IC(\check{Y}^{\theta})$ is a perverse sheaf acted on by $\GG^{\theta}_2$, and its $(\GG^{\theta}_2, \zeta_2)$-isotypic direct summand is identified with 
$\Four_{\psi}(\IC(W_{max}))$.

 Let now $V$ be a rank one local system on $X$ and $E=V\oplus V^*$. By Conjecture~\ref{Con_compatibility_parab_induction}, $QL(\Aut_E)$ should be identified with $(AV)_{\Omega^{3/2}}\otimes\Eis(AV[\dim\Bun_1])$. In Section~4.3 we prove Theorem~\ref{Th_4}, which is our main result. It calculates the complex
\begin{equation}
\label{complex_K_E_after_transfrom_on_RCov}
\cS_{V\oplus V^*}=
(AV)_{\Omega^{3/2}}\otimes\Four_{\psi}\tilde\nu_B^*\Eis(AV)[\dimrel(\tilde\nu_B)+\dim\Bun_1-\dim\RCov^d]\mid_{\RCov^d}
\end{equation}

  To spell the answer, we need more notation. For a point $(\cE, s_2)$ of $\RCov^d$ view $\cO_X\oplus\cE$ as an $\cO_X$-algebra with the product given by $s_2$ and let $Y=\Spec(\cO_X\oplus \cE)$. 
Then $Y$ is a smooth projective curve, the projection $\phi: Y\to X$ is a degree 2 covering ramified exactly at $D$. Let $\cE_0$ denote the $\mu_2$-antiinvariants in $\phi_!\Qlb$. This way $\RCov^d$ can be seen as the stack classifying $D\in {^{rss}X^{(d)}}$ and a degree two smooth covering $\phi: Y\to X$ ramified exactly over $D$. 
 
  Theorem~\ref{Th_4} claims that the fibre of (\ref{complex_K_E_after_transfrom_on_RCov}) at $(\cE, s_2)\in\RCov^d$ is the central value of the L-function
$$
\oplus_{\theta\ge 0} \RG(X^{(\theta)}, (V^*\otimes \cE_0)^{(\theta)})[\theta]
$$
tensored by some 1-dimensional space $A_{\cE}$ such that $A_{\cE}^2\,\iso\, \Qlb$. If $V^*\otimes\cE_0$ is nontrivial then the latter is the exteriour algebra of $\H^1(X, V^*\otimes \cE_0)$. In particular, 
(\ref{complex_K_E_after_transfrom_on_RCov}) is a local system on $\RCov^d$ for $d>0$. The occurrence of quadratic L-functions in the Whittaker coefficients of Eisenstein series for $\wt\SL_2$ was discovered in 1937 by Maas \cite{M}, our Theorem~\ref{Th_4} is a geometric analog of his result. 

 Note also that there is a large recent letterature on rapidly developping theory of Weyl group multiple Dirichlet series (see \cite{BBF, CG} or \cite{B} for a survey). Hopefully, our Theorem~\ref{Th_4} is an example of a geometrization of such a Dirichlet series.
  
  We also discuss the conjectural functional equation for $\Eis$ in Remark~\ref{Rem_functional_equation}. According to the functional equations for $\Eis$ and for $\Eis_{\GL_1}^H$, for $n=1$ we expect that $\delta:\Bun_1\to\Bun_1$ in Conjecture~\ref{Con_compatibility_parab_induction} sends $L$ to $L\otimes\Omega^{1/2}$. 
  
  Using Theorem~\ref{Th_4} we determine some residues of $\Eis(AV[\dim\Bun_1])$ in Section~4.3.3. Namely, given $(\cE, s_2)\in \RCov^0$ let $\cE_0$ be the corresponding rank one and order two local system on $X$. Then $\Eis(A\cE_0[\dim\Bun_1])$ is unbounded, and it admits $\tilde\sigma_{\cE}^*\Aut$ as a residue. Here $\tilde\sigma_{\cE}$ is the automorphism of $\Bunt_G$ defined in Section~5.1.

 In Section~4.3.4 we calculate the constant terms of $\Eis(K)$ for $K\in \D(\Bun_1)$, the calculation is based on Theorem~\ref{Th_3}.
 
\medskip 
\noindent
0.3.5  For any local system $V$ on $X$ let $CL_V^d$ denote complex on $\RCov^d$ whose fibre at $(\cE, s_2)$ is
$$
\oplus_{\theta\ge 0} \RG(X^{(\theta)}, (V\otimes\cE_0)^{(\theta)})[\theta]
$$
Here $\cE_0$ is the $\mu_2$-antiinvariants in $\phi_!\Qlb$ for the covering $\phi:Y\to X$ given by $(\cE, s_2)$. The notation $CL$ stands for \select{central value of the L-function}. For example, if $V$ is irreducible and $d>0$ then $CL_V^d$ is a local system. 

  Here is a conjectural characterization of (\ref{sheaf_Whitt_we_seek_for}) for any $\SL_2$-local system $E$ on $X$. It is supported by 
Theorem~\ref{Th_4} on one hand, and by the uniqueness of the Whittaker models for $\wt\SL_2$ on the other hand. It is also consistent with the result of Waldspurger, who showed in \cite{W} that the Whittaker coefficients of the cuspidal automorphic functions on $\wt\SL_2$ are essentially square roots of quadratic twists of L-functions attached to cuspidal forms on $\GL_2$. Another calculation of the square of this Whittaker coefficient in the classical setting is given in \cite{P}. 
  
\begin{Con} 
\label{Con_existence_cS}
For any $\SL_2$-local system $E$ on $X$ and any $d\ge 0$ there is a complex $\cS_E^d\in \D(\RCov^d)$ equipped with an isomorphism 
$$
(\cS_E^d)^{\otimes 2}\,\iso\, CL^d_E
$$
Besides, $\cS_E^d[\dim\RCov^d]$ is Verdier self-dual. The complex (\ref{sheaf_Whitt_we_seek_for}) identifies canonically with $\cS_E^d$. 
\end{Con}

 For $E=V\oplus V^*$ Conjecture~\ref{Con_existence_cS} is true and follows from Theorem~\ref{Th_4}.

\begin{Rem} 
Of course, we also expect that $\cS^d_E$ well depends on $E$. In the D-modules setting assuming $k$ algebraically closed of characteristic zero, 
one has the stack $LS_{\SL_2}$ of $\SL_2$-local systems on $X$, and $\cS^d_E$ should naturally form a complex on $\RCov^d\times LS_{\SL_2}$ as $E$ varies. 
\end{Rem}

\medskip

  In Section~4.4 we also propose Conjecture~\ref{Con_2} giving a relation between the functor $QL$ and the theta-lifting functors $F_G, F_H$. To the difference from $QL$, neither $F_G$ nor $F_H$ is expected to be an equivalence. 
  
 We conjecture also that if $E$ is a rank one local system on $X$ with $E^2$ nontrivial then the contribution of each connected component of $\Bun_1$ to $\Eis(AE[\dim\Bun_1])$ is an irreducible perverse sheaf. We prove a part of this claim in Appendix~A. Here we observe the following important phenomenon. In the situation of the quantum geometric Langlands the Hecke eigensheaves, which are irreducible perverse sheaves, can not be normalized by requiring that their first Whittaker coefficient is "one". Finally, we discuss the conjectural functor $QL$ in the case of genus zero in Section~A.2. 
  
\bigskip

\centerline{\scshape 1. Whittaker category for $\wt\SL_2$}

\bigskip\noindent
1.1 In Section~1 we assume $n=1$ everywhere (except Sections~1.3 and 1.4.1, where $n\ge 1$ can be arbitrary).

 Fix an effective divisor $D\ge 0$ and a line bundle $\cE$ on $X$ equipped with $\cE^2\,\iso\, \cO(-D)$. Let $B\subset G$ be the parabolic subgroup scheme over $X$ preserving $\Omega$. Let $N\subset B$ be the unipotent radical, so $N\,\iso\, \Omega$.  
 
 Let $R$ denote the group scheme $\Omega(-D)$ on $X$. Then $\Bun_R$ is the stack classifying exact sequences $0\to \cE\otimes\Omega\to M\to \cE^{-1}\to 0$ on $X$. This is a twisted version of $\Bun_N$. 
 
Let $\Bunb_R$ be the stack classifying $M\in\Bun_G$ and $s: \cE\otimes\Omega\hook{} M$ an inclusion of coherent sheaves. Let 
$$
\Bunb_{\tilde R}=\Bunb_R\times_{\Bun_G}\Bunt_G\;\;\;\;\;\mbox{and}\;\;\;\;\;\;
\Bun_{\tilde R}=\Bun_R\times_{\Bun_G}\Bunt_G
$$ 

 Fix $x\in X$. Let also $_{x,\infty}\Bunb_R$ be the stack classifying $M\in \Bun_G$ and a section $s: \cE\otimes\Omega\hook{} M(\infty x)$. This is an ind-algebraic stack. Let also 
$$
_{x,\infty}\Bunb_{\tilde R}={_{x,\infty}\Bunb_R}\times_{\Bun_G}\Bunt_G
$$

 The canonical inclusion $R\hook{} \Omega$ yields a morphism of extension of scalars $\Bun_R\to\Bun_{\Omega}$. Let $\ev_{\Omega}:\Bun_{\Omega}\to\A^1$ be the map sending an exact sequence $0\to\Omega\to ?\to \cO_X\to 0$ to the corresponding element of $\H^1(X,\Omega)=k$. Let $\cL=\ev_{\Omega}^*\cL_{\psi}$. We will define the category 
$\P^{\cL}(_{x,\infty}\Bunb_{\tilde R})$ and its derived version $\D^{\cL}(_{x,\infty}\Bunb_{\tilde R})$ following the pattern of (\cite{L2}, Section~2.2.3).
 As in Section~0.3.1, requiring that $\epsilon\in\mu_2$ acts as $-1$, we will further get the categories 
$$
\P_-(_{x,\infty}\Bunb_{\tilde R})\subset \D_-(_{x,\infty}\Bunb_{\tilde R})
$$ 
 
\medskip\noindent
1.2 \select{Definition of Whittaker category} \   We need the following version of the relative determinant. Given vector bundles $M,M'$ on $X$, an effective divisor $D_1$ on $X$ and an isomorphism
$$
\tau: M\mid_{X-D_1}\,\iso\, M'\mid_{X-D_1},
$$
one defines a $\ZZ/2\ZZ$-graded line $\det(M:M')$ as follows. For any lower modification $V\subset M\cap M'$ over $D_1$, set
$$
\det(M:M')=\det\H^0(X, M/V)\otimes\det\H^0(X, M'/V)^{-1}
$$ 
Then $\tau$ yields a $\ZZ/2\ZZ$-graded isomorphism
\begin{equation}
\label{iso_detRG_modification}
\det\RG(X, M')\,\iso\, \det\RG(X,M)\otimes \det(M':M)
\end{equation}

For each $d\ge 0$ let $_{(d)}\cY$ be the stack classifying $D_1\in X^{(d)}$, $(M, s)\in {_{x,\infty}\Bunb_{R}}$ such that $x\notin D_1$, and the inclusion $s: \cE\otimes\Omega\hook{} M(\infty x)$ is a subbundle in a neighbourhood of $D_1$. For a point of $_{(d)}\cY$ over the formal neighbourhood of $D_1$ one gets an exact sequence
\begin{equation}
\label{ext_cE^*_by_EotimesOmega}
0\to \cE\otimes\Omega\to M\to \cE^{-1}\to 0
\end{equation}
It can be seen as an $R$-torsor $\cF_R$ over a formal neighbourhood $X_{D_1}$ of $D_1$ in $X$. Set also $_{(d)}\tilde\cY={_{(d)}\cY}\times_{\Bun_G}\Bunt_G$. 
 
 Let $_{(d)}\cH\cY$ be the stack classifying: $D_1\in X^{(d)}$, $(M, s)\in {_{x,\infty}\Bunb_R}$, $(M', s')\in {_{x,\infty}\Bunb_R}$, an isomorphism $\tau: M\mid_{X-D_1}\,\iso\, M'\mid_{X-D_1}$ of $G$-torsors
such that the diagram commutes
\begin{equation}
\label{diag_for_d_HY}
\begin{array}{ccc}
\cE\otimes\Omega & \hook{s} & M(\infty x)\mid_{X-D_1} \\
 & \searrow\lefteqn{\scriptstyle s'} & \downarrow\lefteqn{\scriptstyle\tau}\\
 && M'(\infty x)\mid_{X-D_1},
\end{array}
\end{equation}
it is required that $(D_1, M,s), (D_1, M', s')\in {_{(d)}\cY}$. The diagram (\ref{diag_for_d_HY}) yields a $\ZZ/2\ZZ$-graded isomorphism $\det(M: M')\,\iso\, k$, and in turn (\ref{iso_detRG_modification}) yields an isomorphism
\begin{equation}
\label{iso_detRG_M_and_M'}
\det\RG(X, M)\,\iso\, \det\RG(X,M')
\end{equation}
for a point of $_{(d)}\cH\cY$. 

  Let $_{(d)}\cH\tilde\cY$ be the stack over $_{(d)}\cH\cY$ classifying the same data together with one-dimensional vector space $\cB$ equipped with isomorphisms $\sigma: \cB^2\,\iso\, \det\RG(X,M)$ and $\sigma': \cB^2\,\iso\, \det\RG(X,M')$ compatible with (\ref{iso_detRG_M_and_M'}). 

 We get a diagram
$$
_{(d)}\tilde\cY
 \;\getsup{h^{\la}_{\tilde\cY}}\; {_{(d)}\cH\tilde\cY} \;\toup{h^{\ra}_{\tilde\cY}}\; {_{(d)}\tilde\cY},
$$
where $h^{\la}_{\tilde\cY}$ (resp., $h^{\ra}_{\tilde\cY}$) sends the above point to $(D_1, M,s, \cB)$ (resp., to $(D_1, M',s', \cB)$). In this way $_{(d)}\cH\tilde\cY$ is naturally a groupoid over $_{(d)}\tilde\cY$.
 
  Let $\ev_{\cH\cY}: {_{(d)}\cH\cY}\to\Bun_{\Omega}$ be the following map. Write $\cF_R$ and $\cF'_R$ for the $R$-torsors over $X_{D_1}$ associated to $(D_1, M,s)$ and $(D_1, M', s')$ as above. The isomorphism $\tau$ yields an isomorphism $\cF_R\,\iso\, \cF'_R$ over the punched formal neighbourhood $X^0_{D_1}=X_{D_1}-D_1$. Let $\cF^1_R$ be the $R$-torsor over $X_{D_1}$ such that the sum of $\cF_R$ with $\cF^1_R$ is canonically identified with $\cF'_R$, so $\cF^1_R$ is trivialized over $X^0_{D_1}$. We also denote by $\cF^1_R$ the $R$-torsor on $X$ obtained by gluing $\cF^1_R$ with the trivial $R$-torsor $\cF^0_R$ over $X-D_1$ via the above trivialization. Finally, $\ev_{\cH\cY}$ sends the above point of $_{(d)}\cH\cY$ to the extension of scalars of $\cF^1_R$ via $R\hook{}\Omega$. 
  
  Denote by $\ev_{\cH\tilde\cY}$ the composition ${_{(d)}\cH\tilde\cY}\to {_{(d)}\cH\cY}\toup{\ev_{\cH\cY}}\Bun_{\Omega}$. Set $h^{\la}=\pr\comp h^{\la}_{\tilde\cY}$, $h^{\ra}=\pr\comp h^{\ra}_{\tilde\cY}$, where $\pr: {_{(d)}\tilde\cY}\to {_{x,\infty}\Bunb_{\tilde R}}$ is the projection forgetting $D_1$. We will use the fact  that $\ev_{\cH\tilde\cY}$ is a homomorphism of groupoids, that is, it sends the composition of arrows to the sum of the corresponding $\Omega$-torsors in $\Bun_{\Omega}$. 

 Given $d\ge 1$, as in \cite{GN}, a \select{generic Hecke correspondence over $X^{(d)}$} is a stack locally of finite type $Y$ equipped with a morphism $\alpha: Y\to {_{(d)}\cH\tilde\cY}$ such that in the diagram
$$
\begin{array}{ccccc}
_{x,\infty}\Bunb_{\tilde R} & \getsup{h^{\la}_Y} & Y & \toup{h^{\ra}_Y} & {_{x,\infty}\Bunb_{\tilde R}}\\
\Arrowvert && \downarrow\lefteqn{\scriptstyle\alpha} && \Arrowvert\\
_{x,\infty}\Bunb_{\tilde R} & \getsup{h^{\la}} & {_{(d)}\cH\tilde\cY} & \toup{h^{\ra}} & {_{x,\infty}\Bunb_{\tilde R}}
\end{array}
$$
the maps $h^{\la}_Y$, $h^{\ra}_Y$ are smooth. For such correspondence we let $\tilde h^{\ra}_Y=h^{\ra}_{\tilde\cY}\comp\alpha$ and $\tilde h^{\la}_Y=h^{\la}_{\tilde\cY}\comp\alpha$.

 The generic Hecke correspondence is called \select{trivial} if it is equipped with a decomposition of $\alpha$ as $Y\to {_{(d)}\tilde\cY}\to {_{(d)}\cH\tilde\cY}$, where the second map is the diagonal. A morphism of generic Hecke correspondences over $X^{(d)}$ is a morphism $\beta: Y_1\to Y_2$ such that the diagram commutes
$$
\begin{array}{ccc}
Y_1 & \toup{\beta} & Y_2\\
 & \searrow & \downarrow\\
 &&  _{(d)}\cH\tilde\cY
\end{array}
$$ 
 
  A perverse sheaf $K$ on $_{x,\infty}\Bunb_{\tilde R}$ is called \select{generic Hecke $\cL$-equivariant for $d$} if for each generic Hecke correspondence $\alpha:  Y\to {_{(d)}\cH\tilde\cY}$ over $X^{(d)}$, $K$ is equipped with an isomorphism
$$
I_Y: h^{\ra *}_Y K\,\iso\, (h^{\la *}_Y K)\otimes \alpha^*\ev_{\cH\tilde\cY}\cL
$$
satisfying the following conditions. First, for any morphism $\beta: Y_1\to Y_2$ of generic Hecke correspondences over $X^{(d)}$, it is required that 
$$
I_{Y_1}=\beta^*(I_{Y_2})
$$
Second, consider a stack locally of finite type $Y$, and for $i\in \ZZ/3\ZZ$, generic Hecke correspondences $Y_i$, and morphisms $\beta_i: Y\to Y_i$ such that the compositions $h^{\la}_{Y_i}\comp \beta_i$, $h^{\ra}_{Y_i}\comp \beta_i$ are smooth, and the diagram commutes
$$
\begin{array}{ccc}
Y & \toup{\beta_{i+1}} & Y_{i+1}\\
\downarrow\lefteqn{\scriptstyle \beta_i} && \downarrow\lefteqn{\scriptstyle \tilde h^{\ra}_{Y_{i+1}}}\\
Y_i & \toup{\tilde h^{\la}_{Y_i}} & {_{(d)}\tilde\cY}
\end{array}
$$
Then it is required that the isomorphism
$$
\beta_3^*(I_{Y_3})\comp \beta_2^*(I_{Y_2})\comp \beta_1^*(I_{Y_1})
$$
be the identity. For the trivial generic Hecke correspondence it is required that $I_Y$ be the identity morphism. 

Say that $K$ on $_{x,\infty}\Bunb_{\tilde R}$ is \select{generic Hecke $\cL$-equivariant} if it is generic Hecke $\cL$-equivariant for each $d\ge 1$, and satisfies the factorization property. That is, for a divisor $D_1$, which is a disjoint union of $D_2$ and $D_3$, the isomorphism $I_Y$ is the product of the corresponding isomorphisms for $D_2$ and $D_3$. This makes sense as the groupoids $_{(d)}\cH\tilde\cY$ satisfy this factorization property. The generic Hecke $\cL$-equivariant structure on $K$ is determined by the isomorphisms $I_Y$ for $d=1$. 

 The category $\P^{\cL}(_{x,\infty}\Bunb_{\tilde R})$ is defined as the category of generic Hecke $\cL$-equivariant sheaves. Let $\P^{\cL}_-(_{x,\infty}\Bunb_{\tilde R})$ be its full subcategory of perverse sheaves on which $\epsilon$ acts as $-1$. Since the group scheme $R$ is unipotent, the natural functor $\P^{\cL}(_{x,\infty}\Bunb_{\tilde R})\to\P(_{x,\infty}\Bunb_{\tilde R})$ is fully faithful.
 
\medskip\noindent
1.3  \select{Genuine spherical sheaves} \  In Sections~1.3 and 1.4.1 we assume $n\ge 1$. Write $\Gr_{G,x}$ for the ind-scheme classifying $M\in\Bun_G$ together with a trivialization $M\,\iso\, M_0\mid_{X-x}$. Let $\cL$ be the line bundle on $\Gr_{G,x}$ with fibre $\det(M_{0,x}: M_x)$ as in \cite{L1}. Let $\wt\Gr_{G,x}$ denote the $\mu_2$-gerb of square roots of $\cL$. 
In (\cite{L1}, Section~8.1) we considered the category $\Sph(\wt\Gr_{G,x})$ of $G(\cO_x)$-equivariant perverse sheaves on $\wt\Gr_{G,x}$ on which $\epsilon$ acts as $-1$ called \select{genuine spherical sheaves}. 

 Remind the tensor category $\Sph(\wt\Gr_{G,x})^{\flat}$ introduced in (\cite{L1}, Definition~6). Remind that first one equips $\Sph(\wt\Gr_{G,x})$ with a convolution product, associativity constraint and a commutativity constraint (via fusion product) as in (\cite{L1}, Section~8.3). Then one introduces a new $\ZZ/2\ZZ$-grading on $\Sph(\wt\Gr_{G,x})$, namely, an irreducible object $K\in \Sph(\wt\Gr_{G,x})$ is even (resp., odd) if it appears as a direct summand in an even (resp., odd) tensor power of $\cA_{\alpha}$. Here $\alpha$ is the highest weight of the standard representation of the Langlands dual group $\check{G}=\SO_{2n+1}$. For each dominant coweight $\lambda$ of $G$, $\cA_{\lambda}\in \Sph(\wt\Gr_{G,x})$ is defined on (p. 438, Section~8.1). 

 One then defines $\Sph(\wt\Gr_{G,x})^{\flat}$ as the category of even objects in 
$\Sph(\wt\Gr_{G,x})\otimes \Vect^{\epsilon}$, where $\Sph(\wt\Gr_{G,x})$ is considered with this new $\ZZ/2\ZZ$-grading. By (\cite{L1}, Theorem~3), one has a canonical equivalence of tensor categories 
$$
\Sph(\wt\Gr_{G,x})^{\flat}\,\iso\, \Rep(\Sp_{2n})
$$

\medskip\noindent
1.4.1 \select{Hecke functors} \  Let $_x\cH_G$ be the Hecke stack classifying $G$-torsors $\cF_G,\cF'_G$ and an isomorphism $\beta:\cF_G\,\iso\, \cF'_G\mid_{X-x}$. We have the diagram
$$
\Bun_G\;\getsup{p}\; {_x\cH_G}\;\toup{p'}\;\Bun_G,
$$
where $p$ (resp., $p'$) sends the above point to $\cF_G$ (resp., $\cF'_G$). Let $_x\tilde\cH_G$ be the stack obtained from $\Bunt_G\times\Bunt_G$ by the base change $_x\cH_G\toup{p,p'}\Bun_G\times\Bun_G$. Denote by $\tilde p, \tilde p'$ the corresponding projections
$$
\Bunt_G\;\getsup{\tilde p}\; {_x\tilde\cH_G}\; \toup{\tilde p'}\;\Bunt_G
$$
The stack $_x\tilde\cH_G$ classifies $(M,\cB), (M',\cB')\in\Bunt_G$ together with an isomorphism of $G$-torsors $\beta: M\,\iso\, M'\mid_{X-x}$. 

 Denote by $\Bun_{G,x}$ the stack classifying $\cF_G$ on $X$ together with a trivialization $\nu:\cF_G\,\iso\, \cF^0_G\mid_{D_x}$, here $D_x$ is the formal neighbourhood of $x$. Set $\Bunt_{G,x}=\Bunt_G\times_{\Bun_G}\Bun_{G,x}$. 
 
 Denote by $\gamma$ (resp., by $\gamma'$) the isomorphism $\Bun_{G,x}\times_{G(\cO_x)}\Gr_{G,x}\,\iso\, {_x\cH_G}$ such that the projection to the first term corresponds to $p$ (resp., to $p'$). One has a canonical isomorphism
\begin{equation}
\label{iso_cA_over_xcH_G}
\gamma'^*p^*\cA \,\iso\, \cA\tboxtimes \cL^{-1}
\end{equation}
This yields a $G(\cO_x)$-torsor $\Bunt_{G,x}\times \wt\Gr_{G,x}\to {_x\tilde\cH_G}$ extending the $G(\cO_x)$-torsor
$$
\Bun_{G,x}\times\Gr_{G,x}\to \Bun_{G,x}\times_{G(\cO_x)}\Gr_{G,x}\toup{\gamma'}{_x\cH_G}
$$
We get the cartesian square
\begin{equation}
\label{diag_for_H_operators_1}
\begin{array}{ccc}
\Bunt_{G,x}\times \wt\Gr_{G,x} & \toup{\pr} &\Bunt_{G,x}\\
\downarrow && \downarrow\\
_x\tilde\cH_G & \toup{\tilde p'} & \Bunt_G
\end{array}
\end{equation}
So, for $\cS\in \Sph(\wt\Gr_{G,x})$, $\cT\in \D(\Bunt_G)$ we can form their twisted exteriour product $(\cT\tboxtimes\cS)^r\in \D(_x\tilde\cH_G)$ for the above torsor. The Hecke functor $\H^{\ra}: \Sph(\wt\Gr_{G,x})\times\D(\Bunt_G)\to \D(\Bunt_G)$ is given by
$$
\H^{\ra}(\cS,\cT)=\tilde p_!((\cT\tboxtimes\cS)^r)
$$
This is a right action, that is, $\H^{\ra}(\cS_2, \H^{\ra}(\cS_1,\cT))\,\iso\, \H^{\ra}(\cS_1\ast\cS_2,\cT)$, where $\ast$ denotes the convolution.

\begin{Lm} 
\label{Lm_compatibility_pm_with_Hecke}
The functor $\H^{\ra}(\cA_{\alpha},\cdot)$ sends $\H^{\ra}(\cA_{\alpha},\cdot): \D_{=}(\Bunt_G)\to \D_{\mp}(\Bunt_G)$ and $\H^{\ra}(\cA_{\alpha},\cdot): \D_{\mp}(\Bunt_G)\to \D_{=}(\Bunt_G)$.
\end{Lm}
\begin{Prf}
Let $_x\cH^{\alpha}_G$ be the substack of $_x\cH_G$ given by the property that $\cF_G$ and $\cF'_G$ are in the position $\alpha$ at $x$. Let $_x\wt\cH^{\alpha}_G={_x\cH^{\alpha}_G}\times{_x\cH_G}\; {_x\wt\cH_G}$. In (\cite{L1}, Lemma~14) we constructed the isomorphisms
$$
\kappa, \kappa': {_x\wt\cH^{\alpha}_G}\,\iso\, (\Bunt_G\times_{\Bun_G}\,{_x\cH^{\alpha}_G})\times B(\mu_2),
$$
where we used $p: \cH^{\alpha}_H\to\Bun_G$ (resp., $p': \cH^{\alpha}_H\to\Bun_G$) in the fibred product, and the projection to the fist term corresponds to $\tilde p$ (resp., to $\tilde p'$). 
For $K\in \D(\Bunt_G)$ by definition, 
$$
\H^{\ra}(\cA_{\alpha}, K)=\tilde p_!(\tilde p^{\prime *}K\otimes \kappa^*W)[2n+1]
$$
Remind that $\kappa$ followed by the projection onto $B(\mu_2)$ sends a point of $_x\wt\cH^{\alpha}_G$ to $\cB_0$, where
$$
\cB'=\cB\otimes ((M+M')/M')^*\otimes \cB_0
$$ 
Let $\tau$ be the 2-automorphism of $_x\wt\cH^{\alpha}_G$ acting as $-1$ on $M,M'$, and as 1 on $\cB, \cB'$. The image of $\tau$ under both $\tilde p$ and $\tilde p'$ is $\bar\epsilon$. Assume that $\bar\epsilon$ acts on $K\in \D(\Bunt_G)$ as $a\in \{1,-1\}$. Since $\tau$ acts on $\kappa^*W$ as $-1$, we see that $\tau$ acts on $\tilde p^{\prime *}K\otimes \kappa^*W$ as $-a$. So, $\bar\epsilon$ acts on $\H^{\ra}(\cA_{\alpha}, K)$ as $-a$.
\end{Prf} 
 
\medskip\noindent 
1.4.2 Assume $n=1$ again. Let $Z={_x\tilde\cH_G}\times_{\Bunt_G}{_{x,\infty}\Bunb_{\tilde R}}$, where we used $\tilde p'$ to define the fibre product. Let $p'_Z: Z\to {_{x,\infty}\Bunb_{\tilde R}}$ be the second projection. The stack $Z$ classifies collections:
$(M,\cB), (M',\cB')\in\Bunt_G$, an isomorphism of $G$-torsors $\beta: M\,\iso\, M'\mid_{X-x}$, and $s': \cE\otimes\Omega\hook{} M'(\infty x)$. 

 Let $p_Z: Z\to {_{x,\infty}\Bunb_{\tilde R}}$ be the map sending the above point to $(M,\cB, s)$, where $s$ is the composition 
$$
\cE\otimes\Omega\hook{s'} M'(\infty x)\toup{\beta^{-1}} M(\infty x)
$$ 
We get the diagram, where both squares are cartesian
$$
\begin{array}{ccccc}
_{x,\infty}\Bunb_{\tilde R} & \getsup{p_Z} & Z & \toup{p'_Z} & _{x,\infty}\Bunb_{\tilde R}\\
\downarrow && \downarrow && \downarrow\\
\Bunt_G & \getsup{\tilde p} & {_x\tilde\cH_G} & \toup{\tilde p'} & \Bunt_G
\end{array}
$$
Let $_{x,\infty}\Bunb_{\tilde R,x}={_{x,\infty}\Bunb_{\tilde R,x}}\times_{\Bun_G}\Bun_{G,x}$. As for (\ref{diag_for_H_operators_1}), we get the cartesian square
$$
\begin{array}{ccc}
_{x,\infty}\Bunb_{\tilde R,x}\times \wt\Gr_{G,x} & \toup{\pr} & _{x,\infty}\Bunb_{\tilde R,x}\\
\downarrow && \downarrow\\
Z & \toup{p'_Z} & _{x,\infty}\Bunb_{\tilde R},
\end{array}
$$
where vertical arrows are $G(\cO_x)$-torsors. For $\cT\in \D(_{x,\infty}\Bunb_{\tilde R})$, $\cS\in \Sph(\wt\Gr_{G,x})$ we denote by $(\cT\tboxtimes\cS)^r\in \D(Z)$ the corresponding twisted exteriour product. We define 
$$
\H^{\ra}: \Sph(\wt\Gr_{G,x})\times \D(_{x,\infty}\Bunb_{\tilde R})\to \D(_{x,\infty}\Bunb_{\tilde R})
$$
by
$\H^{\ra}(\cS,\cT)=p_{Z !}((\cT\tboxtimes\cS)^r)$. 

 As in Section~1.2, we can define the category $\P^{\cL}(Z)$, and if $\cT\in \P^{\cL}(_{x,\infty}\Bunb_{\tilde R})$, $\cS\in \Sph(\wt\Gr_{G,x})$ then $(\cT\tboxtimes\cS)^r\in \P^{\cL}(Z)$ naturally. So, the complex $\H^{\ra}(\cS,\cT)$
inherits the generic Hecke $\cL$-equivariant structure, each perverse cohomology sheaf of $\H^{\ra}(\cS,\cT)$ lies in $\P^{\cL}(_{x,\infty}\Bunb_{\tilde R})$.
 
\medskip\noindent
1.5 \select{Objects of the Whittaker category} \  Let $\Lambda$ denote the coweight lattice of $G$, let $\check{\Lambda}$ be the weight lattice of $G$. Remind that $\alpha$ denotes the unique coroot, let $\check{\alpha}$ be the corresponding root. As in \cite{L1}, let $\check{\omega}_1$ denote the highest weight of the standard representation of $G$. Let $\Lambda^+$ denote the set of dominant coweights. 

 For $\lambda\in\Lambda$ let $_{x,\lambda}\Bunb_R\subset {_{x,\infty}\Bunb_R}$ be the closed substack given by the property that 
\begin{equation}
\label{map_s_with_lambda_x} 
s: \cE\otimes\Omega\hook{} M(\<\lambda, \check{\omega}_1\>x)
\end{equation}
is regular. Let $_{x,\lambda}\Bun_R\subset {_{x,\lambda}\Bunb_R}$ be the open substack given by the property that (\ref{map_s_with_lambda_x}) is a subbundle over the whole of $X$. Let the stacks
\begin{equation}
\label{stacks_x_lambda_Bunb_R}
_{x,\lambda}\Bun_{\tilde R}\subset {_{x,\lambda}\Bunb_{\tilde R}}
\end{equation}
be obtained by the base change $\Bunt_G\to\Bun_G$. A point of $_{x,\lambda}\Bun_R$ yields an exact sequence 
\begin{equation}
\label{seq_cE_with_divisor_by_cEOmega} 
 0\to \cE\otimes\Omega(-\<\lambda, \check{\omega}_1\>x)\to M\to \cE^{-1}(\<\lambda, \check{\omega}_1\>x)\to 0,
\end{equation} 
hence a section $i_{\lambda}: {_{x,\lambda}\Bun_R}\to {_{x,\lambda}\Bun_{\tilde R}}$ of the $\mu_2$-gerb $_{x,\lambda}\Bun_{\tilde R}\to {_{x,\lambda}\Bun_R}$ given by 
$$
\cB=\det\RG(X, \cE\otimes\Omega(-\<\lambda, \check{\omega}_1\>x))
$$ 
together with the isomorphism $\cB^2\,\iso\, \det\RG(X, M)$ given by (\ref{seq_cE_with_divisor_by_cEOmega}). 

 If $D+2\<\lambda, \check{\omega}_1\>x \ge 0$ then we get a canonical inclusion $\cE^2\otimes\Omega(-2\<\lambda, \check{\omega}_1\>x)\hook{}\Omega$, let $\ev_{x,\lambda}:{_{x,\lambda}\Bun_R}\to \Bun_{\Omega}$ be the corresponding push-forward map. 

 The section $i_{\lambda}$ yields an isomorphism $_{x,\lambda}\Bun_R\times B(\mu_2)\,\iso\, {_{x,\lambda}\Bun_{\tilde R}}$, so we get the local system $(\ev_{x,\lambda}^*\cL)\boxtimes W$ on $_{x,\lambda}\Bun_{\tilde R}$. We denote by $W$ the nontrivial rank one local system on $B(\mu_2)$ corresponding to the inclusion $\mu_2\hook{}\Qlb^*$. Under the condition $D+2\<\lambda, \check{\omega}_1\>x \ge 0$ we denote by $\cB^{\lambda}$ the intermediate extension of $(\ev_{x,\lambda}^*\cL)\boxtimes W[\dim (_{x,\lambda}\Bun_R)]$ under (\ref{stacks_x_lambda_Bunb_R}). By construction, $\cB^{\lambda}\in \P^{\cL}_-(_{x,\infty}\Bunb_{\tilde R})$. 
 
  We define also $\cB^{\lambda,!}$ and $\cB^{\lambda,*}$ by replacing in the definition of $\cB^{\lambda}$ the intermediate extension by the extension by zero and $*$-extension respectively. By (\cite{FGV}, Proposition~3.3.1) the open immersion (\ref{stacks_x_lambda_Bunb_R}) is affine, so $\cB^{\lambda,!}$ and $\cB^{\lambda,*}$ are perverse sheaves lying in $\P^{\cL}_-(_{x,\infty}\Bunb_{\tilde R})$. If $\lambda=0$ then $_{x,0}\Bun_R=\Bun_R$ and $D+2\<\lambda, \check{\omega}_1\>x \ge 0$, so $\cB^0\in \P^{\cL}_-(_{x,\infty}\Bunb_{\tilde R})$. 
 
\medskip\noindent
1.6  From now on for simplicity we assume that $D$ is a reduced divisor (without multiplicities). This readily implies that $\cB^0$ is the extension by zero under $\Bun_{\tilde R}\hook{} \Bunb_{\tilde R}$. For $\mu\in\Lambda$ 
the condition $D+2\<\mu,\check{\omega}_1\>\ge 0$ is equivalent to $\mu\in\Lambda^+$. Moreover, from the generic $\cL$-equivariance we see that if some $\cB^{\mu}$ has a nonzero $*$-fibre at some point of $_{x,\infty}\Bunb_{\tilde R}$ then this point lies in $_{x,\lambda}\Bun_{\tilde R}$ for some $\lambda\in\Lambda^+$ with $\lambda\le \mu$. One also checks that each irreducible object of $\P^{\cL}(_{x,\infty}\Bunb_{\tilde R})$ is identifed with $\cB^{\mu}$ for some $\mu\in\Lambda^+$. 

\begin{Th}
\label{Th_1}
 For each $\lambda\in\Lambda^+$ one has $\H^{\ra}(\cA_{\lambda}, \cB^0)\,\iso\, \cB^{\lambda}$ canonically.
\end{Th}

 Let $\lambda\in\Lambda^+$. Let $Z^{\le \lambda}$ be the stack classifying $(\cE\otimes\Omega\hook{s'} M')\in\Bun_R$, a modification $M$ of $M'$ at $x$ such that $M$ is in the position $\le \lambda$ with respect to $M'$ at $x$ (we denote by $\beta: M\,\iso\, M'\mid_{X-x}$ the corresponding isomorphism), lines $\cB$, $\cB'$ equipped with $\cB^2\,\iso\, \det\RG(X,M)$ and $\cB'^2\,\iso\, \det\RG(X,M')$.  
The perverse sheaf $(\cB^0\tboxtimes \cA_{\lambda})^r$ is the extension by zero under $Z^{\le\lambda}\hook{} Z$. The map $p_Z$ restricts to a morphism
$$
\pi_{\lambda}: Z^{\le\lambda}\to {_{x,\lambda}\Bunb_{\tilde R}}
$$
sending the above collection to $s: \cE\otimes\Omega\hook{} M(\<\lambda, \check{\omega}_1\>x)$, $\cB$ together with $\cB^2\,\iso\, \det\RG(X,M)$. 

 We adopt the argument of (\cite{FGV}, Theorem~4) for our situation. By decomposition theorem, $\H^{\ra}(\cA_{\lambda}, \cB^0)$ is a direct sum of its perverse cohomology sheaves, and each of its perverse cohomology sheaves lies in $\P^{\cL}(_{x,\infty}\Bunb_{\tilde R})$.
 
 For $\nu\in \Lambda^+$ let $K^{\nu}$ denote the $*$-restriction of $\pi_{\lambda !}((\cB^0\tboxtimes \cA_{\lambda})^r)$ to the stratum $_{x,\nu}\Bun_{\tilde R}$. It suffices to prove the following.
 
\begin{Lm} 
\label{Lm_one}
1) For each $0\le \nu\le \lambda$ the complex $K^{\nu}$ lies in perverse degrees $\le 0$.  \\
2) The $0$-th perverse cohomology sheaf of $K^{\nu}$ vanishes unless $\nu=\lambda$, and in the latter case it is canonically isomorphic to $\cB^{\lambda}$.
\end{Lm}

 For $0\le \nu\le \lambda$ let $Z^{\nu, \le\lambda}=\pi_{\lambda}^{-1}(_{x,\nu}\Bun_{\tilde R})$. For $\mu\in\Lambda^+$ with $\mu\le\lambda$ 
let $Z^{\nu, \mu}\subset Z^{\nu, \le\lambda}$ be the locally closed substack given by the property that $M$ is in the position $\mu$ with respect to $M'$ at $x$. 
Denote by $K^{\nu,\mu}$ the $!$-direct image under 
\begin{equation}
\label{map_p_lambda_restricted_nu_mu}
\pi_{\lambda}: Z^{\nu,\mu}\to {_{x,\nu}\Bun_{\tilde R}}
\end{equation}
 of the $*$-restriction of $(\cB^0\tboxtimes \cA_{\lambda})^r$. Lemma~\ref{Lm_one} is reduced to the following.

\begin{Lm} 1) The complex $K^{\nu,\mu}$ lives in perverse degrees $\le 0$, and the inequality is strict unless $\mu=\nu=\lambda$.\\
2) the 0-th perverse cohomology sheaf of $K^{\lambda,\lambda}$ is canonically isomorphic to $\cB^{\lambda}$. 
\end{Lm}
\begin{Prf} \    The inclusion $Z^{\lambda,\lambda}\subset Z^{\le\lambda}$ is open. The map $\pi_{\lambda}: Z^{\lambda,\lambda}\to {_{x,\lambda}\Bun_{\tilde R}}$ is a $\mu_2$-gerb, and the perverse sheaf $(\cB^0\tboxtimes \cA_{\lambda})^r$ is constant along this gerb. This yields the isomorphism $K^{\lambda,\lambda}\,\iso\, \cB^{\lambda}$.
 
 Note that $Z^{\nu,\mu}$ is empty unless $\nu\le\mu$. If $\mu=\nu$ then (\ref{map_p_lambda_restricted_nu_mu}) is a trivial $\mu_2$-gerb. If $\nu<\mu$ then (\ref{map_p_lambda_restricted_nu_mu}) is a fibration, whose fibre is a trivial $\mu_2$-gerb over 
$$
\A^{\<\mu-\nu, \, \check{\omega}_1\>-1}\times\Gm\,\iso\, S^{-\nu}\cap \Gr^{\mu}_G
$$ 
in the standard notation of (\cite{FGV}, Section~7).  Note that $\check{\rho}=\check{\omega}_1$, where $\check{\rho}$ is the half sum of positive roots.

 The projection $p'_Z: Z^{\nu,\mu}\to \Bun_{\tilde R}$ makes $Z^{\nu,\mu}$ into a fibration over $\Bun_{\tilde R}$ with typical fibre $B(\mu_2)\times (S^{\nu}\cap \Gr^{\mu})$. So, the $*$-restriction $(\cB^0\tboxtimes \cA_{\lambda})^r\mid_{Z^{\nu,\mu}}$ is the twisted exteriour product
\begin{equation}
\label{complex_cB_0_tboxtimes_cA_lambda_restr}
\cB^0\tboxtimes \cA_{\lambda}\mid_{B(\mu_2)\times (S^{\nu}\cap \Gr^{\mu})}
\end{equation}
for this projection. Since $\cA_{\lambda}\mid_{B(\mu_2)\times \Gr^{\mu}}$ lives in perverse degrees $\le 0$ and is a constant complex, it may have usual cohomology sheaves only in degrees $\le -\<\mu, 2\check{\rho}\>$, and the inequality is strict unless $\mu=\lambda$.  So, (\ref{complex_cB_0_tboxtimes_cA_lambda_restr}) is placed in usual degrees
$$
\le -\dim\Bun_R-\<\mu, 2\check{\rho}\>,
$$
and the inequality is strict unless $\mu=\lambda$.

 Since the dimension of the fibres of (\ref{map_p_lambda_restricted_nu_mu}) is $\<\mu-\nu, \check{\rho}\>$, the complex $K^{\nu,\mu}$ is placed in usual degrees $\le -\dim\Bun_R-\<\nu, 2\check{\rho}\>$, and the inequality is strict unless $\mu=\lambda$. One checks that 
$$
\dim{_{x,\nu}\Bun_R}=\deg D+\<\nu, 2\check{\rho}\>+1-g=\dim\Bun_R+\<\nu, 2\check{\rho}\>
$$
So, $K^{\nu,\mu}$ is placed in perverse degrees $\le 0$, and the inequality is strict unless $\mu=\lambda$.

 It remains to study the $0$-th perverse cohomology of $K^{\nu,\lambda}$ for $\nu<\lambda$. Over $Z^{\nu,\lambda}$ the $*$-restriction $(\cB^0\tboxtimes \cA_{\lambda})^r\mid_{Z^{\nu,\lambda}}$ is a shifted local system. This local system is always nontrivial. If $x\notin D$ then according to (\cite{FGV}, Proposition~7.1.7) the contribution of $\ev_{x,\lambda}^*\cL$ to this local system is notnrivial, and our $\mu_2$-gerb could only contribute to this local system by tensoring it with another local system of rank one and order two. So, for $\nu<\lambda$ the 0-th perverse cohomology of $K^{\nu,\lambda}$ vanishes. If $x\in D$ then the contribution of $\ev_{x,\lambda}^*\cL$ may be trivial (this is the case for $\nu=0$ and $\lambda=\alpha$, for example), in which case  the $*$-restriction $(\cB^0\tboxtimes \cA_{\lambda})^r\mid_{Z^{\nu,\lambda}}$ is a nontrivial (shifted) local system of rank one and order 2 coming from the $\mu_2$-gerb.
We are done.
\end{Prf}  

\medskip

Theorem~\ref{Th_1} is proved. Now the arguments of (\cite{FGV}, Theorem~3(1) and 3(2)) apply without changes and prove the following.

\begin{Th} 
\label{Th_2}
1) The category $\P^{\cL}(_{x,\infty}\Bunb_{\tilde R})$ is semi-simple. For $\lambda,\mu\in\Lambda^+$ one has $\Ext^1(\cB^{\lambda},\cB^{\mu})=0$, where $\Ext$ is calculated in the category $\P(_{x,\infty}\Bunb_{\tilde R})$ of all perverse sheaves. \\
2) For any $\lambda\in\Lambda^+$ the natural maps 
$$
\cB^{\lambda, !}\to \cB^{\lambda}\to \cB^{\lambda,*}
$$ 
are isomorphisms.
\end{Th} 

\begin{Cor} 
\label{Cor_free_Whit_module}
The functor $\H^{\ra}(\cdot, \cB^0): \Rep(\SL_2)\to \P^{\cL}(_{x,\infty}\Bunb_{\tilde R})$ is an equivalence.
\end{Cor}

 Since $\P^{\cL}(_{x,\infty}\Bunb_{\tilde R})\subset \P(_{x,\infty}\Bunb_{\tilde R})$ is a Serre subcategory (in the strong sense, that is, closed under taking the subquotients and extensions), one can define 
$$
\D^{\cL}(_{x,\infty}\Bunb_{\tilde R})\subset \D(_{x,\infty}\Bunb_{\tilde R})
$$ 
as the full subcategory whose all perverse cohomology sheaves lie in $\P^{\cL}(_{x,\infty}\Bunb_{\tilde R})$.
Similarly for $\D^{\cL}_-(_{x,\infty}\Bunb_{\tilde R})$,  compare with (\cite{G}, Section~4).

\bigskip

\centerline{\scshape 2. Whittaker categories in families}

\bigskip\noindent
2.1 We change the following notations compared to Section~1. We do not assume any more that $D, \cE$ and an isomorphism $\cE^2\,\iso\, \cO(-D)$ are given. These objects will have a new meaning in this section.

Let $\cQ_1$ be the stack classifying $(M,\cB)\in\Bunt_G$, $\cE\in\Bun_1$ equipped with a section $s: \cE\otimes\Omega\hook{} M$. Let $\nu_{\cQ}: \cQ_1\to \Bunt_G$ denote the projection forgetting $\cE,s$. 

Let $\cQ_{2,ex}$ be the stack classifying a point of $\cQ_1$ as above together with a section $s_2: \cE^2\to\cO$. Let $\cQ_2\subset \cQ_{2,ex}$ be the open substack given by the property that $s_2\ne 0$. Let $\pi_{2,ex,1}: \cQ_{2,ex}\to \cQ_1$ be the projection forgetting $s_2$. Let $\pi_{2,1}: \cQ_2\to \cQ_1$ be its restriction to $\cQ_2$. 

 We will define Whittaker categories and Whittaker functors for $\Bunt_G$ in the style of (\cite{G}, Sections~4, 5). Our first aim is to define a full triangulated subcategory $\D^{\cL}(\cQ_{2,ex})\subset \D(\cQ_{2,ex})$, which we think of as the family of Whittaker categories from Section~1. It will be given by the condition that $K\in \D(\cQ_{2,ex})$ lies in $\D^{\cL}(\cQ_{2,ex})$ if all its perverse cohomology sheaves lie in some Serre subcategory $\P^{\cL}(\cQ_{2,ex})\subset\P(\cQ_2)$. (We mean here a Serre subcategory in the strong sense, that is, a strictly full subcategory closed under taking subquotients and taking extensions). 
 
  We could define these categories repeating the definition of generic Hecke $\cL$-equivariance from Section~1.2, but we will give a definition in the style of (\cite{G}, Section~4).
  
\medskip\noindent  
2.2  For any $x\in X$ let $\cQ_1^{x}\subset \cQ_1$ be the open substack given by the property that $s:\cE\otimes\Omega\hook{} M$ is a subbundle in a neighbourhood of $x$. Over the formal neighbourhood $D_x$ of $x$ in $X$ a point of $\cQ_1^x$ yields an exact sequence (\ref{ext_cE^*_by_EotimesOmega}). We consider it as an $R$-torsor $\cF_R$ over $D_x$, where now $R$ denotes the group scheme $\cE^2\otimes\Omega$ on $X$. Let $D_x^*$ denote the punched formal neighbourhood of $x$. 
  
 Let $\cH^x_1$ be the stack classifying a collection $(M,\cB, \cE, s)\in \cQ_1^x$ and an element 
$$
\sigma\in \H^0(D_x^*, R)/\H^0(D_x, R)
$$
The element $\sigma$ gives rise to an $R$-torsor $\cF^1_R$ over $X$ trivialized over $X-x$. Let $\cF'_R$ denote the $R$-torsor over $D_x$ obtained as the sum of $\cF_R$ with $\cF^1_R$. We consider $\cF'_R$ as an exact sequence
$$
0\to \cE\otimes\Omega\to M'\to \cE^{-1}\to 0
$$
over $D_x$ equipped with an isomorphism of $R$-torsors $\tau: \cF_R\,\iso\, \cF'_R\mid_{D_x^*}$. We also denote by $M'$ the $G$-torsor on $X$ obtained as the gluing of $M'\mid_{D_x}$ with $M\mid_{X-x}$ via $\tau$ over $D_x^*$. So, $\tau$ yields an isomorphism of $G$-torsors
$$
\tau: M_{X-x}\,\iso\, M'\mid_{X-x},
$$
hence also a $\ZZ/2\ZZ$-graded isomorphisms (\ref{iso_detRG_modification}) and (\ref{iso_detRG_M_and_M'}) as in Section~1.2. Now (\ref{iso_detRG_M_and_M'}) yields also a $\ZZ/2\ZZ$-graded isomorphism $\cB^2\,\iso\,\det\RG(X, M')$. Let $s'$ denote the composition $\cE\otimes\Omega\hook{s} M \toup{\tau} M'$ then $s': \cE\otimes\Omega\to M'$ is regular over $X$. We get the diagram
\begin{equation}
\label{diag_for_cH_x^1}
\cQ_1^x \;\getsup{h^{\la}_{1}}\; \cH^x_1\; \toup{h^{\ra}_{1}} \;\cQ_1^x,
\end{equation}
where $h^{\la}_{1}$ (resp., $h^{\ra}_{1}$) sends the above point to $(M,\cB, \cE,s)$ (resp., to $(M',\cB, \cE, s')$). In this way $\cH^x_1$ is naturally a groupoid over $\cQ^x_1$. 

 For $i\ge 0$ let $^i\cH^x_1\subset \cH^x_1$ be the closed substack given by the property that 
$$
\sigma\in \cE^2\otimes\Omega(ix)/\cE^2\otimes\Omega(\cO_x)
$$ 

 Let $\cQ^x_{2,ex}=\cQ^x_1\times_{\cQ_1}\cQ_{2,ex}$. Consider the stack $\cH^x_1\times_{\cQ^x_1} \cQ^x_{2,ex}$, where we used $h^{\la}_{1}$ to define the fibred product. In the same way, $\cH^x_1\times_{\cQ^x_1} \cQ^x_{2,ex}$ is a groupoid over $\cQ^x_{2,ex}$, we denote the corresponding diagram by
$$
\cQ_{2,ex}^x\;\getsup{h^{\la}_2} \; \cH^x_1\times_{\cQ^x_1} \cQ^x_{2,ex} \; \toup{h^{\ra}_2}\; \cQ_{2,ex}^x
$$
  
  Let $\ev_{\cQ}^x: \cH^x_1\times_{\cQ^x_1} \cQ^x_{2,ex}\to \Bun_{\Omega}$ be the following map. Applying $s_2$ to $\sigma$, one gets an element 
$$
s_2(\sigma)\in \H^0(D_x^*, \Omega)/\H^0(D_x,\Omega)
$$ 
Let $\cF_{\Omega}$ denote the $\Omega$-torsor on $X$ obtained as the gluing of the trivial $\Omega$-torsors over $D_x$ and $X-x$ via the gluing datum $s_2(\sigma)$. 

 We can now define the full subcategory $\P^{\cL}(\cQ^x_{2,ex})\subset \P(\cQ_{2,ex})$ as the one consisting of those perverse sheaves $K$ on $\cQ_{2,ex}$ for which for any $i\ge 0$ there is an isomorphism
\begin{equation}
\label{iso_for_def_cat_on_cQx_2_ex}
(h_2^{\ra})^*K\,\iso\, (h^{\la}_2)^*K\otimes (\ev_{\cQ}^x)^*\cL
\end{equation}
over $^i\cH^x_1\times_{\cQ^x_1} \cQ^x_{2,ex}$, whose restriction to the unit section is the identity map. For each $i$ and $K\in \P(\cQ_{2,ex})$ there could be at most one such isomorphism. If it exists, this family of isomorphisms is automatically compatible. 

 Since the projections
$$
h^{\la}_2, h^{\ra}_2: {^i\cH^x_1\times_{\cQ^x_1} \cQ^x_{2,ex}}\to \cQ^x_{2,ex}
$$
are affine fibrations, $\P^{\cL}(\cQ^x_{2,ex})$ is a Serre subcategory in $\P(\cQ^x_{2,ex})$. Define $\D^{\cL}(\cQ^x_{2,ex})\subset \D(\cQ^x_{2,ex})$ as the full subcategory of complexes whose all perverse cohomology lie in $\P^{\cL}(\cQ^x_{2,ex})$. By (\cite{G}, Lemma~4.8) for each $K\in \D^{\cL}(\cQ^x_{2,ex})$ and $i\ge 0$ over $^i\cH^x_1\times_{\cQ^x_1} \cQ^x_{2,ex}$ there is an isomorphism (\ref{iso_for_def_cat_on_cQx_2_ex}), whose restriction to the unit section is the identity map.
  
\medskip\noindent  
2.3 If now $\bar x=x_1,\ldots, x_m$ is a collection of pairwise distinct points on $X$, let $\cQ_1^{\bar x}\subset\cQ_1$ be the open substack given by requiring that $s: \cE\otimes\Omega\hook{}M$ is a subbundle in a neighbourhood of each of $x_i$. In a similar way, one defines the Serre subcategory $\P^{\cL}(\cQ^{\bar x}_{2,ex})\subset \P(\cQ^{\bar x}_{2,ex})$ and a full triangulated subcategory $\D^{\cL}(\cQ^{\bar x}_{2,ex})\subset \D(\cQ^{\bar x}_{2,ex})$. 

\medskip\noindent
2.4 For $d\ge 0$ let $^d\cQ_1\subset \cQ_1$ be the locally closed substack given by the property that there is $D_1\in X^{(d)}$ such that $s: \cE\otimes\Omega(D_1)\hook{}M$ is a subbudle. For a point of $^d\cQ_1$ we get an exact sequence on $X$
\begin{equation}
\label{ext_cE*_by_cEOmega_with_D_1}
0\to \cE\otimes\Omega(D_1)\to M\to \cE^{-1}(-D_1)\to 0
\end{equation}

 Let $^d\cQ_1^x={^d\cQ_1}\cap\cQ_1^x$. Let $^d\cQ_{2,ex}^x\subset {^d\cQ_{2,ex}}$ denote the preimages of $^d\cQ_1^x\subset {^d\cQ_1}$ in $\cQ_{2,ex}$. The substack $^d\cQ_{2,ex}^x\subset {\cQ_{2,ex}^x}$ is stable under the action of the groupoid $\cH^x_1\times_{\cQ_1^x}\cQ_{2,ex}^x$, so we can define the categories 
$$
\P^{\cL}(^d\cQ_{2,ex}^x)\subset\D^{\cL}(^d\cQ_{2,ex}^x)
$$ 
as above. The analog of (\cite{G}, Lemma~4.11) holds, namely we get the following.

\begin{Lm} 1) The $*$ and $!$-restrictions map $\D^{\cL}(\cQ^x_{2,ex})$ to $\D^{\cL}(^d\cQ^x_{2,ex})$.\\
2) The $*$ and $!$-direct images map $\D^{\cL}(^d\cQ^x_{2,ex})$ to $\D^{\cL}(\cQ^x_{2,ex})$.\\
3) An object $K\in \D(\cQ^x_{2,ex})$ lies in $\D^{\cL}(\cQ^x_{2,ex})$ iff its $*$-restriction (or, equivalently, its $!$-restriction) lies to $^d\cQ^x_{2,ex}$ lies 
in $\D^{\cL}(^d\cQ^x_{2,ex})$ for each $d\ge 0$. 
\end{Lm}
  
  For a point of $^q\cQ_1$ the exact sequence (\ref{ext_cE*_by_cEOmega_with_D_1}) yields an isomorphism 
$$
\det\RG(X,M)\,\iso\, \det\RG(X, \cE\otimes\Omega(D_1))^{\otimes 2}
$$ 
So, $^d\cQ_1$ can be seen as the stack classifying: $D_1\in X^{(d)}$, $\cE\in\Bun_1$, an exact sequence (\ref{ext_cE*_by_cEOmega_with_D_1}) on $X$, a line $\cB_0$ equipped with $\cB_0^2\,\iso\, k$ giving a point of $B(\mu_2)$. For these data we let $\cB=\cB_0\otimes \det\RG(X, \cE\otimes\Omega(D_1))$ equipped with the induced isomorphism $\cB^2\,\iso\, \det\RG(X,M)$. 
  
  The stack $^d\cQ_{2,ex}^x$ classifies the same data as for $^d\cQ_1$ such that $x\notin D_1$ together with a morphism $s_2:\cE^2\to\cO$. Let 
\begin{equation}
\label{def_closed_part_of_dcQ_2ex}
^d\cQ_{2,ex}^{' x}\subset {^d\cQ_{2,ex}^x}
\end{equation}
be the closed substack given by the property that $s_2$ extends to a regular map $\cE^2(2D_1)\to\cO$ over $X$. Let 
$$
^d\ev^x: {^d\cQ_{2,ex}^{' x}}\to \Bun_{\Omega}
$$ 
be the map sending the above point to the push-forward of (\ref{ext_cE*_by_cEOmega_with_D_1}) by $s_2$. 

 Let $^d\cP^x$ be the stack classifying $\cE\in\Bun_1$, $D_1\in (X-x)^{(d)}$ with a section $s_2: \cE^2(2D_1)\to\cO$, a line $\cB_0$ equipped with $\cB_0^2\,\iso\, k$. Let $^d\pr^x: {^d\cQ_{2,ex}^{' x}}\to {^d\cP^x}$ be the map sending the above point to $(\cE, D_1, s_2,\cB_0)$. 

\begin{Lm} 
\label{Lm_support_condition_for_x}
1) Any object of $\D^{\cL}(^d\cQ_{2,ex}^x)$ is the extension by zero under (\ref{def_closed_part_of_dcQ_2ex}). \\
2) The functor 
$$
K\mapsto (^d\pr^x)^*K\otimes (^d\ev^x)^*\cL[\dimrel(^d\pr^x)]
$$ 
is an equivalence of triangulated categories $\D(^d\cP^x)\,\iso\, \D^{\cL}(^d\cQ_{2,ex}^x)$ exact for the perverse t-structures.
\end{Lm}

 As in (\cite{G}, Corollary~4.14) Lemma~\ref{Lm_support_condition_for_x} implies the following.
 
\begin{Cor}  Let $\bar x=x_1,\ldots, x_m$ be a collection of pairwise distinct points of $X$ and $x$ be one of them. Then the restriction functor maps 
$\D^{\cL}(\cQ_{2,ex}^x)$ to $\D^{\cL}(\cQ_{2,ex}^{\bar x})$.
\end{Cor}

 Finally, define $\P^{\cL}(\cQ_{2,ex})\subset\P(\cQ_{2,ex})$ as the full subcategory of those $K\in \P(\cQ_{2,ex})$ such that for each finite collection $\bar x$ the restriction $K\mid_{\cQ_{2,ex}^{\bar x}}$ lies in $\P^{\cL}(\cQ_{2,ex}^{\bar x})$. As in \cite{G}, this is a Serre subcategory of $\P(\cQ_{2,ex})$. We define 
$$
\D^{\cL}(\cQ_{2,ex})\subset \D(\cQ_{2,ex})
$$ 
as the full subcategory of complexes whose all perverse cohomologies lie in $\P^{\cL}(\cQ_{2,ex})$.
 
\medskip\noindent 
2.5   Note that $\cQ_2\subset\cQ_{2,ex}$ is stable under the action of our groupoids, so we similarly get the categories $\P^{\cL}(\cQ_2)\subset\D^{\cL}(\cQ_2)$. The zero section $s_{1,2}: \cQ_1\to \cQ_{2,ex}$ of $\pi_{2,ex,1}$ is a closed substack stable under the action of our groupoids, so we have the ful subcategory $\P^{\cL}(\cQ_1)\subset \P^{\cL}(\cQ_{2,ex})$, whose objects are those of $\P^{\cL}(\cQ_{2,ex})$ which are etxensions by zero under $s_{1,2}$.
Similarly, we have $\D^{\cL}(\cQ_1)\subset \D^{\cL}(\cQ_{2,ex})$.   

\bigskip

\centerline{\scshape 3. Whittaker functors}

\bigskip\noindent
3.1 As in \cite{G}, we will construct a functor 
\begin{equation}
\label{functor_W_ex}
W_{ex}: \D(\cQ_1)\to\D^{\cL}(\cQ_{2,ex})
\end{equation}
satisfying the following.

\begin{Pp} 
\label{Pp_one}
The functor (\ref{functor_W_ex}) is an equivalent identifying $\P(\cQ_1)$ with $\P^{\cL}(\cQ_{2,ex})$. The quasi-inverse functor is given by 
$K\mapsto (\pi_{2,ex,1})_!K$. Moreover, for each $K\in \D^{\cL}(\cQ_{2,ex})$ the natural map $(\pi_{2,ex,1})_!K\to (\pi_{2,ex,1})_*K$ is an isomorphism.
\end{Pp}

\medskip\noindent
3.2 For each $d\ge 0$ we have the subcategories $\P^{\cL}(^d\cQ_{2,ex})\subset \P(^d\cQ_{2,ex})$ and $\D^{\cL}(^d\cQ_{2,ex})\subset \D(^d\cQ_{2,ex})$ defined as above. As in \cite{G}, we first describe the functor $W_{ex}$ on strata. Namely, for $d\ge 0$ we will define an equivalence of triangulated categories
$$
^dW_{ex}: \D(^d\cQ_1)\,\iso\, \D^{\cL}(^d\cQ_{2,ex})
$$
exact for the perverse t-structures. 

 Remind that $^d\cQ_{2,ex}$ classifies $\cE\in\Bun_1$, $D_1\in X^{(d)}$, $s_2: \cE^2\to\cO$,  a line $\cB_0$ equipped with $\cB_0^2\,\iso\, k$, and an exact sequence (\ref{ext_cE*_by_cEOmega_with_D_1}). Let $^d\cQ'_{2,ex}\subset {^d\cQ_{2,ex}}$ be the closed substack given by the property that $s_2:\cE^2(2D_1)\to\cO$ is regular over $X$. Let $^d\ev: {^d\cQ'_{2,ex}}\to \Bun_{\Omega}$ be the map sending the above point to the push-forward of (\ref{ext_cE*_by_cEOmega_with_D_1}) by $s_2$.
 
  Let $^d\cP$ be the stack classifying $\cE\in\Bun_1$, $D_1\in X^{(d)}$ with $s_2: \cE^2(2D_1)\to\cO$, a line $\cB_0$ equipped with $\cB_0^2\,\iso\, k$. Let $^d\pr: {^d\cQ'_{2,ex}}\to {^d\cP}$ be the map sending the above point to $(\cE, D_1, s_2,\cB_0)$. 
  
\begin{Lm} 
\label{Lm_description_D_cL_d_cQ_2ex}
Any object of $\D^{\cL}(^d\cQ_{2,ex})$ is the extension by zero from $^d\cQ'_{2,ex}$. The functor 
$$
K\mapsto (^d\pr)^*K\otimes (^d\ev)^*\cL[\dimrel(^d\pr)]
$$
is an equivalence $\D(^d\cP)\,\iso\, \D^{\cL}(^d\cQ_{2,ex})$ exact for the perverse t-structures.
\end{Lm}
  
  Let $^d\bar\cP$ be the stack classifying $\cE\in\Bun_1$, $D_1\in X^{(d)}$, a line $\cB_0$ equipped with $\cB_0^2\,\iso\, k$. Then $^d\cQ_1$ and $^d\cP$ are naturally dual (generalized) vector bundles over $^d\bar\cP$. Define $^dW_{ex}$ as the composition
$$
\D(^d\cQ_1)\toup{\Four_{\psi}} \D(^d\cP)\,\iso\, \D^{\cL}(^d\cQ_{2,ex}),
$$
where $\Four_{\psi}$ is the Fourier transform, and the second functor is the equivalence of Lemma~\ref{Lm_description_D_cL_d_cQ_2ex}.   

 Let $^d\pi_{2,ex,1}: {^d\cQ_{2,ex}}\to {^d\cQ_1}$ be the restriction of $\pi_{2,ex,1}$. The argument of (\cite{G}, Lemma~5.3) applies and gives the following.
 
\begin{Lm} The functor $(^d\pi_{2,ex,1})_!: \D^{\cL}(^d\cQ_{2,ex})\to \D(^d\cQ_1)$ is an equivalence quasi-inverse to $^dW_{ex}$. Besides, for $K\in \D^{\cL}(^d\cQ_{2,ex})$ the natural map $(^d\pi_{2,ex,1})_! K\to (^d\pi_{2,ex,1})_* K$ is an isomorphism.
\end{Lm}

\medskip\noindent
3.3 For $x\in X$ let us construct now the functor $W^x_{ex}: \D(\cQ_1^x)\to \D^{\cL}(\cQ^x_{2,ex})$. Remind the diagram (\ref{diag_for_cH_x^1}) and the closed substack $^i\cH^x_1\subset \cH^x_1$ from Section~2.2. Write also $^iE={^i\cH^x_1}$, this is a vector bundle over $\cQ_1^x$ with fibre $\cE^2\otimes\Omega(ix)/\cE^2\otimes\Omega$ over $(M,\cB,\cE, s)$. 

 Let $^i\check{E}\to \cQ^x_1$ denote the dual vector bundle, its fibre over $(M,\cB,\cE, s)$ is $\cE^{-2}/\cE^{-2}(-ix)$ by Serre duality. We get a map $f_i: \cQ_{2,ex}^x\to {^i\check{E}}$ sending $s_2$ to its image under $\H^0(X,\cE^{-2})\to \cE^{-2}/\cE^{-2}(-ix)$. Over the open substack $\cU\subset \cQ^x_1$ given by $\H^0(X, \cE^{-2}(-ix))=0$, the map $f_i$ is a closed immersion. Besides, for $j\ge i$ we get a diagram
$$
\begin{array}{ccc}
\cQ_{2,ex}^x & \toup{f_j} &  {^j\check{E}}\\
& \searrow\lefteqn{\scriptstyle f_i} & \downarrow\lefteqn{\scriptstyle \pr_{j,i}}\\
&& ^i\check{E},
\end{array}
$$
where $\pr_{j,i}$ is a natural projection. Let $^ih^{\ra}_1: {^i\check{E}}\to \cQ^x_1$ denote the restriction of $h^{\ra}_1$, the map $^ih^{\ra}_1$ is smooth. Let $^iW^x_{ex}: \D(\cQ^x_1)\to \D(^i\check{E})$ denote the functor 
$$
K\mapsto \Four_{\psi}(^ih^{\ra}_1)^*K[\dimrel(^ih^{\ra}_1)],
$$
it is exact for the perverse t-structures. From the standard properties of the Fourier transform, we get for $j\ge i$ an isomorphism functorial in $K\in \D(\cQ^x_1)$
\begin{equation}
\label{iso_ji_for_iWx}
(\pr_{j,i})_! {^jW^x_{ex}}(K)\,\iso\, {^iW^x_{ex}}(K)
\end{equation}

\begin{Lm} 
\label{Lm_support_for_iWx}
Let $K\in \D(\cQ^x_1)$ and $i\ge 0$. Over the open substack $\cU\subset \cQ^x_1$ given by $\H^0(X, \cE^{-2}(-ix))=0$, the complex $^iW^x_{ex}(K)$ is the extension by zero under the closed immersion $f_i: \cQ^x_{2,ex}\to {^i\check{E}}$.
\end{Lm}
\begin{Prf} For $d\ge 0$ let $^{d,i}E={^i E}\times_{\cQ^x_1} {^d\cQ^x_1}$, where we used $h^{\la}_1$ to define the fibred product, set also 
$$
^{d,i}\check{E}={^i\check{E}}\times_{\cQ_x^1}{^d\cQ^x_1}
$$
It suffices to show that for any $d\ge 0$ over $\cU$ the $*$-restriction of $^iW^x_{ex}(K)$ to $^{d,i}\check{E}$ is the extension by zero under 
$f_i: {^d\cQ^x_{2,ex}} \to {^{d,i}\check{E}}$.

 The stack $^d\cQ^x_1$ classifies: $D_1\in (X-x)^{(d)}$, $\cE\in\Bun_1$, an exact sequence (\ref{ext_cE*_by_cEOmega_with_D_1}), and a line $\cB_0$ equipped with $\cB_0^2\,\iso\, k$. 
 
 Let $^dZ^x$ be the stack classifying: $D_1\in (X-x)^{(d)}$, $\cE\in\Bun_1$,  exact sequences (\ref{ext_cE*_by_cEOmega_with_D_1}) and 
\begin{equation}
\label{ext_M_1_with_D_1}
0\to \cE\otimes\Omega(D_1)\to M_1\to \cE^{-1}(-D_1)\to 0
\end{equation} 
over $X$, and a line $\cB_0$ equipped with $\cB_0^2\,\iso\, k$. The map  $h^{\ra}_1: {^{d,i} E}\to {^d\cQ^x_1}$ can be written as the composition
$$
{^{d,i} E}\toup{\gamma} {^dZ^x}\toup{\delta} {^d\cQ^x_1},
$$
where $\delta$ is defined as the sum of the extension (\ref{ext_cE*_by_cEOmega_with_D_1}) with (\ref{ext_M_1_with_D_1}). The map $\gamma$ sends the above point of $^d\cQ^x_1$ together with $\sigma\in \cE^2\otimes\Omega(ix)/\cE^2\otimes\Omega$ to the same point of $^d\cQ^x_1$ together with the exact sequence (\ref{ext_M_1_with_D_1}) obtained from the natural map
$$
\alpha: \cE^2\otimes\Omega(2D_1+ix)/\cE^2\otimes\Omega(2D_1)\to \Bun_{\cE^2\otimes\Omega(2D_1)}
$$
Over $\cU$, the map $\alpha$ is surjective. The dual morphism is the restriction
$$
\check{\alpha}: \H^0(X, \cE^{-2}(-2D_1))\to \cE^{-2}/\cE^{-2}(-ix)
$$ 
Clearly, it factors through $f_i: \H^0(X, \cE^{-2})\to \cE^{-2}/\cE^{-2}(-ix)$. Since the Fourier transform exchanges the $*$-inverse image with the $!$-direct image, our assertion follows.
\end{Prf}

\medskip

 Lemma~\ref{Lm_support_for_iWx} together with the isomorphisms (\ref{iso_ji_for_iWx}) yield the desired functor $W^x_{ex}$. Over any substack of finite type $\cV\subset \cQ^x_1$, the complex $W^x_{ex}(K)$ is defined as $^i W^x_{ex}(K)$ for $i$ large enough depending on $\cV$. As in (\cite{G}, Proposition~5.10), one proves the following.
 
\begin{Lm} 
\label{Lm_functor_Wx_ex}
The functor $K\mapsto (\pi_{2,ex,1})_! K$ is quasi-inverse to $W^x_{ex}: \D(\cQ^x_1)\to \D^{\cL}(\cQ^x_{2,ex})$. Moreover, for any $K\in \D^{\cL}(\cQ^x_{2,ex})$ the natural map $(\pi_{2,ex,1})_! K\to (\pi_{2,ex,1})_* K$ is an isomorphism.
\end{Lm}

\medskip\noindent
3.4 One finishes the proof of Proposition~\ref{Pp_one} as in (\cite{G}, Section~5.11). Namely, first one gets the functor $W_{ex}: \P(\cQ_1)\to \P^{\cL}(\cQ_{2,ex})$ as the gluing of the functors $W^x_{ex}$ for all $x\in X$ using Lemma~\ref{Lm_functor_Wx_ex}. This lemma also implies that $(\pi_{2,ex,1})_!$ is an equivalence $\P^{\cL}(\cQ_{2,ex})\to \P(\cQ_1)$. Finally, as in (\cite{G}, Section~5.11), one shows that 
$$
(\pi_{2,ex,1})_! : \D^{\cL}(\cQ_{2,ex})\to \D(\cQ_1)
$$
is also an equivalence. Proposition~\ref{Pp_one} is proved.

\bigskip

\centerline{\scshape 4. Geometric Eisenstein series for $\wt\SL_2$}

\bigskip\noindent
4.1.1 Keep the notations of Section~0.3.1. The stack $\Bun_P$ classifies $\cE\in\Bun_n$ and an exact sequence 
\begin{equation}
\label{seq_cE*_by_cE_otimes_Omega}
0\to \Sym^2 \cE\to ?\to \Omega^{-1}\to 0
\end{equation}
The above sequence yields $0\to \cE\otimes\Omega\to M\to \cE^*\to 0$. 
Let $\nu_P: \Bun_P\to \Bun_G$ be the map sending (\ref{seq_cE*_by_cE_otimes_Omega}) to $M$. Set $\Bun_{\tilde P}=\Bun_P\times_{\Bun_G}\Bunt_G$. 

 Let $\Bunt_P$ be the stack classifying $\cE\in\Bun_n$, $M\in\Bun_G$ and an inclusion $s: \cE\otimes\Omega\hook{} M$ of coherent sheaves whose image is 
isotropic. This notation agrees with the partial compactification of $\Bun_P$ from (\cite{BG}, Section~1.3.6). So, $\Bun_P\subset\Bunt_P$ is the open substack given by requiring that $\cE\otimes\Omega$ is a subbundle. 

 Set $\Bunt_{\tilde P}=\Bunt_P\times_{\Bun_G}\Bunt_G$. Let $j_P: \Bun_{\tilde P}\to \Bunt_{\tilde P}$ be the natural open immersion. 

  Let $\epsilon$ also denote the 2-automorphism of $\Bunt_{\tilde P}$ acting on $(s, \cE, M, \cB)$ so that it acts trivially on $(s, \cE, M)$ and as $-1$ on $\cB$. Let $\D_-(\Bunt_{\tilde P})\subset \D(\Bunt_{\tilde P})$ be the full subcategory of objects on which $\epsilon$ acts as $-1$. 
  
  Let $\bar\epsilon$ be the 2-automorphism of $\Bunt_{\tilde P}$ acting on $(s, \cE, M, \cB)$ so that it acts trivially on $\cB$ and as -1 on $(\cE, M)$. As in Section~0.3.1, one defines the categories $\D_{=}(\Bunt_{\tilde P})$ and $\D_{\mp}(\Bunt_{\tilde P})$, so
$$
\D_-(\Bunt_{\tilde P})\,\iso\, \D_{=}(\Bunt_{\tilde P})\times \D_{\mp}(\Bunt_{\tilde P})
$$
  
 The stack $\Bun_{\tilde P}$ classifies $\cE\in\Bun_n$, an exact sequence (\ref{seq_cE*_by_cE_otimes_Omega}), and a line $\cB$ equipped with $\cB^2\,\iso\,\det\RG(X,M)$. We have the isomorphism 
$$
\Bun_P\times B(\mu_2)\,\iso\, \Bun_{\tilde P}
$$ 
sending $\cE$, $\cB_0$ equipped with $\cB_0^2\,\iso\, k$ and (\ref{seq_cE*_by_cE_otimes_Omega}) to (\ref{seq_cE*_by_cE_otimes_Omega}) with $\cB=\cB_0\otimes \det\RG(X, \cE\otimes\Omega)$ equipped with the induced isomorphism
$\cB^2\,\iso\, \det\RG(X,M)$. Remind that $W$ denotes the local system on $B(\mu_2)$ corresponding to the nontrivial character $\mu_2\hook{}\Qlb^*$. 
Set 
$$
\IC_{P,W}=(j_P)_{!*}(\IC(\Bun_P)\boxtimes W)\;\;\;\;\mbox{and}\;\;\;\;
\IC_{P,W !}=(j_P)_!(\IC(\Bun_P)\boxtimes W)
$$
These are objects of $\D_-(\Bunt_{\tilde P})$. Note that $\bar\epsilon$ acts on $\det\RG(X, \cE\otimes\Omega)$ as $(-1)^{\chi(\cE\otimes\Omega)}$. So, they lie in $\D_{\mp}(\Bunt_{\tilde P})$ over the connected components with $\chi(\cE\otimes\Omega)$ even and in $\D_{=}(\Bunt_{\tilde P})$ over the connected components with $\chi(\cE\otimes\Omega)$ odd.
 
\medskip\noindent 
4.1.2 Consider the diagram
\begin{equation}
\label{diag_Eis_construction}
\Bun_n \,\getsup{\bar \gq}\, \Bunt_{\tilde P}\, \toup{\bar\gp}\, \Bunt_G,
\end{equation}
where $\bar\gq$ (resp., $\bar\gp$) sends the above point to $\cE\in\Bun_n$ (resp., to $(M,\cB)$). The map $\bar\gp$ is compatible with the 2-automorphisms $\epsilon, \bar\epsilon$. Let $\gq: \Bun_{\tilde P}\to \Bun_n$ be the restriction of $\bar\gq$.

 We define the functor $\bar\gq^{!*}: \D(\Bun_n)\to \D(\Bunt_{\tilde P})$ by
$$
\bar\gq^{!*}(K)=\IC_{P,W}\otimes \bar q^*(K)[-\dim\Bun_n]
$$
It takes values in $\D_-(\Bunt_{\tilde P})$. Set also $\gq^{!*}=j_P^*\comp\bar\gq^{!*}$. One easily adopts the proof of (\cite{BG}, Theorem~5.1.5) to our case to get the following.

\begin{Pp} The complexes $\IC_{P,W}$ and $\IC_{P,W !}$ are
ULA with respect to $\bar\gq: \Bunt_{\tilde P}\to \Bun_n$.
\end{Pp}
\begin{Prf} We will give a proof for $n=1$, the case of any $n$ is left to a reader.
Remind that $\alpha$ denotes the unique coroot of $G$. Let $m\ge 2g-1$, let $\cH^?_G$ be the stack classifying $G$-torsors $M,M'$ on $X$, points $x_1,\ldots, x_m\in X^m-\vartriangle$ together with an isomorphism of $G$-torsors $\tau: M\,\iso\, M'\mid_{X-\{x_1,\ldots, x_m\}}$ such that $M'$ is in the position $\alpha$ with respect to $M$ at each $x_i$, and a line $\cB'$ equipped with $\cB^{' 2}\,\iso\, \det\RG(X,M')$. Here $\vartriangle\subset X^m$ is the divisor of diagonals.
 
  For a point of $\cH^?_G$ denote $\bar M=M\cap M'$ then $M'/M\,\iso\, \cO_{x_1+\ldots+x_m}$ non canonically, and there is a line bundle over $\cH^?_G$ whose fibre at the above point is $\det\H^0(X, M/\bar M)$. For a point of $\cH^?_G$
we get a canonical $\ZZ/2\ZZ$-graded isomorphisms 
$$
\det(M:M')\,\iso\, \det\H^0(X, M/\bar M)^{\otimes 2}
$$ 
and $\cB^2\,\iso\,\det\RG(X,M)$, where $\cB=\cB'\otimes \det\H^0(X, M/\bar M)$. The latter isomorpism is induced by (\ref{iso_detRG_modification}). We get a diagram of projections 
$$
\Bunt_G \,\getsup{h^{\la}}\, \cH^?_G \,\toup{h^{\ra}}\, \Bunt_G,
$$
where $h^{\la}$ (resp., $h^{\ra}$) sends the above point to $(M,\cB)$ (resp., to $(M',\cB')$). 

 Set $\bar Z=\cH^?_G\times_{\Bunt_G} \Bunb_{\tilde B}$, where we used $h^{\ra}$ to define the fibred product, let $h^{\ra}_Z:\bar Z\to \Bunb_{\tilde B}$ be the second projection. Write a point of $\bar Z$ as a collection $(M,M',\tau, x_i, \cB)$ as above together with a subsheaf $s':\cE'\otimes \Omega\hook{} M'$. Let 
$h^{\la}_Z: \bar Z\to \Bunb_{\tilde B}$ be the map sending the above point to $(M,\cB, \cE, s)$, where $\cE=\cE'(-\sum_i x_i)$ and $s$ is the composition
$$
\cE'\otimes \Omega\hook{s'} M'\toup{\tau^{-1}} M(D_1)
$$
with $D_1=\sum_i x_i$. Let $Z\subset \bar Z$ be the open substack given by the properties: the map $s': \cE'\otimes\Omega\to M'$ has no zero at each $x_i$, and $s: \cE\otimes\Omega\to M$ has no zero at each $x_i$. 

 Let $\phi: Z\to \Bunb_{\tilde B}\times (X^m-\vartriangle)$ be the map, whose first component is $h^{\la}_Z$, and the second sends a point of $Z$ to the collection $(x_i)\in X^m-\vartriangle$. Then $\phi$ is an isomorphism, because $S^{w_0(\lambda)}\cap \Gr^{\lambda}$ is a point scheme in the notation of (\cite{BG}, Section~5.2.2). 
 
 Let $AJ: (X^m-\vartriangle)\to \Bun_T$ be the map sending $\{x_i\}$ to $\cF_T=\cF^0_T(\sum_i \alpha x_i)$. Write $\bar\phi: Z\to \Bunb_{\tilde B}\times\Bun_T$ for $(\id\times AJ)\comp\phi$. Let $m$ be the composition 
$$
\Bunb_{\tilde B}\times \Bun_T\toup{\bar\gq\times\id} \Bun_T\times\Bun_T\to\Bun_T,
$$
where the second map is the tensor product. We get a commutative diagram
$$
\begin{array}{ccc}
Z & \toup{h^{\ra}_Z} & \Bunb_{\tilde B}\\
\downarrow\lefteqn{\scriptstyle \bar\phi} && \downarrow\lefteqn{\scriptstyle\bar \gq}\\
\Bunb_{\tilde B}\times \Bun_T & \toup{m} & \Bun_T,
\end{array}
$$
where $h^{\ra}_Z$ is smooth and surjective, and $\tilde\phi$ is smooth. 
Since 
$$
(h^{\ra}_Z)^*\IC_{B,W}\,\iso\, \bar\phi^*(\IC_{B,W}\boxtimes \IC(\Bun_T))
$$ 
up to a shift, our assertion for $\IC_{B,W}$
follows as in (\cite{BG}, Sections~5.2.1-5.2.2). The proof for $\IC_{B,W!}$ is similar, use the fact that, for a point of $Z$, $\cE\otimes\Omega\subset M$ is a subbundle iff $\cE'\otimes\Omega\subset M'$ is a subbundle.
\end{Prf}

\begin{Cor} 1) The functor $\bar\gq^{!*}$ is exact for the perverse t-structures and commutes with the Verdier duality.\\
2) If $K\in \P(\Bun_n)$ then $\bar\gq^{!*}(K)\,\iso\, (j_P)_{!*} \gq^{!*}(K)$ naturally.
\end{Cor}
\begin{Prf} Identical to \cite{BG}, Theorem~2.1.2 using the above proposition.
\end{Prf}

\medskip

 Now define $\Eis: \D(\Bun_n)\to \D_-(\Bunt_G)$ by $\Eis(K)=\bar\gp_!\bar\gq^{!*}(K)$ for the diagram (\ref{diag_Eis_construction}). If $K\in \D(\Bun_n)$ is supported by a connected component of $\Bun_n$ with $\chi(E\otimes\Omega)$ even (resp., odd) then $\Eis(K)\in \D_{\mp}(\Bunt_G)$ (resp., $\Eis(K)\in \D_{=}(\Bunt_G)$).
 
\begin{Cor} The functor $\Eis$ commutes with the Verdier duality and maps pure complexes to pure ones.
\end{Cor}
 
\medskip\noindent 
4.1.3 For the rest of Section~4 we assume $n=1$. In this case we also write $\Bunb_B:=\Bunt_B$. It is easy to see that $\Bunb_B$ is smooth. In turn, this implies that $\Bunb_{\tilde B}:=\Bunt_{\tilde B}$ is smooth. Note that $\Bunb_{\tilde B}=\cQ_1$, where $\cQ_1$ is defined in Section~2.1. 

 For any $x\in X$ the perverse sheaf $\IC_{B,W}\mid_{\cQ^x_1}$ is equivariant under the action of the groupoid $\cH^x_1$ of Section~2.2. Indeed, this is true over $\cQ_1^x\cap \Bun_{\tilde B}$, and this property is preserved under the intermediate extension. This shows that $\IC_{B,W}$ is an object of the category $\P^{\cL}(\cQ_1)$ introduced in Section~2.5. Remind that $^d\cQ_1$ is defined in Section~2.4, and $W$ denotes the nontrivial rank one and order 2 local system on $B(\mu_2)$.

 Consider the stack $^d\bar\cP=X^{(d)}\times\Bun_1\times B(\mu_2)$ defined in Section~3.2. A point of $^d\cQ_1$ is given by $(M,\cB)\in\Bunt_G$, $\cE\in\Bun_1, D_1\in X^{(d)}$ and a subbundle $\cE\otimes\Omega(D_1)\hook{} M$. 
Let $^d\bar\pr: {^d\cQ_1}\to {^d\bar\cP}$ be the map sending the above point to $(D_1, \cE, \cB_0)$, where 
$$
\cB_0=\cB\otimes\det\RG(X, \cE\otimes\Omega(D_1))^{-1}
$$ as in Section 2.4. From Lemma~\ref{Lm_description_D_cL_d_cQ_2ex}, we learn that the $*$-restriction of $\IC_{B,W}$ to each stratum $^d\cQ_1\subset\cQ_1$ descends with respect to the projection $^d\bar\pr: {^d\cQ_1}\to {^d\bar\cP}$. 
 
  For $d$ even let $\kappa: X^{(d/2)}\to X^{(d)}$ be the map sending $D$ to $2D$. Let 
$$
^d\cQ_{1,even}={^d\cQ_1}\times_{X^{(d)}} X^{(d/2)},
$$
where we used $\kappa$ to define the fibred product. Let $^d\pr_{even}$ be the composition of projections 
$$
{^d\cQ_{1,even}}\toup{^d\bar\pr} X^{(d/2)}\times\Bun_1\times B(\mu_2)\to B(\mu_2)
$$ 
For the inclusion $^d\cQ_1\subset\cQ_1$ we have $\dim(^d\cQ_1)=\dim\cQ_1-d$ with the abuse of notation meaning that the dimension of a connected component of $\cQ_1$ is compared with the dimension of the stratum $^d\cQ_1$ intersected with this component.

\begin{Th} 
\label{Th_3}
For any $d\ge 0$ the $*$-restriction of $\IC_{B,W}$ to the stratum $^d\cQ_1\subset \Bunb_{\tilde B}$ vanishes unless $d$ is even. For $d$ even, it is the extension by zero under $^d\cQ_{1,even}\hook{} {^d\cQ_1}$ of the complex
$(^d\pr_{even})^*W[\dim\cQ_1-d]$. 
\end{Th}
  
\medskip\noindent
4.2.1 Consider the stack $\Bun_B\times_{\Bun_G}\Bunb_B$. Write a point of this stack as an exact sequence (\ref{ext_cE^*_by_EotimesOmega}) together with a subsheaf $\tilde t: \cE_1\otimes\Omega\hook{} M$ with $\cE_1\in\Bun_1$. Let $\cX_B\subset \Bun_B\times_{\Bun_G}\Bunb_B$ be the open substack given by the propety that $\tilde t$ does not factor through $\cE\otimes\Omega$. A point of $\cX_B$ gives rise to the diagram
\begin{equation}
\label{diag_point_of_cX_B}
\begin{array}{cccc}
0\to \cE\otimes\Omega\to M & \!\to &\cE^{-1} & \to 0\\
& \!\nwarrow\lefteqn{\scriptstyle \tilde t} & \uparrow\lefteqn{\scriptstyle t}\\
&& \cE_1\otimes\Omega
\end{array}
\end{equation}
Let $\cX_{\tilde B}\subset \Bun_B\times_{\Bun_G}\Bunb_{\tilde B}$ be defined by the same condition. 

 Consider the stack $\cZ_B$ classifying $\cE,\cE_1\in\Bun_1$ and an inclusion $t: \cE_1\otimes\Omega\hook{} \cE^{-1}$. Write $t^*: \cE\hook{} \cE_1^*\otimes \Omega^*$ for the dual map. Let $K_t$ denote the cokernel of $t^*: \cE^2\otimes\Omega\hook{}\cE\otimes\cE_1^{-1}$. Let $p_{\cX}: \cX_B\to \cZ_B$ be the map sending the above point of $\cX_B$ to $(\cE, \cE_1, t)$. By (\cite{LL}, Section~8.2, Example 1), the map $p_{\cX}$ is a vector bundle whose fibre over a point of $\cZ_B$ is $\H^0(X, K_t)$. 
  
 Let $\nu_B: \Bun_B\to\Bun_G$ be the map sending (\ref{ext_cE^*_by_EotimesOmega}) to $M$. Let $\Bun_B^0\subset \Bun_B$ be the open substack given by $\H^0(X, \cE^2\otimes\Omega^2)=0$. The restriction $\nu_B^0: \Bun_B^0\to\Bun_G$ of $\nu_B$ is smooth. 
 
 Denote by $\cX_B^0\subset \cX_B$ the open substack given by $\H^0(X, \cE^2\otimes\Omega^2)=0$, similarly for $\cX_{\tilde B}^0$. 
% Let $q_{\cZ}: \cZ_B\to \Bun_1$ be the map sending a point $(\cE, \cE_1,t)$ to $\cE$.
We have an isomorphism 
\begin{equation}
\label{trivialization_gerb_over_cX_B} 
 \cX_B\times B(\mu_2) \,\iso\, \cX_{\tilde B}
\end{equation} 
sending $\cB_0$ equipped with $\cB_0^2\,\iso\, k$ and (\ref{diag_point_of_cX_B}) to (\ref{diag_point_of_cX_B}) and $\cB=\cB_0\otimes \det\RG(X, \cE\otimes\Omega)$ equipped with the isomorphism $\cB^2\,\iso\, \det\RG(X,M)$ induced by (\ref{ext_cE^*_by_EotimesOmega}).
 
\medskip\noindent
4.2.2 Let $\theta\ge 0$, we may think of $\theta$ as the element $\theta\alpha$ of the $\ZZ_+$-span of positive coroots. Fix for a moment a line bundle $\cE$ on $X$. Define $Z^{\theta}$ and $\tilde Z^{\theta}$ as 
$$
Z^{\theta}=\cX_B\times_{\cZ_B} X^{(\theta)}\;\;\;\;\mbox{and}\;\;\;\; \tilde Z^{\theta}=\cX_{\tilde B}\times_{\cZ_B} X^{(\theta)},
$$
where the map $X^{(\theta)}\to \cZ_B$ sends $D$ to $(\cE, \cE_1=\cE^{-1}\otimes\Omega^{-1}(-D), t)$, here $t:\cE_1\otimes\Omega\hook{} \cE^{-1}$ is the canonical inclusion. We will write $Z^{\theta}_{\cE}$, $\tilde Z^{\theta}_{\cE}$ if we need to express the dependence on $\cE$.

 The scheme $Z^{\theta}$ is the usual Zastava space for $G$ considered in \cite{BFGM}. Let $\pi_{\theta}: Z^{\theta}\to X^{(\theta)}$ denote the projection, this
is a vector bundle with fibre $\cE^2\otimes\Omega(D)/\cE^2\otimes\Omega$ over $D\in X^{(\theta)}$. The fibre of the dual vector bundle $\check{Z}^{\theta}\to X^{(\theta)}$ over $D$ is $\cE^{-2}/\cE^{-2}(-D)$. The isomorphism (\ref{trivialization_gerb_over_cX_B}) yields a trivialization 
$$
\epsilon^{\theta}: Z^{\theta}\times B(\mu_2)\,\iso\, \tilde Z^{\theta} 
$$
As in \cite{BFGM}, denote by $Z^{\theta}_{max}\subset Z^{\theta}$ the open subscheme given by the property that $\tilde t: \cE_1\otimes\Omega\to M$ is a subbundle. A point 
$$
(D\in X^{(\theta)},\sigma\in \cE^2\otimes\Omega(D)/\cE^2\otimes\Omega)\in Z^{\theta}
$$ 
lies in $Z^{\theta}_{max}$ iff for any effective divisor $D'<D$ we have $\sigma\notin \cE^2\otimes\Omega(D')/\cE^2\otimes\Omega$.

 Over $Z^{\theta}_{max}$ we have another trivialization of the gerb $\tilde Z^{\theta}\to Z^{\theta}$, the difference between the two is given by the following local system $W_{max}$ on $Z^{\theta}_{max}$. For a point of $Z^{\theta}_{max}$ we get canonical $\ZZ/2\ZZ$-graded isomorphisms
$$
\det\RG(X, \cE^{-1}(-D))^2\,\iso\,\det\RG(X, M)\,\iso\,\det\RG(X, \cE\otimes\Omega)
$$
and
$$
\det\RG(X, \cE^{-1}(-D))\otimes\det\H^0(X, \cE^{-1}/\cE^{-1}(-D))\,\iso\, \det\RG(X, \cE^{-1})\,\iso\, \det\RG(X, \cE\otimes\Omega)
$$
They yield a $\ZZ/2\ZZ$-graded trivialization 
$
\det\H^0(X, \cE^{-1}/\cE^{-1}(-D))^2\,\iso\, k
$. Let $Y^{\theta}_{max}\to Z^{\theta}_{max}$ denote the $\mu_2$-torsor whose fibre consists of elements of 
$$
\det\H^0(X, \cE^{-1}/\cE^{-1}(-D))
$$ 
of square one. Define $W_{max}$ as the extension of scalars of this torsor under $\mu_2\hook{}\Qlb^*$. Let $\IC(W_{max})$ be the intermediate extension of $W_{max}[\dim Z^{\theta}]$ to $Z^{\theta}$. 

As in \cite{BFGM}, $\epsilon^{\theta}_*(\IC(W_{max})\boxtimes W)$ is a local model for the perverse sheaf $\IC_{B,W}$. 

  For $\theta_i\ge 0$ with $\theta_1+\theta_2=\theta$ let $X^{\theta_1,\theta_2}_{disj}\subset X^{(\theta_1)}\times X^{(\theta_2)}$ be the open subscheme of divisors $D_1,D_2$ such that $D_1\cap D_2=\emptyset$. The usual factorization property of $Z^{\theta}$ is the natural isomorphism
\begin{equation}
\label{iso_factorization_Z_theta}
Z^{\theta}\times_{X^{(\theta)}} \;X^{\theta_1,\theta_2}_{disj}\,\iso\, (Z^{\theta_1}\times Z^{\theta_2})\times_{X^{(\theta_1)}\times X^{(\theta_2)}} \;X^{\theta_1,\theta_2}_{disj}
\end{equation}
Clearly, $Z^{\theta}_{max}$ satisfies the same factorization property. 
\begin{Lm} 
\label{Lm_first_factorization_property}
1) There is a canonical isomorphism $W_{max}\mid_{X^{\theta_1,\theta_2}_{disj}}\,\iso\, (W_{\max}\boxtimes W_{\max})$ over
$$
Z^{\theta}_{max}\times_{X^{(\theta)}} \;X^{\theta_1,\theta_2}_{disj}\,\iso\, (Z^{\theta_1}_{max}\times Z^{\theta_2}_{max})\times_{X^{(\theta_1)}\times X^{(\theta_2)}} \;X^{\theta_1,\theta_2}_{disj}
$$
2) Over (\ref{iso_factorization_Z_theta}) we have canonically $\IC(W_{max})\,\iso\, \IC(W_{max})\boxtimes \IC(W_{max})$. 
\end{Lm}
\begin{Prf} The $\mu_2$-torsor $Y^{\theta}_{max}$ after the base change by $X^{\theta_1,\theta_2}_{disj}$ is obtained from $Y^{\theta_1}_{max}\times Y^{\theta_2}_{max}$ by the extension of scalars $\mu_2\times\mu_2\to \mu_2$ given by the product in $\mu_2$.
\end{Prf}

\medskip

 In the case $\theta=1$ the fibre of $Z^{\theta}\to X^{(\theta)}=X$ over $x$ is the geometric fibre $\cE^2_x$. In this case let $Y\to X$ be the total space of the vector bundle $\cE$ on $X$. Let $\nu_Y: Y\to Z^{\theta}$ be the map over $X$ sending $v$ to $v^2$. Then $\IC(W_{max})$ is this case identifies with $\mu_2$-antiinvariants in $\nu_{Y !}\Qlb$. In particular, the $*$-restriction of $\IC(W_{max})$ under the zero section $X\to Z^{\theta}$ vanishes.
 
  For $D=\theta x$ the fibre $\SSS^{\theta}_x$ of $Z^{\theta}_{max}$ over $D$ consists of elements $\sigma\in \cE^2\otimes\Omega(\theta x)/\cE^2\otimes\Omega$ that do not lie in 
$\cE^2\otimes\Omega((\theta-1) x)/\cE^2\otimes\Omega$. There is a canonical $\ZZ/2\ZZ$-graded isomorphism
$$
\det\H^0(X, \cE^{-1}/\cE^{-1}(-\theta x))\,\iso\, \cE_x^{-\theta}\otimes \cO({\scriptstyle\frac{(1-\theta)\theta}{2}})_x
$$ 
The geometric fibre of $\sigma$ is a nonzero element $\sigma_x\in \cE^2\otimes\Omega(\theta x)_x\,\iso\, \cE^2_x\otimes \cO((\theta-1)x)_x$. The covering $Y^{\theta}_{max}\to Z^{\theta}_{max}$ over $\SSS^{\theta}_x$ consists of elements 
$$
\gamma\in  \cE_x^{-\theta}\otimes \cO({\scriptstyle\frac{(1-\theta)\theta}{2}})_x
$$ 
such that $\gamma^2=\sigma_x^{-\theta}$. We see that the $*$-restriction of $W_{max}$ to $\SSS^{\theta}_x$ is trivial (resp., nontrivial) for $\theta$ even (resp., for $\theta$ odd). 

 One can describe explicitely the restriction of $W_{max}$ under the projection 
$Z^{\theta}_{max}\times_{X^{(\theta)}} {^{rss}X^{(\theta)}}\to Z^{\theta}_{max}$. 
Given $D=\sum_i x_i\in {^{rss}X^{(\theta)}}$, we have a canonical $\ZZ/2\ZZ$-graded isomorphism 
$$
\det\H^0(X, \cE^{-1}/\cE^{-1}(-D))\,\iso\, \otimes_i \cE^{-1}_{x_i}
$$ 
The fibre of the $\mu_2$-torsor 
$$
Y^{\theta}_{max}\times_{X^{(\theta)}} {^{rss}X^{(\theta)}}\to Z^{\theta}_{max}\times_{X^{(\theta)}} {^{rss}X^{(\theta)}}
$$
over a given $D=\sum_i x_i$, $\sigma=(\sigma_i)\in \oplus_i \cE^2_{x_i}$ consists of those elements $\gamma\in \otimes_i \cE_{x_i}$ such that $\gamma^{-2}=\otimes \sigma_i\in \otimes_i \cE^2_{x_i}$. Here the tensor product is taken in the $\ZZ/2\ZZ$-graded sense, so is independent of the order of $x_i$.
We have a group scheme $\GG^{\theta}$ over $X^{(\theta)}$ with fibre $(\cO/\cO(-D))^*$ over $D$.

\begin{Pp}  
\label{Pp_monodromy_important}
Consider the $S_{\theta}$-Galois covering
$$
Z^{\theta}_{max}\times_{X^{(\theta)}} {^{rss} X^{\theta}}\to Z^{\theta}_{max}\times_{X^{(\theta)}} {^{rss} X^{(\theta)}}
$$ 
By the factorization property of Lemma~\ref{Lm_first_factorization_property}, the restriction of $W_{max}$ under this map identifies with $\boxtimes_{i=1}^{\theta} W_{max}$ over $(\prod_{i=1}^{\theta} Z^1_{max})\times _{X^{\theta}}{^{rss}X^{\theta}}$. In addition, the corresponding descent data with respect to $S_{\theta}$ is trivial, that is, there are no $S_{\theta}$-monodromy in the latter isomorphism.
\end{Pp}
  
\begin{Lm} 
\label{Lm_square_over_line_bundle}
For any $D\in X^{(\theta)}$ there is a canonical $\ZZ/2\ZZ$-graded isomorphism 
$$
\det\RG(X, \cE^2\otimes\Omega(D)/\cE^2\otimes\Omega)\otimes\det\RG(X, \Omega(D)/\Omega)\,\iso\, \det\RG(X, \cE^{-1}/\cE^{-1}(-D))^2
$$
It can be seen as an isomorphism of the corresponding $\ZZ/2\ZZ$-graded line bundles on $X^{(\theta)}$. 
\end{Lm}
\begin{Prf}
Apply \cite{L3}, Lemma~1, i) and ii).
\end{Prf} 
  
\medskip  
\begin{Prf}\select{of Proposition~\ref{Pp_monodromy_important}} 

\medskip\noindent
Consider the line bundle over $^{rss}X^{(\theta)}$ whose fibre over $D=\sum_i x_i$ is $\det\RG(X, \Omega(D)/\Omega)$. This is the line bundle $\cO$ over $^{rss}X^{\theta}$ equipped with the descent data with respect to $^{rss}X^{\theta}\to {^{rss}X^{(\theta)}}$ 
given by the sign representation of $S_{\theta}$.

 Consider the line bundle on $^{rss}X^{(\theta)}$ with fibre $\det\RG(X, \cE^2\otimes\Omega(D)/\cE^2\otimes\Omega)$. This is the line bundle $\boxtimes_{i=1}^{\theta} \cE^2$ on $^{rss}X^{\theta}$ with the descent data given by the sing representation of $S_{\theta}$. 
 
  Now consider the line bundle $(\sym_*(\boxtimes_{i=1}^{\theta} \cE^2))^{S_{\theta}}$ on $^{rss}X^{(\theta)}$, where $\sym: X^{\theta}\to X^{(\theta)}$ is the sum of divisors. Let $\cV^{\theta}$ denote the total space of this line bundle with zero section removed. By Lemma~\ref{Lm_square_over_line_bundle}, we get a $\mu_2$-local system $W\cV^{\theta}$ over $\cV^{\theta}$ classifying $D\in {^{rss} X^{(\theta)}}, \sigma\in \cV^{\theta}$ and $\gamma\in\det\RG(X, \cE^{-1}/\cE^{-1}(-D))$ with $\gamma^2=\sigma$.
    
  There is a morphism 
$$
\delta: Z^{\theta}_{max}\times_{X^{(\theta)}} {^{rss} X^{(\theta)}}\to \cV^{\theta}
$$
over $^{rss} X^{(\theta)}$ sending $\sigma=(\sigma_i)\in \oplus_i \cE^2_{x_i}$ to 
$\otimes_i \sigma_i$. Then $\delta^*(W\cV^{\theta})$ identifies with the restriction of $W_{max}$ to $Z^{\theta}_{max}\times_{X^{(\theta)}}$.

 We have the natural map $\xi: (\prod_{i=1}^{\theta} \cV^1)\times_{X^{\theta}}{^{rss}X^{\theta}}\to \cV^{\theta}$ given by the product, and $\xi^*W\cV^{\theta}\,\iso\, \boxtimes_{i=1}^{\theta} (W\cV^1)$. The $S_{\theta}$-minodromy in the latter isomorphism is trivial. Our assertion follows now from the commutative diagram
$$ 
\begin{array}{ccc}
Z^{\theta}_{max}\times_{X^{(\theta)}} {^{rss} X^{\theta}} & \iso & (\prod_{i=1}^{\theta} Z^1_{max})\times_{X^{\theta}} {^{rss} X^{\theta}}\\
\downarrow\lefteqn{\scriptstyle \delta\times\id} && \downarrow\\
\cV^{\theta}\times_{X^{(\theta)}} {^{rss} X^{\theta}} & \gets &
(\prod_{i=1}^{\theta} \cV^1)\times_{X^{\theta}}{^{rss}X^{\theta}}
\end{array}
$$ 
\end{Prf}  

\medskip\noindent
4.2.3 %This subsection is not used for the proof of our main results and may be skipped. 
For any $\theta\ge 0$ let $\check{Y}^{\theta}\to X^{(\theta)}$ be the vector bundle with fibre $\cE^{-1}/\cE^{-1}(-D)$ over $D$. Let $\pi_{\check{Y}}: \check{Y}^{\theta}\to \check{Z}^{\theta}$ be the map over $X^{(\theta)}$ sending a section $v\in \cE^{-1}/\cE^{-1}(-D)$ to $v^2\in \cE^{-2}/\cE^{-2}(-D)$. 
 
 The morphism $\pi_{\check{Y}}$ is equivariant with respect to the group scheme $\GG^{\theta}$. It is understood that $g\in (\cO/\cO(-D))^*$ acts on $\cE^{-1}/\cE^{-1}(-D)$ as $g$, and on $\cE^{-2}/\cE^{-2}(-D)$ as $g^2$. The morphism $\pi_{\check{Y}}$ is affine, so $\pi_{\check{Y} !}\IC(\check{Y}^{\theta})$ is placed in perverse degrees $\ge 0$ and is $\GG^{\theta}$-equivariant.
 
\begin{Pp} 
\label{Pp_small_map}
The morphism $\pi_{\check{Y}}: \check{Y}^{\theta}\to \check{Z}^{\theta}$ is small, and $\pi_{\check{Y} !}\IC(\check{Y}^{\theta})$ is a perverse sheaf. The $*$-restriction of $\pi_{\check{Y} !}\IC(\check{Y}^{\theta})$ to the complement of $\check{Z}^{\theta}\times_{X^{(\theta)}} {^{rss}X^{(\theta)}}$ is placed in perverse degrees $<0$. 
\end{Pp}
\begin{Prf} 
For a partition $U(\theta)=(n_1\ge\ldots\ge n_k\ge 1)$ of $\theta$ consider the 
scheme $X^{U(\theta)}={^{rss}X^k}$ and the map
map $X^{U(\theta)}\to X^{(\theta)}$ sending $\{x_1,\ldots, x_k\}$ to 
$D=\sum_i n_i x_i$. Here $^{rss}X^k\subset X^k$ is the complement to all diagonals. Let $\gB(\theta)$ be a datum of $U(\theta)$ together with a collection 
$(m_1, \ldots, m_k)$, where $0\le m_i\le n_i$ for each $i$. 

 Given $\gB(\theta)$, let $\check{Z}^{\gB(\theta)}$ be the locally closed subscheme in $\check{Z}^{\theta}\times_{X^{(\theta)}} X^{U(\theta)}$ classifying $D\in X^{U(\theta)}$ and $v=(v_i)\in \prod_i \cE^{-2}/\cE^{-2}(-n_ix_i)$ such that $$
v_i\in \cE^{-2}(-m_ix_i)/\cE^{-2}(-n_ix_i)\;\;\;\;\mbox{and}\;\;\;\;
v_i\notin \cE^{-2}(-(m_i+1)x_i)/\cE^{-2}(-n_ix_i)
$$
The dimension of $\check{Z}^{\gB(\theta)}$ is $k+\sum_i (n_i-m_i)$ and $\dim Z^{\theta}=2\theta$. The fibre of 
$$
\pi_{\check{Y}}\times\id: \check{Y}^{\theta}\times_{X^{(\theta)}} X^{U(\theta)}\to \check{Z}^{\theta}\times_{X^{(\theta)}} X^{U(\theta)}
$$ 
over a point of $\check{Z}^{\gB(\theta)}$ is of the form $\prod_i S_i$. If $m_i=n_i$ then $S_i$ identifies with $\AA^{n_i/2}$ (resp., with $\AA^{(n_i-1)/2}$) for $n_i$ even (resp., for $n_i$ odd). If $0\le m_i<n_i$ then $S_i$ is
empty unless $m_i$ is even, and for $m_i$ even $S_i\,\iso\, \AA^{m_i/2}\sqcup \AA^{m_i/2}$.

 If the fibre of $\pi_{\check{Y}}\times\id$ over a point of $\check{Z}^{\gB(\theta)}$ is nonempty then the dimension of the fibre is $\le \sum_i \frac{m_i}{2}$. So, $2\dim(fibre)=\sum_i m_i\le \codim(\check{Z}^{\gB(\theta)})=\theta-k+\sum_i m_i$, and the inequality is strict unless $k=\theta$. 
\end{Prf} 

\medskip

 Note that $\pi_{\check{Y} !}\IC(\check{Y}^{\theta})$ is acted on by $\mu_2$. It surjects onto the intermediate extension of its restriction to $\check{Z}^{\theta}\times_{X^{(\theta)}} {^{rss}X^{(\theta)}}$, because it does not admit any perverse quotient sheaves supported over the complement of $\check{Z}^{\theta}\times_{X^{(\theta)}} {^{rss}X^{(\theta)}}$. 
 
\medskip\noindent
4.2.4 Consider the homomorphism $\zeta: \GG^{\theta}\to\Gm$
of group schemes over $X^{(\theta)}$ sending $g\in (\cO/\cO(-D))^*$ to the induced automorphism of $\det\H^0(X, \cO/\cO(-D))$. We have the local system over $\Gm$ corresponding to the covering $\Gm\to \Gm$, $g\mapsto g^2$ and the identity character $\mu_2\hook{}\Qlb^*$. Let $W_{\GG}$ denote the restriction of this local system under $\zeta$.
Then $W_{\GG}$ is a character local system on $\GG^{\theta}$, its $*$-restriction under the unit section $X^{(\theta)}\to \GG^{\theta}$ is naturally trivialized. 
 
 Consider the action of $\GG^{\theta}$ on $Z^{\theta}$ such that $g\in \GG^{\theta}$ acts on $v\in \cE^2\otimes\Omega(D)/\cE^2\otimes\Omega$ as $gv$. The perverse sheaf $\IC(W_{max})$ is $(\GG^{\theta}, W_{\GG})$-equivariant. Indeed, this holds over the locus $Z^{\theta}\times_{X^{(\theta)}}{^{rss}X^{(\theta)}}$, and this property is preserved under the intermediate extension.
 
  For $D\in X^{(\theta)}$ let $\GG^{\theta}_D$ be the fibre of $\GG^{\theta}$ over $D$.  If $D=\theta x\in X^{(\theta)}$ for some $x\in X$ then the $*$-restriction of $W_{\GG}$ to $\GG^{\theta}_D$ is trivial iff $\theta$ is even. 
 
% The $*$-restriction of $\IC(W_{max})$ to $Z^{\gB(\theta)}$ is $(\GG^{\theta}, W_{\GG})$-equivariant and so, by Remark~\ref{Rem_one} ii), is as an inverse image of some sheaf $K^{\cB(\theta)}$ from $X^{U(\theta)}$ tensored by some fixed $(\GG^{\theta}, W_{\GG})$-equivariant local system on $Z^{\gB(\theta)}$. Moreover $K^{\gB(\theta)}$ satisfies the factorization property. So, to understand
% $\IC_{B,W}$ it suffices to calculate $K^{\cB(\theta)}$ over the main diagonal, that is, for the partition $U(\theta)=(\theta)$.   

 For $\theta'\le \theta$ let
$$
_{\theta'}Z^{\theta}=Z^{\theta-\theta'}_{max}\times X^{(\theta')}
$$ 
As in \cite{BFGM}, this is a locally closed subscheme in $Z^{\theta}$ classifying diagrams (\ref{diag_point_of_cX_B}) and a divisor $D'\in X^{(\theta')}$ such that $\tilde t: \cE_1\otimes\Omega(D')\hook{} M$ is a subbundle, and $\div(\cE^{-1}/\cE_1\otimes\Omega(D'))\in X^{(\theta-\theta')}$. 

  Let $\gs^{\theta}: X^{(\theta)}\to Z^{\theta}$ denote the zero section, it identifies with $_{\theta}Z^{\theta}$. From the $\Gm$-equivariance, where $g\in \Gm\subset \GG^{\theta}$ acts as $g^2$, as
in (\cite{BFGM}, Proposition~5.2) we get an isomorphism 
$$
\pi_{\theta !} \IC(W_{max})\,\iso\, \gs^{\theta !}\IC(W_{max}),
$$ 
By purity as in (\cite{BFGM}, Corollary~5.5), the latter complex is a direct sum of irreducible perverse sheaves. 

\begin{Pp} 
\label{Pp_restriction_zero_section}
The complex $\pi_{\theta !} \IC(W_{max})$ vanishes if and only if $\theta$ is odd. For $\theta$ even there is an isomorphism over $X^{(\theta)}$
\begin{equation}
\label{iso_for_Pp_restriction_zero_section}
\pi_{\theta !} \IC(W_{max})\,\iso\, \kappa_!\Qlb,
\end{equation}
where $\kappa: X^{(\theta/2)}\to X^{(\theta)}$ sends $D$ to $2D$. 
\end{Pp}

 Proposition~\ref{Pp_restriction_zero_section} implies an isomorphism $\gs^{\theta *}\IC(W_{max})\,\iso\, \kappa_!\Qlb[\theta]$. Moreover, for $\theta'\le \theta$ it yields an isomorphism for the $*$-restriction
\begin{equation}
\label{iso_IC_Wmax_on_stratum}
\IC(W_{max})\mid_{_{\theta'}Z^{\theta}}\,\iso\, \left\{
\begin{array}{ll}
\pr_1^*W_{max}\boxtimes \kappa_!\Qlb[2\theta-\theta'], & \mbox{for}\;\; \theta'\;\; \mbox{even}\\
0, &  \mbox{for}\;\; \theta'\;\; \mbox{odd},
\end{array}
\right.
\end{equation}
here now $\kappa: X^{(\theta'/2)}\to X^{(\theta')}$ for $\theta'$ even.

\medskip
\begin{Prf}\select{of Proposition~\ref{Pp_restriction_zero_section}} \  
Since $\IC(W_{max})$ is $(\GG^{\theta}, W_{\GG})$-equivariant, the complex $\gs^{\theta !}\IC(W_{max})$ is also $(\GG^{\theta}, W_{\GG})$-equivariant by functoriality. However, the action of $\GG^{\theta}$ on the zero section is trivial.
So, if $\theta$ is odd then the $*$-restriction of $\pi_{\theta !}\IC(W_{max})$ to the main dagonal $\delta: X\hook{} X^{(\theta)}$ will vanish. Indeed, for any $D\in X^{(\theta)}$ the $*$-restriction of $W_{\GG}$ to $\GG^{\theta}_D$ will be nontrivial.

 Assume $\theta>0$ even and our result known for $\theta'<\theta$. We argue by induction on $\theta$. From the factorization property of $\IC(W_{max})$ and $(\GG^{\theta}, W_{\GG})$-equivariance, we also see that $\pi_{\theta !} \IC(W_{max})$ vanishes outside the image of $\kappa$. 
 
 Given a decomposition $\theta=\theta_1+\theta_2$ with $\theta_i>0$ even, after the base change $X^{\theta_1, \theta_2}_{disj}\to X^{(\theta)}$ we get an isomorphism by the induction hypothesis and 2) of Lemma~\ref{Lm_first_factorization_property}
\begin{equation}
\label{iso_outside_diag_for_Pp}
\pi_{\theta !}\IC(W_{max})\mid_{X^{\theta_1, \theta_2}_{disj}}\,\iso\, (\kappa\times\kappa)_!\Qlb\mid_{X^{\theta_1, \theta_2}_{disj}},
\end{equation}
where $\kappa\times\kappa: X^{(\theta_1/2)}\times X^{(\theta_2/2)}\to X^{(\theta_1)}\times X^{(\theta_2)}$. By Proposition~\ref{Pp_monodromy_important}, the isomorphism (\ref{iso_outside_diag_for_Pp}) has no $S_2$-monodromy with respect to the permutation of two divisors in the case $\theta_1=\theta_2$. 
  
  By purity, for $\theta>2$ there is a pure complex $K^{\theta}$ on $X$ such that for the main diagonal $\vartriangle: X\to X^{(\theta)}$ we have
\begin{equation}
\label{iso_for_theta_at_least_4}
\pi_{\theta !} \IC(W_{max})\,\iso\, \kappa_!\Qlb\oplus \vartriangle_!\! K^{\theta}
\end{equation}   
For $\theta=2$ we don't know yet that $\kappa_!\Qlb$ is a direct summand in $\pi_{\theta !} \IC(W_{max})$, as (\ref{iso_outside_diag_for_Pp}) is vacuous in this case.
 
 First, we show that $\vartriangle^*\!\pi_{\theta !} \IC(W_{max})$ is placed in perverse degrees $\le 1$, and $^p\H^1$ of this complex identifies with $\Qlb$. To do so, calculate this complex using the stratification of $Z^{\theta}$ by $_{\theta'}Z^{\theta}$ for $0\le \theta'\le \theta$. Let $i_{\theta'}$ denote the composition
$$
_{\theta'}Z^{\theta}\times_{X^{(\theta)}} X \hook{} Z^{\theta}\times_{X^{(\theta)}} X\hook{} Z^{\theta}
$$
The complex $\pi_{\theta !}i_{\theta'}^*\IC(W_{max})$ is calculated using (\ref{iso_IC_Wmax_on_stratum}) by the induction hypothesis. It vanishes for $\theta'$ odd. For $\theta'<\theta$ even it identifies with the direct image with compact support of $\Qlb[2\theta-\theta']$ under the projection $Z^{\theta-\theta'}_{max}\times_{X^{(\theta-\theta')}} X\to X$. The latter is a constant complex placed in usual degrees $\le -\theta'$, so in perverse degrees $\le 0$ for $\theta'\ne 0$. The complex $i_{\theta}^*\IC(W_{max})$ by definition is placed in perverse degrees $<0$. So, $\vartriangle^*\!\pi_{\theta !} \IC(W_{max})$ is placed in perverse degrees $\le 1$, and only $_0Z^{\theta}$ contributes to its 1st perverse cohomology sheaf. We get an isomorphism 
$$
^p\H^1(\vartriangle^*\!\pi_{\theta !} \IC(W_{max}))\,\iso\, \Qlb
$$ 

 For $\theta=2$ the complex $\pi_{\theta !} \IC(W_{max})$ is the extension by zero under $\kappa: X\to X^{(2)}$, it is placed in perverse degrees $\le 1$ and 
$$
^p\H^1(\pi_{\theta !} \IC(W_{max}))\,\iso\, \kappa_!\Qlb
$$ 
Since $\gs^!\IC(W_{max})$ is placed in perverse degrees $>0$, the isomorphism (\ref{iso_for_Pp_restriction_zero_section}) follows for $\theta=2$.

 For $\theta>2$ the isomorphism (\ref{iso_for_theta_at_least_4}) yields
$\vartriangle^*\!\pi_{\theta !} \IC(W_{max})\,\iso\, \Qlb\oplus K^{\theta}$. We see that $K^{\theta}$ is placed in perverse degrees $\le 0$. On the other side, $\gs^!\IC(W_{max})\,\iso\, \kappa_!\Qlb\oplus\!\vartriangle_! \!K^{\theta}$ is placed in perverse degrees $>0$, so $K^{\theta}=0$. We are done.
\end{Prf}
  
\medskip\noindent
Theorem~\ref{Th_3} is a straightforward consequence of Proposition~\ref{Pp_restriction_zero_section}.

\medskip\noindent
4.2.5 Let $\Four_{\psi}(\IC(W_{max}))\in \P(\check{Z}^{\theta})$ denote the Fourier transform of $\IC(W_{max})$. Since $\IC(W_{max})$ is $(\GG^{\theta}, W_{\GG})$-equivariant, $\Four_{\psi}(\IC(W_{max}))$ is also $(\GG^{\theta}, W_{\GG})$-equivariant, where now $g\in G^{\theta}$ acts on $v\in \cE^{-2}/\cE^{-2}(-D)$ as $gv$. 

  Let $\GG^{\theta}_2$ be the kernel of $\GG^{\theta}\to\GG^{\theta}$, $g\mapsto g^2$. Remind the map $\pi_{\check{Y}}$ defined in Section~4.2.3. Now $\GG^{\theta}_2$ is a group scheme over $X^{(\theta)}$, it acts naturally on $\pi_{\check{Y} !}\IC(\check{Y}^{\theta})$ acting trivially on $\check{Z}^{\theta}$. 

  For a line bundle $\cL$ on $X$ and an effective divisor $D$ on $X$ write $(\cL/\cL(-D))_{max}$ for the open subscheme of $\cL/\cL(-D)$ classifying $v\in \cL/\cL(-D)$ such that for any $0<D_1\le D$ one has $v\notin \cL(-D_1)/\cL(-D)$.
Denote by 
$$
\check{Z}^{\theta}_{max}\subset \check{Z}^{\theta}
$$ 
the open subscheme classifying $D\in X^{(\theta)}, v\in (\cE^{-2}/\cE^{-2}(-D))_{max}$.
Let $\check{Y}^{\theta}_{max}=\pi^{-1}_{\check{Y}}(\check{Z}^{\theta}_{max})$. The morphism 
\begin{equation}
\label{torsor_for_checkZ_max}
\pi_{\check{Y}}: \check{Y}^{\theta}_{max}\to \check{Z}^{\theta}_{max}
\end{equation}
is finite, this is a torsor under $\GG^{\theta}_2$. 
%Note that 
%\begin{equation}
%\label{Galois_covering_over_checkZ_theta}
%\pi_{\check{Y}}: \check{Y}^{\theta}_{max}\times_{X^{(\theta)}} {^{rss}X^{(\theta)}} \to \check{Z}^{\theta}_{max}\times_{X^{(\theta)}} {^{rss}X^{(\theta)}}
%\end{equation} 
%is a Galois covering with Galois group scheme $\GG^{\theta}_2$.
 Let $\zeta_2$ denote the composition $\GG^{\theta}_2\hook{} \GG^{\theta}\toup{\zeta}\Gm$, it takes values in $\mu_2$ actually. Let $\check{W}_{max}$ be the 
local system on $\check{Z}^{\theta}_{max}$ obtained from (\ref{torsor_for_checkZ_max}) by the extension of scalars via $\zeta_2$. Though the group scheme $\GG^{\theta}_2$ is not flat over $X^{(\theta)}$, the extension of scalars makes sense. 

\begin{Pp} 
\label{Pp_calculation_Four_IC(W_max)}
There is a 1-dimensional vector space $A_{\cE}$ depending on $\cE\in\Bun_1$ and a canonical isomorphism over $\check{Z}^{\theta}_{\cE}$
\begin{equation}
\label{iso_Four_IC_Wmax_answer}
A_{\cE}\otimes \IC(\check{W}_{max})\,\iso\,\Four_{\psi}(\IC(W_{max}))
\end{equation} 
Besides, $A_{\cE}^2\,\iso\, \Qlb$ canonically.
\end{Pp} 
\begin{Prf}
Let $E$ be the local system on $\Gm$ corresponding to the covering $\Gm\to\Gm$, $x\to x^2$ and the identity character $\mu_2\hook{}\Qlb^*$. If $\tilde E$ is the intermediate extension of $E[1]$ under $\Gm\hook{}\A^1$ then $\Four_{\psi}(\tilde E)\,\iso\,\tilde E$, see (\cite{L1}, Propositions~1 and 3). 

  From this observation and the factorization property of $W_{max}$ we obtain the isomorphism (\ref{iso_Four_IC_Wmax_answer}) after restriction to $\check{Z}^{\theta}\times_{X^{(\theta)}}{^{rss}X^{\theta}}$. The fact that the $S_{\theta}$-monodromy for $^{rss}X^{\theta}\to {^{rss}X^{(\theta)}}$ after the Fourier transform is still trivial follows from Lemma~\ref{Lm_Fourier_and_S2-monodromy} below.
\end{Prf}

\medskip
\noindent
4.2.6 The group $S_2$ acts on $\A^2$ by permuting the two coordinates. Write $\A^2/S_2$ for the stack quotient for this action. The Fourier transform can be seen as an auto-equivalence $\Four_{\psi}: \D(\A^2/S_2)\,\iso\, \D(\A^2/S_2)$. 
Given $K\in \P(\A^1)$, consider $K\boxtimes K$ on $\A^2$ with the trivial equivariant structure for the above action of $S_2$, the resulting perverse sheaf on $\A^2/S_2$ will also be denoted $K\boxtimes K$ by abuse of notation.

\begin{Lm}
\label{Lm_Fourier_and_S2-monodromy}
For any perverse sheaf $K\in\A^1$ the functor $\Four_{\psi}: \D(\A^2/S_2)\,\iso\, \D(\A^2/S_2)$ sends $K\boxtimes K$ to $\Four_{\psi}(K)\boxtimes\Four_{\psi}(K)$. That is, the Fourier transform preserves the trivial $S_2$-monodromy for the projection $\A^2\to \A^2/S_2$. 
\end{Lm}
\begin{Prf}
We have an isomorphism over $\A^2/S_2$
$$
\Four_{\psi}(K\boxtimes K)\,\iso\, (\Four_{\psi}(K)\boxtimes\Four_{\psi}(K))\otimes W',
$$ 
where $W'$ is a rank one local system on $B(S_2)$ that we have to determine. Clearly, $W'$ is independent of $K$. We must prove that $W'$ is trivial. Take $K=\Qlb[1]$. Let $p: \A^2/S_2\to B(S_2)$ be the projection and $^tp: B(S_2)\to \A^2/S_2$ the zero section. Since the Fourier transform exchanges $p_!$ and $(^tp)^*$, for this $K$ we get $W'\,\iso\,\Qlb$.
\end{Prf}

\medskip\noindent
4.2.7 Remind that $\GG^{\theta}_2$ acts on $\pi_{\check{Y}!}\IC(\check{Y}^{\theta})$ naturally. In this subsection we prove the following.

\begin{Pp} 
\label{Pp_zeta_2_summand_great}
The isotypic direct summand of $\pi_{\check{Y}!}\IC(\check{Y}^{\theta})$ on which $\GG^{\theta}_2$ acts by $\zeta_2$ identifies canonically with $\IC(\check{W}_{max})$.
\end{Pp}

First, we check that the perverse sheaf $\IC(\check{W}_{max})$ is essentially "the same" as $\IC(W_{max})$. 

 For $0\le \theta'\le\theta$ let $_{\theta'}\check{Z}^{\theta}\subset \check{Z}^{\theta}$ be the locally closed subscheme classifying $D'\in X^{(\theta')}$, $D''\in X^{(\theta-\theta')}$ and $v \in (\cE^{-2}(-D')/\cE^{-2}(-D'-D''))_{max}$. Note that $_0\check{Z}^{\theta}=\check{Z}^{\theta}_{max}$. Let 
$$
_{\theta'}\check{Z}^{\theta}_{even}\subset {_{\theta'}\check{Z}^{\theta}}
$$ 
be the closed subscheme given by the condition that $D'$ is divisible by 2. 
Let $_{\theta'}\check{Y}^{\theta}$ be the scheme classifying $D'\in X^{(\theta')}$ divisible by 2, $D''\in X^{(\theta-\theta')}$ and 
$
v \in \cE^{-1}(-D'/2)/\cE^{-1}(-D'/2-D'')_{max}
$.
Let 
\begin{equation}
\label{map_theta'_pi_checkY}
_{\theta'}\pi_{\check{Y}}: {_{\theta'}\check{Y}^{\theta}}\to {_{\theta'}\check{Z}^{\theta}_{even}}
\end{equation}
be the map sending $(D', D'', v)$ as above to $(D', D'', v^2)$. As above, we have the group schemes $\GG^{\theta-\theta'}_2\subset
\GG^{\theta-\theta'}$ over $X^{(\theta-\theta')}$ and the homomorphism $\zeta_2: \GG^{\theta-\theta'}_2\to \mu_2$. The map (\ref{map_theta'_pi_checkY}) is a torsor under $\GG^{\theta-\theta'}_2$ (viewed this time as a group scheme over ${_{\theta'}\check{Z}^{\theta}_{even}}$). 

 Denote by $_{\theta'}\check{W}_{max}$ the local system on ${_{\theta'}\check{Z}^{\theta}_{even}}$ obtained from (\ref{map_theta'_pi_checkY}) as the extension of scalars via $\zeta_2: \GG^{\theta-\theta'}_2\to \mu_2\hook{}\Qlb^*$. Note that if $\theta'=\theta$ is even then $_{\theta}\check{W}_{max}=\Qlb$. 
 
  By its definition, $\IC(\check{W}_{max})$ satisfies the factorization property as in Lemma~\ref{Lm_first_factorization_property}. So, we can repeat the proof of Proposition~\ref{Pp_restriction_zero_section} for $\check{Z}^{\theta}$ instead of $Z^{\theta}$ and obtain the following.
  
\begin{Pp} 
\label{Pp_IC_checkW_description}
Let $0\le \theta'\le\theta$. The $*$-restriction of $\IC(\check{W}_{max})$ to $_{\theta'}\check{Z}^{\theta}$ vanishes unless $\theta'$ is even. For $\theta'$ even it identifies canonically with the extension by zero under $_{\theta'}\check{Z}^{\theta}_{even}\hook{} {_{\theta'}\check{Z}^{\theta}}$ of $_{\theta'}\check{W}_{max}[2\theta-\theta']$. $\square$
\end{Pp}

\begin{Prf}\select{of Proposition~\ref{Pp_zeta_2_summand_great}.}
Over $\check{Z}^{\theta}_{max}$, the perverse sheaf $\IC(\check{W}_{max})$ coincides by definition with $(\GG^{\theta}_2, \zeta_2)$-invariants in $\pi_{\check{Y}!}\IC(\check{Y}^{\theta})$. By Proposition~\ref{Pp_small_map}, $\pi_{\check{Y}!}\IC(\check{Y}^{\theta})$ surjects onto the intermediate extension of its restriction to $\check{Z}^{\theta}_{X^{(\theta)}} {^{rss}X^{(\theta)}}$, hence a surjective morphism of perverse sheaves 
\begin{equation}
\label{map_surjection_proving_iso}
(\pi_{\check{Y}!}\IC(\check{Y}^{\theta}))^{(\GG^{\theta}_2, \zeta_2)}\to \IC(\check{W}_{max})
\end{equation} 
We check this is an isomorphism over each geometric point. Use notations as in Proposition~\ref{Pp_small_map}, let $\gB(\theta)$ be a datum of a partition $U(\theta)=(n_1\ge\ldots\ge n_k\ge 1)$ of $\theta$ together with a collection $(m_1,\ldots, m_k)$, where $0\le m_i\le n_i$ for all $i$. Let $D=\sum_i n_i x_i$ with $x_i\in X$ pairwise different. Set $D'=\sum_i m_ix_i$, $\theta'=\deg D'$ and $D''=D-D'$.  Let 
$$
v=(v_i)\in (\cE^{-2}(-D')/\cE^{-2}(D))_{max},
$$
here $v_i\in \cE^{-2}(-m_ix_i)/\cE^{-2}(-n_i x_i)$. So, $(D', D'', v)\in {_{\theta'}\check{Z}^{\theta}}$.

 The fibre of $\pi_{\check{Y}}$ over this point is of the form $\prod_i S_i$, where $S_i$ is the scheme classifying $u_i\in \cE^{-1}/\cE^{-1}(-n_i x_i)$ with $u_i^2=v_i$ in $\cE^{-2}/\cE^{-2}(-n_i x_i)$. 
 
 If $m_i=n_i$ is odd then $D'$ is not divisible by 2, and $\mu_2(\cO/\cO(-n_ix_i))$ acts nontrivially on $\det\H^0(X, \cO/\cO(-n_ix_i))$, so the contribution of $S_i$ to the LHS of (\ref{map_surjection_proving_iso}) is zero. 
 
 If $m_i=n_i$ is even then $\mu_2(\cO/\cO(-n_ix_i))$ acts trivially on $\det\H^0(X, \cO/\cO(-n_ix_i))$, and the contribution of $S_i$ to the LHS of (\ref{map_surjection_proving_iso}) is $\Qlb[-m_i]$.
 
 If $m_i<n_i$ with $m_i$ odd then $S_i$ is empty and both sides of (\ref{map_surjection_proving_iso}) at this point are zero. So, we may assume all $m_i$ even. Then the $*$-restriction of $_{\theta'}\check{W}_{max}$ to $\prod_i (\cE^{-2}(-m_ix_i)/\cE^{-2}(-n_ix_i))_{max}$ is the exteriour product $\boxtimes_i E_i$, where $E_i$ is a rank one local system on $(\cE^{-2}(-m_ix_i)/\cE^{-2}(-n_ix_i))_{max}$. The local system $E_i$ is trivial iff $n_i$ is even.
On the hand, $-1\in \mu_2(\cO/\cO(-n_ix_i))$ acts on $\det\H^0(X, \cO/\cO(-n_ix_i))$ as $(-1)^{n_i}$. Taking into account Proposition~\ref{Pp_IC_checkW_description},  our assertion follows.
\end{Prf}

\medskip\noindent
4.3. \select{Whittaker coefficients of Eisenstein series} 

\medskip\noindent  
4.3.1 Main purpose of this section~4.3 is to prove Theorem~\ref{Th_4} below.
Let $\tilde\nu_B:\Bun_B\to\Bunt_G$ be the map sending (\ref{ext_cE^*_by_EotimesOmega}) to $(M,\cB)$ with $\cB=\det\RG(X, \cE\otimes\Omega)$ equipped with the induced automorphism $\cB^2\,\iso\, \det\RG(X,M)$. Let $\cS_B$ be the stack classifying $\cE\in\Bun_1$ and a section $s_2: \cE^2\to \cO_X$. Let $\Four_{\psi}: \D(\Bun_B)\to \D(\cS_B)$ be the Fourier transform.

  For $d\ge 0$ let $\RCov^d$ be the stack classifying $D_1\in {^{rss}X^{(d)}}$, $\cE\in\Bun_1$ and an isomorphism $s_2: \cE^2\,\iso\, \cO(-D_1)$. We have a natural open immersion $\RCov^d\hook{} \cS_B$. For a point of $\RCov^d$ given by $(\cE, s_2)$ we consider $\cO_X\oplus\cE$ as a $\cO_X$-algebra and let $Y=\Spec(\cO_X\oplus \cE)$. Then $Y$ is a smooth projective curve, and the projection  $\phi: Y\to X$ is a degree two covering ramified exactly over $D_1$, see (\cite{L4}, 7.7.2). We denote by $\cE_0$ the $\mu_2$-antiinvariants in $\phi_!\Qlb$, this is a local system over $X-D_1$ extended by zero to $X$. 

 As $(\cE, s_2)\in \RCov^d$ varies, the coverings $Y$ form a family that we denote $\phi_{univ}: Y_{univ}\to X\times\RCov^d$. We denote by $\cE_{0,univ}$ the $\mu_2$-antiinvariants in $(\phi_{univ})_!\Qlb$. 
 
  Given a rank one local system $E$ on $X$, write $AE$ for the corresponding automorphic local system on $\Bun_1$. Let $AJ^d: X^{(d)}\to \Bun_1$ be the Abel-Jacobi map sending $D$ to $\cO(D)$. Then $AE$ is equipped for each $d\ge 0$ with an isomorphism $(AJ^d)^*AE\,\iso\, E^{(d)}$.
 
 Let $E$ be a rank one local system on $X$. Set $K_E=\Eis(AE[\dim\Bun_T])$.
Our aim is to calculate the complex
\begin{equation}
\label{complex_WC_over_RCov_d}
\Four_{\psi}\tilde\nu_B^*K_E[\dimrel(\tilde\nu_B)]\mid_{\RCov^d}
\end{equation}
 
\medskip\noindent
4.3.2 Write $Z^{\theta}_{\Bun_T}\to \Bun_T$ for the relative version of $Z^{\theta}$ over $\Bun_T$, so the fibre of this map over $\cE$ is $Z^{\theta}_{\cE}$. Similarly for the vector bundle $\check{Z}^{\theta}_{\Bun_T}\to X^{(\theta)}\times\Bun_T$. Let 
$$
f_B: Z^{\theta}_{\Bun_T}\to \Bun_B\times X^{(\theta)}
$$ 
be the map sending (\ref{diag_point_of_cX_B}) to the exact sequence (\ref{ext_cE^*_by_EotimesOmega}) together with $D=\div(\cE^{-1}/\cE_1\otimes\Omega)$. This is a morphism of generalized vector bundles over $\Bun_T\times X^{(\theta)}$ given by $\cE^2\otimes\Omega(D)/\cE^2\otimes\Omega\to \H^1(X, \cE^2\otimes\Omega)$. The dual map over $\Bun_T\times X^{(\theta)}$ is denoted
$$
\check{f}_B: \cS_B\times X^{(\theta)}\to \check{Z}^{\theta}_{\Bun_T},
$$
it sends $(\cE, s_2, D)$ to $(D, \cE, v\in \cE^{-2}/\cE^{-2}(-D))$, where $v$ is the image of $s_2$ under the transpose map 
$$
\H^0(X, \cE^{-2})\to \cE^{-2}/\cE^{-2}(-D)
$$  
 
 The perverse sheaves $\IC(W_{max})$ over $Z^{\theta}_{\cE}$ as $\cE$ varies in $\Bun_T$ naturally form a family, which is a perverse sheaf on $Z^{\theta}_{\Bun_T}$ that also will be denoted $\IC(W_{max})$ by abuse of notation. Similarly, $\IC(\check{W}_{max})$ can be seen as a perverse sheaf on $\check{Z}^{\theta}_{\Bun_T}$. As $\cE\in\Bun_1$ varies $A_{\cE}$ becomes a self-dual rank one local system $A$ on $\Bun_T$, and now Proposition~\ref{Pp_calculation_Four_IC(W_max)} over $\check{Z}^{\theta}_{\Bun_T}$ establishes an isomorphism
$$
A\otimes \IC(\check{W}_{max})\,\iso\,\Four_{\psi}(\IC(W_{max}))
$$
 
  Let $\cS_B^0\subset\cS_B$ be the open substack given by $s_2\ne 0$. Let $\nu_{\cS}: \cS_B\to \Bun_1$ be the map sending $(\cE, s_2)$ to $\cE$. 
   
\begin{Lm} Let us describe 
$$
\Four_{\psi}\tilde\nu_B^*K_E[\dimrel(\tilde\nu_B)]\mid_{\cS_B^0}
$$
over the connected component of $\cS_B^0$ given by $2\deg\cE=-d$.
The contribution of the connected component of $\Bunb_{\tilde B}$ given by $\frac{d}{2}-\deg(\cE_1\otimes\Omega)=\theta$ to the latter complex identifies canonically with
$$
(AE)^{-1}_{\Omega}\otimes\nu_{\cS}^*(AE^{-1}\otimes A)\otimes \pr_{1 !}(\pr_2^*(E^*)^{(\theta)}\otimes(\check{f}_B)^*\IC(\check{W}_{max}))[\dimrel(\check{f}_B)]\mid_{\cS^0_B},
$$
where $\pr_1: \cS_B\times X^{(\theta)}\to \cS_B$ and $\pr_2: \cS_B\times X^{(\theta)}\to X^{(\theta)}$ are the projections.
\end{Lm}
\begin{Prf} The Fourier transform exchanges $f_{B !}$ and $(\check{f}_B)^*[\dimrel(\check{f}_B)]$. For a point (\ref{diag_point_of_cX_B}) of $Z^{\theta}_{\Bun_T}$ with $D=\div(\cE^{-1}/\cE_1\otimes\Omega)$ we get
$AE_{\cE^{-1}\otimes\Omega^{-1}(-D)}\,\iso\, AE_{\cE^{-1}\otimes\Omega^{-1}}\otimes (E^{(\theta)}_D)^*$. The locus of $\Bun_B\times_{\Bun_G}\Bunb_{\tilde B}$ given by the property that $\cE_1\otimes\Omega\to M$ factors through $\cE\otimes\Omega\hook{} M$ does not contribute to the Fourier coefficients over $\cS^0_B$. 
\end{Prf}

\medskip

 Let $s_B: \RCov^d\times X^{(\theta)}\to \check{Z}^{\theta}_{\Bun_T}$ be the restriction of $\check{f}_B$.

\begin{Lm} For $\theta\ge 0$ we have canonically 
$$
s_B^*\IC(\check{W}_{max})[\dimrel(\check{f}_B)]\;\iso\; 
\cE_{0,univ}^{(\theta)}[\dim(\RCov^d\times X^{(\theta)})]
$$
\end{Lm}

 This proves the following.
\begin{Th} 
\label{Th_4}
The fibre of (\ref{complex_WC_over_RCov_d}) at $(\cE, s_2: \cE^2\,\iso\, \cO(-D_1))\in \RCov^d$ is the central value of the L-function
\begin{equation}
\label{central_value_Lfunction_for_RCov}
\oplus_{\theta\ge 0} \RG(X^{(\theta)}, (E^*\otimes \cE_0)^{(\theta)})[\theta]
\end{equation}
tensored by $A_{\cE}\otimes AE_{\cE^{-1}\otimes\Omega^{-1}}[\dim \RCov^d]$. Here $\cE_0$ is the sheaf of $\mu_2$-antiinvariants in $\phi_!\Qlb$ for the covering $\phi: Y\to X$ given by $(\cE, s_2)$. If $d>0$ then (\ref{complex_WC_over_RCov_d}) is a shifted local system on $\RCov^d$.
\end{Th}

%  Remind that $\DD(\Eis(AE[\dim\Bun_T]))\,\iso\,\Eis(AE^*[\dim\Bun_T])$, this is compatible with the above answer.

 If $E^*\otimes \cE_0$ is nontrivial then
$$
\RG(X^{(\theta)}, (E^*\otimes \cE_0)^{(\theta)})[\theta]\,\iso\, \wedge^{\theta} \H^1(X, E^*\otimes \cE_0),
$$
this is always the case for $d>0$. Note that $\chi(X, E^*\otimes\cE_0)=2-2g-d$, so the sum in (\ref{central_value_Lfunction_for_RCov}) is over $0\le \theta\le 2g-2+d$. 

\begin{Rem} 
\label{Rem_functional_equation}
i) The functional equation should say that 
$$
\Eis(AE[\dim\Bun_T])\,\iso\, \cC_E\otimes\Eis(AE^{-1}[\dim\Bun_T]),
$$ 
where $\cC_E$ is a 1-dimensional space eventually depending on the local system $E$. This should come from an isomorphism for any point $(\cE, s_2)$ of $\RCov^d$
\begin{multline*}
AE^2_{\cE^{-1}\otimes\Omega^{-1}}\otimes
\wedge^{\theta} \H^1(X, E^*\otimes \cE_0) \,\iso\, \cC_E\otimes\wedge^{2g-2+d-\theta} \H^1(X, E\otimes \cE_0)=\\
\cC_E\otimes\wedge^{\theta}\H^1(X, E\otimes\cE_0)^*\otimes \det\H^1(X, E\otimes \cE_0)=\cC_E\otimes\wedge^{\theta}\H^1(X, E^*\otimes\cE_0)\otimes \det\H^1(X, E\otimes \cE_0)
\end{multline*}
This would follow from an isomorphism 
\begin{equation}
\label{iso_desired_def_cC_E}
AE^2_{\cE^{-1}\otimes\Omega^{-1}}\,\iso\,  \cC_E\otimes\det\H^1(X, E\otimes \cE_0)
\end{equation}
Note that $AE^2_{\cE^{-1}}\,\iso\, AE_{\cO(D_1)}$. By Lemma~\ref{Lm_detRG_X_E_for_rank_one} below,
$$
\det\RG(X,E)\otimes (AE)_{\Omega}\,\iso\,\Qlb
$$ 
and $\det\RG(Y, \phi^*E)\otimes (A(\phi^*E))_{\Omega_Y}\,\iso\, \Qlb$.
Note that $\Omega_Y\,\iso\, (\phi^*\Omega)(D_1)$, we may view $D_1$ as a divisor on $Y$ as well as on $X$. One has canonically $A(\phi^*E)\,\iso\, N^*(AE)$, where $N: \Bun_{1,Y}\to\Bun_1$ is the norm map, and $N(\Omega_Y)\,\iso\, \Omega^2(D_1)$. Combining the above, we get an isomorphism
$$
\det\RG(X, E\otimes\cE_0)\otimes (AE)_{\Omega(D_1)}\,\iso\, \Qlb,
$$
so $\det\H^1(X, E\otimes \cE_0)\,\iso\, (AE)_{\Omega(D_1)}$. To get (\ref{iso_desired_def_cC_E}) we have to set $\cC_E=(AE)^{-3}_{\Omega}$. 

 Assume $E^2$ nontrivial. If we knew that the contribution of each connected component of $\Bunb_{\tilde B}$ to $\Eis(E[\dim\Bun_T])$ is an irreducible perverse sheaf then our calculation of Whittaker coefficients over $\RCov^d$ would imply the functional equation (the corresponding isomorphism over $\RCov^d$ would extend via the intermediate extension to the whole of $\cS_B$).
 
\medskip\noindent
ii) Consider the automorphism $w$ of $\Bun_1$ sending $\cE$ to $\cE^*\otimes\Omega^{-3}$. Our calculation predicts the form of the functional equation in general. Namely, we conjecture an isomorphism functorial in $K\in \D(\Bun_T)$
$$
\Eis(w^*K)\,\iso\, \Eis(K)
$$
for $K\in\D(\Bun_1)$ satisfying the following regularity assumption, 
the same as for the functor $\Eis_{\GL_1}^H$ in Conjecture~\ref{Con_compatibility_parab_induction}. For $\cL\in\Bun_1$ let $r_{\cL}: \Bun_1\to\Bun_1$ be the map $\cB\mapsto \cL\otimes\cB^2$. Then $K\in \D(\Bun_1)$ is regular for $\Eis_{\GL_1}^H$ iff for any $\cL\in\Bun_1$ one has $\RG_c(\Bun_1, r^*_{\cL}K)=0$.

If the functional equation holds for $K$ which is a bounded complex on each connected component of $\Bun_T$ then $\Eis(K)$ is also bounded.
In particular, the functional equation does not hold for $\Eis(AE)$, where $E$ is a rank one local system on $X$ such that $E^2\,\iso\, \Qlb$, because $\Eis(AE)$ is unbounded in this case. If $g=0$ then there is no functional equation for $\Eis$. 
 The local analog of this functional equation is the property of the principal series representations of the metaplectic group (\cite{WTG}, Proposition~2.1).
 
  Note that the functional equation for $\Eis_{\GL_1}^H$ claims that for $K\in\D(\Bun_1)$ regular one has $\Eis_{\GL_1}^H(\bar w^*K)\,\iso\, \Eis_{\GL_1}^H(K)$, where $\bar w:\Bun_1\,\iso\,\Bun_1$ sends $\cB$ to $\cB^{-1}\otimes\Omega^{-2}$. 
\end{Rem}
\begin{Rem}
\label{Rem_two_remarks_LS}
 The cup product equips $\H^1(X,\Qlb)$ with the structure of a symplectic vector space, so $\det\RG(X,\Qlb)\,\iso\, \Qlb$ canonically. Similarly, $\det\RG(Y,\Qlb)\,\iso\,\Qlb$, so $\det\RG(X, \cE_0)\,\iso\,\Qlb$ canonically (we ignore the Tate twists).
\end{Rem}

\begin{Lm} 
\label{Lm_detRG_X_E_for_rank_one}
Let $E$ be a rank one local system on $X$. One has a canonical isomorphism $\det\RG(X,E)\,\iso\, (AE)^*_{\Omega}$.
\end{Lm}
\begin{Prf} By Remark~\ref{Rem_two_remarks_LS}, we may assume $E$ nontrivial. Then $\det\RG(X,E)$ is canonically the dual of the 1-dimensional space $\RG(X^{(2g-2)}, E^{(2g-2)})[2g-2]$. Clearly, we may assume $g\ge 2$. 
Let $\uPic^{2g-2} X$ be the Picard scheme classifying line bundles of degree $2g-2$, write $AJ: X^{(2g-2)}\to \uPic^{2g-2} X$ for the map $D\mapsto \cO_X(D)$. Let $i: \Spec k\to \uPic^{2g-2} X$ be the point $\Omega$. One has an exact triangle 
$$
\RG(\PP^{g-2},\Qlb)[2g-2]\to (AJ)_!\Qlb[2g-2] \to i_*\Qlb
$$ 
on $\uPic^{2g-2} X$. One has canonically $E^{(2g-2)}\,\iso\, AJ^*(AE)$, and $\RG(\uPic^{2g-2} X, AE)=0$. So, the above triangle yields an isomorphism
$$
\RG(X^{(2g-2)}, E^{(2g-2)})[2g-2]\,\iso\, (AE)_{\Omega}
$$
\end{Prf}

\medskip\noindent
4.3.3 \select{Some residues of $\Eis$}

\medskip\noindent
From (\ref{central_value_Lfunction_for_RCov}) we also see at least some residues of $\Eis(AE[\dim\Bun_1])$. Namely, if $E^2$ is trivial, then $E$ defines a $\mu_2$-torsor over $X$ that we see as a point $(\cE, \cE^2\,\iso\, \cO_X)$ of $\RCov^0$. 
For the corresponding automorphism $\tilde\sigma_{\cE}$ of $\Bunt_G$ defined in Section~5.1 below, the perverse sheaf $\tilde\sigma_{\cE}^*\Aut$ appears as a direct summand in $\Eis(AE[\dim\Bun_1])$. This is a residue of $\Eis(AE[\dim\Bun_1])$, which is unbounded in this case as we see from Theorem~\ref{Th_4}.

\medskip\noindent
4.3.4 \select{Constant terms of $\Eis$}

\medskip\noindent
Remind the map $\tilde\nu_B:\Bun_B\to\Bunt_G$ from Section~4.3.1. Let 
$\nu_T: \Bun_B\to\Bun_T$ be the natural map. Let $\Bun^d_B$ be the connected component of $\Bun_B$ classifying exact sequences (\ref{ext_cE^*_by_EotimesOmega}) with $\deg(\cE)=d$. Write $\Bun_1^d$ for the connected component of $\Bun_1$ classifying $\cE\in\Bun_1$ of degree $d$.

Let $CT: \D(\Bunt_G)\to\Bun_T$ be the functor 
$$
CT(K)=\nu_{T !}\tilde\nu_B^*K[\dimrel(\tilde\nu_B)]
$$ 
Write $CT^d(K)$ for the contribution of $\Bun^d_B$ to the above functor, so $CT(K)=\oplus_{d\in\ZZ} CT^d(K)$. 

\begin{Def} Let $K\in \D(\Bun_1^{d_1})$, $d,d_1\in\ZZ$ with $d-d_1$ even.
\begin{itemize}  
\item For $d_1\le 2-2g-d$ define $F_0(K)\in \D(\Bun_1^d)$ as follows. Let $r=1-g-\frac{d+d_1}{2}$. Let $a: X^{(r)}\times \Bun_1^d\to \Bun^{d_1}_1$ be the map sending $(D, \cE)$ to $\cE^{-1}\otimes\Omega^{-1}(-2D)$. Set $F_0(K)=\pr_{2 !}a^*K[1-g]$, where $\pr_2: X^{(r)}\times \Bun_1^d\to \Bun_1^d$ is the projection.
  
\item For $d_1\le d$ define $F_1(K)\in \D(\Bun_1^d)$ as follows. Let $\bar r=(d-d_1)/2$. Let $b: X^{(\bar r)}\times\Bun_1^d\to \Bun_1^{d_1}$ be the map $(D, \cE)\mapsto \cE(-2D)$. Set $F_1(K)=\pr_{2 !}b^*K[2-2g+d-d_1]$, here $\pr_2: X^{(\bar r)}\times \Bun_1^d\to \Bun_1^d$ is the projection.
\end{itemize}
\end{Def}
 
   The following result is compatible with the form of the functional equation from Remark~\ref{Rem_functional_equation}. 
  
\begin{Pp} Let $K\in \D(\Bun_1^{d_1})$. The complex
$K_{d,d_1}:=CT(\Eis(K))\mid_{\Bun_1^d}$ vanishes unless $d-d_1$ is even. In the latter case it is as follows.
\begin{itemize}
\item[1)] If $d_1> \max\{d, 2-2g-d\}$ then $K_{d,d_1}=0$.
\item[2)] If $d<d_1\le 2-2g-d$ then $K_{d,d_1}\,\iso\, F_0(K)$ canonically.

\item[3)] If $d\ge d_1 > 2-2g-d$ then $K_{d,d_1}\,\iso\, F_1(K)$ canonically.

\item[4)] If $d_1\le \min\{d, 2-2g-d\}$ then there is an exact triangle $F_0(K)\to K_{d,d_1}\to F_1(K)$ over $\Bun_1^d$. 
\end{itemize}
\end{Pp}
\begin{Prf}
We have to integrate with respect to the composition $\Bun^d_B\times_{\Bun_G}\Bunb_{B}^{d_1}\to \Bun^d_B\to \Bun_1^d$. Write a point of $\Bun^d_B\times_{\Bun_G}\Bunb_{B}^{d_1}$ as the diagram (\ref{diag_point_of_cX_B}). \\
1) In this case $\Hom(\cE_1, \cE)=\Hom(\cE_1\otimes\Omega, \cE^{-1})=0$.\\
2) In this case $\Hom(\cE_1, \cE)=0$, and there remains the direct image of $\IC(W_{max})\otimes\tau^*K[1-g]$ under $Z^{\theta}_{\Bun_1^d}\to \Bun_1^d$ for $\theta=2-2g-d-d_1$, here $\tau: \Bun^d_B\times_{\Bun_G}\Bunb_{B}^{d_1}\to \Bun^{d_1}_1$ is the projection. Applying Proposition~\ref{Pp_restriction_zero_section}, we get the desired answer.

\noindent
3) In this case $\Hom(\cE_1\otimes\Omega, \cE^{-1})=0$, so $\tilde t$ factors as $\tilde t: \cE_1\otimes\Omega\hook{} \cE\otimes\Omega$, and our claim follows from Theorem~\ref{Th_3}.

\noindent
4) as in the above cases, stratify $\Bun^d_B\times_{\Bun_G}\Bunb_{B}^{d_1}$ by the property that $\tilde t$ factors as $\cE_1\otimes\Omega\hook{} \cE\otimes\Omega$ or not, calculate the direct image with respect to this stratification.
\end{Prf}

\medskip\noindent
4.4 \select{Whittaker sheaves for $\wt\SL_2$}

\medskip\noindent
Let $E$ be a rank one local system on $X$. If we want $\Eis(AE[\dim\Bun_1])$
to depend only on $E\oplus E^*$, we have to choose a square root of $\Omega$, which we avoided before. This is analogous to a choice of such square root in (\cite{BG}, Section~2.2.1) for Eisenstein series for reductive groups.
  
Fix once and for all a line bundle $\Omega^{1/2}$, whose square is identified with $\Omega$. Set 
\begin{equation}
\label{sheaf_cK_E_plus_E*}
\cK_{E\oplus E^*}=(AE)_{\Omega^{3/2}}\otimes \Eis(AE[\dim\Bun_T])
\end{equation}
Then according to the functional equation, $\cK_{E\oplus E^*}$ should depend only on the $\SL_2$-local system $E\oplus E^*$. 

\begin{Rem} Here is a consistency check of the normalization in (\ref{sheaf_cK_E_plus_E*}). In (\cite{BG}, Section~2.1.7) a notion of a twisted $W$-action on $\Bun_T$ was introduced, here $T$ is a maximal torus of a reductive group, and $W$ the corresponding Weyl group. For the maximal torus $T_H=\GL_1$ of $H=\SO_3$ the nontrivial element of the Weyl group sends $\cB\in\Bun_1$ to $\cB^*\otimes\Omega^{-2}$ under this twisted action. As in (\cite{BG}, Section~2.1.7), given a $\check{T}_H$-local system $E_{\check{T}_H}$, one defines the automorphic sheaf $\Aut_{E_{\check{T}_H}}$ on $\Bun_{T_H}$ as $(AE)_{\Omega}\otimes AE_{\check{T}_H}[\dim\Bun_{T_H}]$. For the map $\delta: \Bun_1\to \Bun_1$, $\cB\mapsto \cB\otimes\Omega^{1/2}$ from Conjecture~\ref{Con_compatibility_parab_induction} the fibre of $\delta^*\Aut_{E_{\check{T}_H}}$ at $\cB\in\Bun_1$ is $(AE)_{\Omega^{3/2}}\otimes (AE)_{\cB}$. 
\end{Rem}

\smallskip

 Remind the 2-automorphism $\bar\epsilon$ on $\RCov^d$ introduced in Section~0.3.3. Define the local system $\cS^d_{E\oplus E^*}$ on $\RCov^d$ as (\ref{complex_WC_over_RCov_d}) tensored by $(AE)_{\Omega^{3/2}}[-\dim\RCov^d]$. It is $\ZZ/2\ZZ$-graded by the action of $\bar\epsilon$. 
 
  Remind that for a local system $V$ on $X$ we introduced a complex $CL^d_V$ on $\RCov^d$ in Section~0.3.5.

\begin{Cor} Over $\RCov^d$ one has a canonical (automatically $\ZZ/2\ZZ$-graded) isomorphism 
$$
(\cS_{E\oplus E^*})^{\otimes 2}\,\iso\, CL^d_{E\oplus E^*} 
$$
\end{Cor}  

 Based on this result, in particular, we proposed Conjecture~\ref{Con_existence_cS} in Section~0.3.5. Remind the theta-lifting functors $F_G$ and $F_H$ from Section~0.3.2. For a $\SL_2$-local system $E$ on $X$ remind our notation $\Aut_E\in \D(\Bun_H)$ from Section~0.3.3. 
     
  Assume Conjecture~\ref{Con_existence_cS} true just for the point $\cE=\cO$ of $\RCov^0$, and denote by $SQ_E$ the fibre of $\cS_E^0$ at $\cO\in \RCov^0$. So, $SQ_E$ is $\ZZ/2\ZZ$-graded by the action of $\bar\epsilon$, and can be viewed as a complex on $B(\mu_2)$. Then we get a $\ZZ/2\ZZ$-graded isomorphism
\begin{equation}
\label{iso_SQ_E_defining}
(SQ_E)^{\otimes 2}\,\iso\, \oplus_{\theta\ge 0} \RG(X^{(\theta)}, E^{(\theta)})[\theta]
\end{equation}

Write $\und{\RCov}^d$ for the coarse moduli space of $\RCov^d$. The $\mu_2$-gerb $\RCov^d\to \und{\RCov}^d$ is trivial, its trivialization yields a projection $\und{p}: \RCov^d\to B(\mu_2)$.

 Remind that each $K\in \D_-(\Bunt_G)$ is $\ZZ/2\ZZ$-graded by the action of $\bar\epsilon$. The following conjecture is motivated by our calculation of geometric Waldspurger periods (Theorem~5 and Conjecture~3 in \cite{L3}), on one hand, and by a relation between $F_G, F_H$ and the first Whittaker coefficients functors for $\wt\SL_2$ and for $H$, on the other hand (we are planning to discuss this relation in a separate paper). 
 
\begin{Con} 
\label{Con_2}
Let $E$ be an irreducible $\SL_2$-local system on $X$. Then 

\smallskip\noindent
i) $\cK_E:=QL(\Aut_E)$ is a perverse sheaf, whose each $\ZZ/2\ZZ$-parity piece is irreducible, and $\DD(\cK_E)\,\iso\, \cK_E$.

\smallskip\noindent
ii) The complex $F_HF_G(\Aut_E)$ is isomorphic to $\Aut_E$ tensored by
$\oplus_{\theta\ge 0} \RG(X^{(\theta)}, E^{(\theta)})[\theta]$. 

\smallskip\noindent
iii) Assume Conjecture~\ref{Con_existence_cS} true for the point $\cE=\cO$ of $\RCov^0$. For a suitable map $\und{p}: \RCov^d\to B(\mu_2)$ there is an isomorphism over $\RCov^d$
$$
(\Four_{\psi}\tilde\nu_B^*F_G(\Aut_E))\mid_{\RCov^d} \,\iso\, \und{p}^*SQ_E\otimes 
(\Four_{\psi}\tilde\nu_B^*\cK_E)\mid_{\RCov^d}
$$
\end{Con}
  
  The above conjecture suggests the following. For an irreducible $\SL_2$-local system $E$ on $X$ one wants to construct the local system $\cS^d_E$ of Conjecture~\ref{Con_existence_cS}. First, construct the local system $\und{p}^*SQ_E\otimes\cS_E^d$ via the theta-lifting as 
$$
\Four_{\psi}\tilde\nu_B^*F_G(\Aut_E)[\dimrel(\tilde\nu_B)-\dim\RCov^d]\mid_{\RCov^d}
$$
Then to get  $\cS_E^d$ it would suffice to get rid of the factor $SQ_E$. 
 
\begin{Rem} i) For an irreducible $\SL_2$-local system $E$ on $X$ Conjecture~\ref{Con_QL} would imply an isomorphism $\R\Hom(\cK_E, \cK_E)\,\iso\, \R\Hom(\Aut_E, \Aut_E)$, which should in turn yield the Rallis inner product formula for $F_G(\Aut_E)$ (see \cite{WTG}, Section~7.7). The latter expresses $\RG_c(\Bunt_G, F_G(\Aut_E)\otimes\DD F_G(\Aut_E))$ via $\RG(\Bun_H, \Aut_E\otimes\DD(\Aut_E))$.

\smallskip\noindent
ii) Let $E$ be an irreducible $\SL_2$-local system on $X$. If we were working over complex numbers with $\D$-modules, then $\H^1(X,E)$ would carry a pure Hodge structure of weight 1, this structure yields a natural candidate for $SQ_E$. 

\smallskip\noindent
iii) If $E$ is an $\SL_2$-local system on $X$ with $\H^0(X,E)=0$ then $\H^1(X, E)$ can be seen as a torsor over $\Spec k$ under $\SO_{4g-4}$. The datum of $SQ_E$ together with a $\ZZ/2\ZZ$-graded isomorphism (\ref{iso_SQ_E_defining}) is equivalent to a datum of a lifting of this torsor to a $\Spin_{4g-4}$-torsor. This was already observed in (\cite{L3}, Remark~14 and Conjecture~5). If we were working with $\cD$-modules in characteristic zero, we could consider the moduli stack $LS_{\SL_2}$ of $\SL_2$-local systems on $X$. Then the corresponding $\Spin_{4g-4}$-torsor can be seen as a torsor over the open substack of $LS_{\SL_2}$ given by $\H^0(X, E)=0$. Over the locus of this substack, where $E=V\oplus V^*$, it should coincide with the torsor given by the complex $\cS^0_{V\oplus V^*}$ defined in (\ref{complex_K_E_after_transfrom_on_RCov}). 
\end{Rem}

\bigskip
\centerline{\scshape 5. First Whittaker coefficients}

\bigskip\noindent
5.1 Use the notation of Section~0.3.1. Note that $\Aut_g\in \D_{\mp}(\Bunt_G)$ and $\Aut_s\in\D_{=}(\Bunt_G)$. 

Remind the stack $\RCov^0$ introduced in Section~0.3.3. Given $(\cE, s_2)\in\RCov^0$ let $\sigma_{\cE}:\Bun_G\,\iso\,\Bun_G$ be the isomorphism sending $M$ to $M\otimes\cE$. 
\begin{Lm} The trivialization $\cE^2\,\iso\, \cO_X$ yields a uniquely defined
$\ZZ/2\ZZ$-graded isomorphism of line bundles on $\Bun_G$
$$
\cA^{-1}\otimes\sigma_{\cE}^*\cA\otimes\;\iso\; \frac{\det\RG(X, \cE)^{2n}}{\det\RG(X, \cO)^{2n}}
$$
\end{Lm}
\begin{Prf}
By (\cite{L3}, Lemma~1), for $M\in \Bun_G$ there is a canonical $\ZZ/2\ZZ$-graded isomorphism 
$$
\det\RG(X, M\otimes\cE)\;\iso\; \frac{\det\RG(X,M)\otimes \det\RG(X, \cE)^{2n-1}\otimes\det\RG(X, \cE\otimes\Omega^n)}{\det\RG(X, \Omega^n)\otimes\det\RG(X, \cO)^{2n-1}}
$$
and 
$$
\frac{\det\RG(X, \cE\otimes\Omega^n)\otimes\det\RG(X, \cE)^{n-1}}{\det\RG(X, \cE\otimes\Omega)^n}\;\iso\; \frac{\det\RG(X,\Omega^n)}{\det\RG(X, \cO)}
$$ 
The trivialization $\cE^2\,\iso\,\cO_X$ yields a $\ZZ/2\ZZ$-graded isomorphism $\det\RG(X, \cE\otimes\Omega)\,\iso\, \det\RG(X, \cE)$. Combining the above, we obtain 
$$
\frac{\det\RG(X, M\otimes\cE)}{\det\RG(X, M)}\,\iso\, \frac{\det\RG(X, \cE)^{2n}}{\det\RG(X, \cO)^{2n}}
$$ 
\end{Prf}  

\medskip

 We extend $\sigma_{\cE}$ to an isomorphism $\tilde\sigma_{\cE}: \Bunt_G\,\iso\, \Bunt_G$ sending $(M,\cB)$, where $\cB$ is a line equipped with $\cB^2\,\iso\, \det\RG(X,M)$ to $(M\otimes\cE, \cB')$, where 
$$
\cB'=\cB\otimes \det\RG(X, \cE)^n\otimes \det\RG(X,\cO)^{-n}
$$ 
is equipped with the isomorphism $\cB^{' 2}\,\iso\, \det\RG(X, M\otimes\cE)$ of the above lemma. 

 Note that $\RCov^0$ is a group stack, which identifies with $\Bun_{\mu_2}$ naturally. One checks that the above gives an action of $\Bun_{\mu_2}$ on $\Bunt_G$. Note that $\tilde\sigma_{\cE}^*$ is compatible with the 2-automorphisms $\epsilon, \bar\epsilon$. So, $\tilde\sigma_{\cE}^*$ preserves each of the categories $\D_{=}(\Bunt_G)$ and $\D_{\mp}(\Bunt_G)$.
 
\begin{Lm} 
\label{Lm_17}
Given $\cE\in\Bun_{\mu_2}$ let $\cE_0$ denote the corresponding rank one and order two local system on $X$. Then the functors $\D_-(\Bunt_G)\to \D_-(X\times\Bunt_G)$ given by 
$$
K\mapsto \H^{\ra}(\cA_{\alpha}, \tilde\sigma_{\cE}^*K)\;\;\;\;\mbox{and by}\;\;\;\; K\mapsto (\id\times\tilde\sigma_{\cE})^*\H^{\ra}(\cA_{\alpha}, K)\otimes\pr_1^*\cE_0
$$ are naturally isomorphic. Here $\pr_1: X\times\Bunt_G\to X$ is the projection.
So, if $E$ is a $\Sp_{2n}$-local system and $K$ is a $E$-Hecke eigensheaf then $\tilde\sigma_{\cE}^*K$ is equipped with the structure of a $E\otimes\cE_0$-Hecke eigensheaf.
\end{Lm}
\begin{Prf}
Let $W$ be the local system on $B(\mu_2)$ corresponding to the character $\mu_2\hook{}\Qlb^*$. The $\mu_2$-torsor $\cE$ on $X$ defines a morphism $X\to B(\mu_2)$, and the inverse image of $W$ under this map identifies with $\cE_0$. Consider the stack $\wt\cH^{\alpha}_G$ from (\cite{L1}, Section~9.2), we have the isomorphisms $\kappa, \kappa'$ given by (\cite{L1}, Lemma~14). 
Remind that $\wt\cH^{\alpha}_G$ classifies $(x, M,\cB, M',\cB')$, where $M$ is in the position $\alpha$ with respect to $M'$ at $x$, and $\cB^2\,\iso\,\det\RG(X, M)$, $\cB^{' 2}\,\iso\,\det\RG(X, M')$. Let $\sigma_{\alpha,\cE}: \wt\cH^{\alpha}_G\to \wt\cH^{\alpha}_G$ be the morphism sending the above point to 
$$
(M\otimes\cE, \bar\cB, M'\otimes\cE, \bar\cB'),
$$
where $\bar\cB=\cB\otimes \det\RG(X, \cE)^n\otimes\det\RG(X,\cO)^{-n}$ and $\bar\cB'=\cB'\otimes \det\RG(X, \cE)^n\otimes\det\RG(X,\cO)^{-n}$. 
Let 
$$
\bar\sigma_{\alpha,\cE}: (\Bunt_G\times_{\Bun_G} \cH^{\alpha}_G)\times B(\mu_2)\to (\Bunt_G\times_{\Bun_G} \cH^{\alpha}_G)\times B(\mu_2)
$$
be the morphism sending $(x, M,\cB, M', \cB_0)$ with $\cB_0^2\,\iso\, k$ and $\cB^2\,\iso\, \det\RG(X,M)$ to 
$$
(x, M\otimes\cE, \bar\cB, M'\otimes\cE, \cB_0\otimes\cE_x),
$$
where as above $\bar\cB=\cB\otimes \det\RG(X, \cE)^n\otimes\det\RG(X,\cO)^{-n}$.

Our claim follows from the fact that the diagram
$$
\begin{array}{ccc}
\wt\cH^{\alpha}_G & \toup{\kappa} & (\Bunt_G\times\Bun_G \cH^{\alpha}_G)\times B(\mu_2)\\
\downarrow\lefteqn{\scriptstyle \sigma_{\alpha,\cE}} && \downarrow\lefteqn{\scriptstyle \bar\sigma_{\alpha,\cE}}\\
\wt\cH^{\alpha}_G & \toup{\kappa} & (\Bunt_G\times\Bun_G \cH^{\alpha}_G)\times B(\mu_2)
\end{array}
$$
is canonically 2-commutative, and $K$ changes under the action of $B(\mu_2)$ by $W$. In particular, the inverse image of $K$ under the map $X\times\Bunt_G\to \Bunt_G$, $(x, M,\cB)\mapsto (M, \cB\otimes\cE_x)$ identifies with $\cE_0\boxtimes K$.
\end{Prf}

\medskip

 Consider the split group $\bar H=\GSpin_{2n+1}$, it fits into an exact sequence $1\to \Gm\to \bar H\to H\to 1$. So, $\Bun_1$ acts naturally on $\Bun_{\bar H}$, and the extension of scalars map $\nu_{\bar H}:\Bun_{\bar H}\to \Bun_H$ is invariant under this action. Given $(\cE,s_2)\in \RCov^0$, let $\cE_0$ be the $\mu_2$-antiinvariants in $\phi_!\Qlb$ for the corresponding covering $\phi: Y\to X$. Let $\nu_1: \Bun_{\bar H}\to \Bun_1$ be the extension of scalars map. Clearly, $\nu_1^*A\cE_0$ is equivariant under the above action of $\Bun_1$ on $\Bun_{\bar H}$, so there is a local system $\cE_{0,H}$ on $\Bun_H$ equipped with an isomorphism $\nu_{\bar H}^*\cE_{0,H}\,\iso\, \nu_1^*A\cE_0$ over $\Bun_{\bar H}$. For $\gamma\in \pi_1(\bar H)$ let $\bar\gamma\in \pi_1(H)$ be its image, let $\Bun_{\bar H}^{\gamma}$ be the corresponding connected component of $\Bun_{\bar H}$ and similarly for $H$. Since the fibres of $\nu_{\bar H}: \Bun_{\bar H}^{\gamma}\to \Bun_H^{\bar\gamma}$ are connected, $\cE_{0,H}$ is defined up to a unique isomorphism (cf. \cite{G}, Lemma~4.8). 

 According to Lemma~\ref{Lm_17}, the functor $\tilde\sigma_{\cE}^*$ on $\D_-(\Bunt_G)$ should correspond via the equivalence of Conjecture~\ref{Con_QL} to the tensoring by $\cE_{0,H}$ on $\D(\Bun_H)$. 
 
\medskip\noindent
5.2 Assume $n=1$. The first Whittaker coefficients for $\Bunt_G$ are the functors indexed by the points of $\RCov^0$. Given $\cE$ with a trivialization $\cE^2\,\iso\,\cO_X$, let 
\begin{equation}
\label{functor_W_cE_first_Whit}
W_{\cE}: \D(\Bunt_G)\to \D(\Spec k)
\end{equation}
 be the following functor. Let $R=\cE^2\otimes\Omega$, this is a group scheme on $X$. The stack $\Bun_R$ classifies exact sequences (\ref{ext_cE^*_by_EotimesOmega}) on $X$. 
Let $\tilde\nu_R: \Bun_R\to\Bunt_G$ be the map sending (\ref{ext_cE^*_by_EotimesOmega}) to $(M,\cB)$, where $\cB=\det\RG(X, \cE\otimes\Omega)$ is equipped with the $\ZZ/2\ZZ$-graded isomorphism $\cB^2\,\iso\,\det\RG(X,M)$ induced by (\ref{ext_cE^*_by_EotimesOmega}). 
Then (\ref{functor_W_cE_first_Whit}) is defined as the $*$-fibre of $\Four_{\psi}\tilde\nu_R^*[\dimrel(\tilde\nu_R)]$ at the identity.

\begin{Pp} For any $\cE\in\RCov^0$ there is an isomorphism of functors $W_{\cO}\comp \tilde\sigma_{\cE}^*\,\iso\, W_{\cE}$ from $\D(\Bunt_G)$ to $\D(\Spec k)$.
\end{Pp}
\begin{Prf}
Let $R=\Omega$ and $R'=\cE^2\otimes\Omega$. The map $\tilde\nu_R$ sends an exact sequence 
\begin{equation}
\label{ext_cO_by_Omega}
0\to \Omega\to M\to \cO_X\to 0
\end{equation}
on $X$ to $(M,\cB)$, where $\cB=\det\RG(X,\Omega)$. Then $\tilde\sigma_{\cE}(M,\cB)=(M\otimes\cE, \cB')$, where $\cB'=\cB\otimes \det\RG(X,\cE)\otimes\det\RG(X,\cO)^{-1}\,\iso\, \det\RG(X, \cE\otimes\Omega)$ is equipped with $\cB^{' 2}\,\iso\,\det\RG(X, M\otimes\Omega)$. So, the diagram is naturally 2-commutative
$$
\begin{array}{ccc}
\Bun_{\Omega} & \toup{\nu_{\cE}} & \Bun_{\cE^2\otimes\Omega}\\
\downarrow\lefteqn{\scriptstyle \tilde\nu_R} && \downarrow\lefteqn{\scriptstyle \tilde\nu_{R'}}\\
\Bunt_G & \toup{\tilde\sigma_{\cE}} & \Bunt_G,
\end{array}
$$
where $\nu_{\cE}$ is the map tensoring the exact sequence (\ref{ext_cO_by_Omega}) by $\cE$. 
\end{Prf}

\medskip
\begin{Rem} The functors $W_{\cE}$ are pairwise non-isomorphic for different points $\cE\in\RCov^0$. Indeed, for the theta-sheaf $\Aut$ we have $W_{\cE}(\Aut)\ne 0$ iff $\cE\,\iso\, \cO$. Given $\cE\in\RCov^0$ let $K=\tilde\sigma_{\cE *}\Aut$. Then for $\cE'\in\RCov^0$ we get $W_{\cE'}(K)\ne 0$ iff $\cE\,\iso\, \cE'$.
\end{Rem}
  
\bigskip

\centerline{\scshape Appendix A. More about Eisenstein series on $\Bunt_G$}  

\bigskip\noindent
A.1 Keep notations of Section~4. Remind the Shatz stratification of $\Bun_G$.
Let $Shatz^0\subset\Bun_G$ be the open substack of semi-stable $G$-bundles. So, $M\in Shatz^0$ iff for any rank one subsheaf $L\subset M$ one has $\deg L\le g-1$. For $d>0$ let $Shatz^d\subset\Bun_G$ be the stack classifying $\cE\in\Bun_1^{d+1-g}$ and an exact sequence $0\to\cE\otimes\Omega\to M\to \cE^{-1}\to 0$. The map $Shatz^d\to\Bun_G$ sending this sequence to $M$ is a locally closed immersion. Besides, $Shatz^d$ for $d\ge 0$ form a stratification of $\Bun_G$ by locally closed substacks. The stack $Shatz^d$ is irreducible of dimension $2g-2-2d$ for $d>0$, and $\dim Shatz^0=\dim\Bun_G=3g-3$.

  Write $\Bunb_B^d$ for the connected component of $\Bunb_B$ given by $\deg \cE=d$ for $(\cE\otimes\Omega\hook{} M)\in\Bunb_B$. The stack $\Bunb^d_B$ is smooth irreducible of dimension $-2d$.
  
  For $d>0$ the image of $\bar\gp^{d+1-g}: \Bunb^{d+1-g}_B\to\Bun_G$ is the closure $\ov{Shatz}^d$ of $Shatz^d$. Let $Shatz^d_{\tilde G}$ (resp., $\ov{Shatz}^d_{\tilde G}$) denote the preimage of $Shatz^d$ (resp., of $\ov{Shatz}^d$) under $\Bunt_G\to\Bun_G$. Let $p_{Sh}: Shatz^d\to \Bun_1^{d+1-g}$ be the map sending  (\ref{ext_cE^*_by_EotimesOmega}) to $\cE$.
  
   For $d>0$ we have an isomorphism 
$$
Shatz^d\times B(\mu_2)\,\iso\, Shatz^d_{\tilde G}
$$ 
sending $(\cE, \cB_0, \cB_0^2\,\iso\, k)$ and the exact sequence (\ref{ext_cE^*_by_EotimesOmega}) to (\ref{ext_cE^*_by_EotimesOmega}) and $\cB=\cB_0\otimes \det\RG(X, \cE\otimes\Omega)$ with the induced isomorphism $\cB^2\,\iso\, \det\RG(X, M)$. For $d>0$ and a rank one local system $E$ on $X$ let $\IC(AE, d)$ denote the intermediate extension of 
$(p_{Sh}^*AE)\boxtimes W[\dim Shatz^d]$ 
under $Shatz^d_{\tilde G}\to \Bunt_G$. Note that $\IC(AE, d)$ lies in $\D_{\mp}(\Bunt_G)$ (resp., in $\D_{=}(\Bunt_G)$) for $d$ even (resp., for $d$ odd). For $d>0$ the map 
$$
\bar\gp^{d+1-g}: \Bunb^{d+1-g}_{\tilde B}\to \ov{Shatz}^d_{\tilde G}
$$ 
is an isomorphism over $Shatz^d_{\tilde G}$. 
 
\begin{Pp} 
\label{Pp_A1}
Let $E$ be a rank one local system on $X$ and $d>0$. Consider the complex 
\begin{equation}
\label{complex_Eis_d+1-g}
\bar\gp^{d+1-g}_!(\IC_{B,W}\otimes \bar q^*(AE))
\end{equation} 
If $E^2$ is not trivial then it is canonically isomorphic to $\IC(AE, d)$. If $E^2\,\iso\,\Qlb$ then (\ref{complex_Eis_d+1-g}) identifies with 
$$
\oplus_{k\ge 0} \IC(AE, d+2k)
$$
\end{Pp}
\begin{Prf}
We know already that $\IC(AE, d)$ appears in (\ref{complex_Eis_d+1-g})
with multiplicity one. Consider a point $(\cE\otimes\Omega\hook{} M)$ of $Shatz^r$ for some $r>d$. Set $\cB=\det\RG(X, \cE\otimes\Omega)$, so $(M,\cB)\in Shatz^r_{\tilde G}$ naturally. Let $K$ denote the $*$-fibre of (\ref{complex_Eis_d+1-g}) at this point.

The fibre of $\bar\gp^{d+1-g}$ over $(M,\cB)$ identifies with 
$X^{(r-d)}$, namely to $D\in X^{(r-d)}$ there corresponds the subsheaf $\cE\otimes\Omega(-D)\subset M$. From Theorem~\ref{Th_3} we learn that
$K=0$ unless $r-d$ is even, and for $r-d$ even we get $K\,\iso\, (AE)_{\cE}\otimes\RG(X^{(\frac{r-d}{2})}, (E^{-2})^{(\frac{r-d}{2})})[2g-2-d-r]$. The codimension of $Shatz^r$ in $\ov{Shatz}^d$ is $2(r-d)$. If $\H^2(X, E^{-2})=0$, we see that the $*$-restriction of (\ref{complex_Eis_d+1-g}) is placed in perverse degrees $<0$. The first assertion follows.

 If $E=\Qlb$ then the same calculation shows that $\IC(A\Qlb, r)$ appears in (\ref{complex_Eis_d+1-g}) with multiplicity one, and the second assertion follows. 
\end{Prf}  

\medskip\noindent
A.2 Assume $g=0$ and $n=1$. We give some precisions about the functor $QL$ (cf. Conjecture~\ref{Con_compatibility_parab_induction}) in this case.
Write $\IC_{G,d}$ for $\IC(A\Qlb, d)$ for $d>0$. Let also $\IC_{G,0}=\Aut_g$ for $d=0$. From Proposition~\ref{Pp_A1} we get the following.

\begin{Lm} Let $d\ge 0$. The $*$-restriction of $\IC_{G,d}$ to $Shatz^r_{\tilde G}$ vanishes unless $r-d$ is even. The $*$-restriction of $\IC_{G,d}$ to $Shatz^{d+2k}_{\tilde G}$ identifies with $(\Qlb\boxtimes W)[\dim Shatz^d -2k]$. \QED
\end{Lm}
  
   For $d\ge 0$ let $Shatz^d_H\subset\Bun_H$ be the locally closed substack classifying $H$-torsors isomorphic to the push-forward of $\cO(1)$ by $\frac{d}{2}\alpha: \Gm\to H$, where $\alpha$ is the positive coroot of $H$. So, the corresponding orthogonal vector bundle is $\cO(d)\oplus\cO\oplus\cO(-d)$. For $a\in\ZZ/2\ZZ$ this gives the Shatz straification of $\Bun_H^a$ by $Shatz^d_H$ with $a=d\!\mod \! 2$. Write $\IC_{H,d}$ for the $\IC$-sheaf of $Shatz^d_H$ on $\Bun_H$.
Write $\ov{Shatz}^d_H$ for the closure of $Shatz^d_H$ in $\Bun_H$.   
Note that $Shatz^0_H$ (resp., $Shatz^1_H$) is the open Shatz stratum in $\Bun_H^0$ (resp., in $\Bun_H^1$).
   
Let $x\in X$. Denote by $_x\H_H: \D(\Bun_H)\to\D(\Bun_H)$ the Hecke functor corresponding to the standard representation of $\check{H}\,\iso\, \SL_2$.
The surjective homomorphism $\GL_2\to H$ gives a morphism of extension of scalars $\Bun_2\to\Bun_H$. Taking the inverse images of these sheaves on $\Bun_2$ and applying Hecke functors for $\GL_2$, one easily proves the following (the details are left to a reader). 
   
\begin{Pp} 1) For any $d\ge 0$ one has $\IC_{H, d}\,\iso\, \Qlb[\dim Shatz^d_H]$ over $\ov{Shatz}^d_H$.\\
2) One has canonically 
$$
_x\H_H(\IC_{H, 1})\,\iso\, \IC_{H,0}[1]+\IC_{H,0}[-1]\;\;\; \mbox{and}\;\;\;\; _x\H_H(\IC_{H,0})\,\iso\, \IC_{H,1}[1]+\IC_{H,1}[-1]
$$ 
3) Let $d\ge 2$. Then one has canonically  
$$
_x\H_H(\IC_{H, d})\,\iso\, \IC_{H, d-1}+\IC_{H, d+1}
$$
\end{Pp}
 
\begin{Cor} Assume $g=0$ and $n=1$. The functor $QL$ given by $QL(\IC_{H,d})=\IC_{G,d}$ for $d\ge 0$ extends uniquely to an equivalence of the category of pure complexes of weight zero on $\Bun_H$ with the full subcategory of $\D_-(\Bunt_G)$ of pure complexes of weight zero. This functor is compatible with $\ZZ/2\ZZ$-gradings and commutes with the action of $\Rep(\SL_2)$ on both sides.
\end{Cor}

The example of $g=0$ shows that it is natural to expect that $QL$ is exact for the perverse t-structures for any $g$ and $n$.

\end{document}